\documentclass{amsart}
\usepackage{amsfonts,amssymb,amsmath,formel,graphicx, bm, amscd}
\usepackage[all]{xy}

\usepackage[backref=page]{hyperref}
\hypersetup{colorlinks=true,citecolor=blue,linkcolor=blue,urlcolor=blue,
pdfstartview=FitH, pdfauthor=Valentin Blomer and J\"org Br\"udern and Per Salberger,pdftitle=On a certain senary cubic form}

\def\cal{\mathcal}
\def\ver{{\rm ver\,}}
\def\PP{\Bbb{P}}

\newtheorem{thm}{Theorem}

\newtheorem{lem}{Lemma}
\numberwithin{equation}{section}

\begin{document}

\author{Valentin Blomer}
\author{J\"org Br\"udern}
\author{Per Salberger}
\address{Mathematisches Institut, Bunsenstr. 3-5, 37073 G\"ottingen, Germany} \email{blomer@uni-math.gwdg.de}\email{bruedern@uni-math.gwdg.de}
\address{Mathematical Sciences, Chalmers University of Technology, SE-412 96 G\"oteborg, Sweden} \email{salberg@chalmers.se}

\title{On a certain senary cubic form}

\thanks{The first author was supported  by an ERC Strating grant and a Volkswagen Lichtenberg Fellowship. }


\begin{abstract}  A strong form of the Manin-Peyre conjecture with a power saving error term is proved for a certain cubic fourfold. 
 \end{abstract}


\maketitle

%\tableofcontents

\section{Introduction}
\subsection{Principal results}

The main goal of this paper is the verification of predictions due to Manin and Peyre concerning the distribution of rational points on the cubic   fourfold in ${\Bbb P}^5$ defined by
the equation
\be{1}
 W : \quad  x_1 y_2 y_3 + x_2 y_1 y_3 + x_3 y_1 y_2 =0.
\ee
We shall derive an asymptotic formula not only in line with the aforementioned expectation, but of strength sufficient to obtain an analytic continuation for the associated height zeta function beyond the region of absolute convergence. Along the way, we construct a crepant resolution of its singularities and determine the universal torsor, thus providing a comprehensive picture of the arithmetic and algebraic properties of the   fourfold defined by \eqref{1}.  The properties of this fourfold are sufficiently distinct from those among the small stock of cubic fourfolds for which the Manin-Peyre conjectures are already known, to require treatment by a new analytic toolbox. \\

More notation is required for precise statements of our results as well as comments on the methods involved.  If a   point in  ${\Bbb
  P}^5$ is represented by $(x_1,x_2,x_3,y_1,y_2,y_3)\in{\Bbb Z}^6$
with coprime coordinates, then $H(\textbf{x}, \textbf{y}) = \max (|x_j|,|y_j|)^3$ is a natural anticanonical 
height function on $W(\Bbb{Q})$.  The coordinate 3-planes where at least one of $y_1, y_2, y_3$ vanishes are accumulating subsets of $W$ with exceptionally many points. %There are
%exceptionally many points on the intersection with the union of
%coordinate hyperplanes described by
%\be{1a}
%$y_1y_2y_3 =0$. %\ee
For instance, choosing $y_1 = y_2 = 0$ and integral $x_1, x_2, x_3, y_3$ subject only to a coprimality condition, we find  about  $P^{4/3}$  rational points $(\textbf{x}, \textbf{y}) \in \Bbb{P}^5$   satisfying the equations \eqref{1} and $y_1y_2y_3 = 0$ as well as the height condition
\be{2}
  H(\textbf{x}, \textbf{y}) \leq P. 
\ee
%With a little more care  it is easy to see that there are
%\begin{displaymath}
 % \frac{48}{\zeta(4)} P^{4/3} + O(P)
%\end{displaymath}
%rational points   of height not exceeding $P$ satisfying the equations \eqref{1} and $y_1y_2y_3 = 0$. 
%To see
 % this, we count solutions of the system of equations \rf{1}, \rf{1a}
%in rational integers $x_j, y_j$ that lie within the box 
%\be{2}
%|x_j|\le P, \, |y_j|\le P \quad (1\le j\le 3);
%\ee
%here $P$ is a large real parameter. Choosing $y_1 = y_2 = 0$, there are no constraints on the remaining 4 variables, and we find immediately $\asymp P^4$ solutions. With a little more care one can show that there are
%First consider solutions with
 %$y_1=0$. Then, by
%\rf{1}, one has  $x_1 y_2 y_3=0$ so that at least one of $x_1, y_2,
%y_3$ must also vanish. However, when
%$y_1$ and only one of $x_1, y_2, y_3$ are zero, then \rf{1} already
%holds. This determines all solutions of \rf{1} with $y_1=0$. 
%By permuting the indices, one finds that \rf{1} and \rf{1a} hold
%simultaneously if and only if  two of $y_1, y_2, y_3$ are zero,
%or 
% $x_j = y_j =0$ holds for some
%$j$.  There are no
%constraints on the remaining four variables in all cases, whence it
%follows that there are
%$$
%96 P^4 + O(P^3)
%$$
%integral solutions to \rf{1}, \rf{1a} and  \rf{2}. %With a little
%more care, one even finds that the exact number in \rf{3} is a
%polynomial of degree $4$ in $P$, at least when $P$ is an integer.
%\medskip
  Our principal result concerns the density of the
  rational points on the Zariski open subset $W^{\circ}$ defined by \eqref{1} and 
\be{nonzero}
  y_1y_2y_3 \not = 0.
\ee
Let
 $N(P )$ denote the number of rational points on $W^{\circ}$ satisfying 
 % of sextuples $(x_1,
%x_2, x_3, y_1, y_2, y_3)\in\ZZ^6$ 
satisfying  \rf{2}. %, and let $V^*(P)$ be the number of those solutions
%counted by $V(P)$ that satisfy the coprimality condition 
%$(x_1;x_2;x_3;y_1;y_2;y_3)=1$. 
\begin{thm}\label{thm1}
 There is a constant $\del >0$ and a real
  polynomial $Q$ of degree $4$ such
  that
\begin{displaymath}
N(P )= P \, Q(\log P) + O(P^{1-\del}).
\end{displaymath}
The leading coefficient of $Q$ equals
\be{pey*}
\frac{1}{324}(\pi^2 + 24\log 2 - 3 )
\prod_p \Big(1-\f1p\Big)^5
\Big(1+\f{5}{p}+\f{6}{p^2}+\f{5}{p^3}+\f{1}{p^4}\Big).
%\lim_{P\to\infty} \f{N(P)}{P^3(\log P)^4} = 
%\frac{1}{2}(\pi^2 + 24\log 2 - 3)  \prod_p \Big(1-\f{9}{p^2}+\f{16}{p^3}-\f{9}{p^4}+\f{1}{p^6}\Big) = 0.57\ldots.
\ee
\end{thm}
It is possible to provide a numerical value for $\delta$,   but we have made no effort to optimize it. A careful estimation would produce a reasonable size for $\delta$. \\

%The geometry of smooth cubic fourfolds has been investigated extensively, for instance in \cite{Has}. 
The cubic fourfold $W$ defined by \rf{1} is singular. We note, however, that any
singular solution of \rf{1} %must have $y_j=0$ for at least two values
%of $j$, and so 
is contained in   $W \setminus W^{\circ}$. We also observe  that $W^{\circ}$ is completely covered by rational planes 
% One sees directly that $W$ contains rational planes in abundance, given by 
$x_j = a_jy_j$ ($1 \leq j \leq 3$), which is in sharp contrast to the case of smooth cubic fourfolds.  % whenever $a_1 + a_2 + a_3 = 0$. It can be shown, however, that the number of points on these projective planes counted by $N(P )$ is only $O(P)$, hence the majority of rational points lies outside these linear subvarieties.\\

For a better understanding of the geometry of  a singular variety it is necessary to construct a desingularization. In the case of the variety $W$ defined in \eqref{1}, we are in the convenient situation that there exists a \emph{crepant} resolution. Recall that a resolution $f: X\rightarrow W$ of a normal variety with invertible canoncial sheaf $\omega_W$  is said to be crepant if  $f^{\ast}\omega_W = \omega_X$ (see \cite{Re1}).  Specifically, let
  $X \subset \Bbb{P}^{5} \times \Bbb{P}^{2} \times \Bbb{P}^{2}$ be the tri-projective variety with tri-homogeneous coordinates $(x_1, x_2, x_3, y_1, y_2, y_3; Y_1, Y_2, Y_3; Z_1, Z_2, Z_3)$ defined by the equations
\begin{displaymath}
\begin{split}
  & x_1 Z_1 + x_2 Z_2 + x_3Z_3 = 0,\\
  &    y_iY_j - y_jY_i = 0, \quad 1 \leq i < j \leq 3,\\
  & Y_1Z_1= Y_2Z_2 = Y_3Z_3.
\end{split}  
 \end{displaymath}
 In Chapter \ref{secres} we shall prove a more general version of the following result. 
 \begin{thm}\label{thmres} The projection $\Bbb{P}^5 \times \Bbb{P}^2 \times \Bbb{P}^2 \rightarrow \Bbb{P}^5$ restricts to a   crepant resolution $f : X \rightarrow W$. 
 \end{thm}

%[Manin conjecture: various conjectures, known cases, resolution, Theorem 1 agrees] 
Manin has put forward a fundamental conjecture relating the geometry of a projective variety to the arithmetic of its rational points  (\cite{FMT},  \cite{BM}). Originally this  
conjecture was formulated  for smooth Fano varieties. The number of $\log$-powers in an asymptotic formula for the density of rational points is one off the rank of the Picard group, and Peyre  \cite{P} (see also \cite{Sa1})  
has suggested a formula for the leading constant. This has been  generalized to large classes of singular Fano varieties by Batyrev and Tschinkel in   \cite{BT}. Since the resolution in Theorem \ref{thmres} is crepant, $X$ is an ``almost Fano" variety in the sense of \cite[Def.\ 3.1]{P2} (cf.\ Lemma \ref{almostfano} below). In particular, Peyre's ``Formule empirique" \cite[5.1]{P2} is expected to predict the asymptotic behaviour of the counting function $N(P )$.   We will discuss Peyre's   formula in detail in Chapter \ref{maninpeyre}, and we show  in Theorem \ref{thmpey} that it agrees with our Theorem \ref{thm1}.  In particular, the 
 Tamagawa constant may be interpreted as an adelic volume of the  universal torsor  over the crepant resolution which we describe  in Chapter \ref{descent}. 

% For singular Fano varieties with canoncial   singularites there is still a corresponding 
%conjecture by Batyrev-Tschinkel \cite{BT}, but the Picard numbers are
%not the ones given by the singular variety itself, but rather by a suitable desingularisation.

Significant progress on the Manin-Peyre conjecture has been made for %del Pezzo and Ch\^atelet
surfaces. Using a variety of methods from analytic number theory, it has been verified in a number of cases by Browning, de la Br\`eteche, Derenthal, Peyre and others. The important papers \cite{BaBr, dlB1,  BBD, BBP} and the  references in   \cite{Br} will guide the reader into the realm of the extensive research literature. There are few definitive results on Fano threefolds, and the remarkable paper \cite{dlB2} on the Segre cubic illustrates the additional difficulties that may appear. For cubic fourfolds, we do not know of any asymptotic formulas except for toric varieties related to the hypersurface $x_1x_2x_3 = y_1y_2y_3$.  %The cubic fourfold \eqref{1} under consideration is very different from all of these cases, and we need to introduce a new toolbox of techniques to establish the Manin-Peyre conjecture in this case.  
 %Very few   explicit results of the type of Theorem \ref{thm1} are known for higher dimensional hyper surfaces, save for toric varieties such as the   cubic fourfold $x_1x_2x_3 = y_1y_2y_3$. 
%and Theorem \ref{thm1} seems to be the first instance where the density of rational points on a cubic fourfold is considered. In fact, the authors are not aware of any result in this direction where the ratio of  the number of (affine) variables in the defining equation and the
%degree is  $2$  (save for
%the classical theory of quarternary quadratic forms). 
The investigation of higher-dimensional varieties is  an interesting testing ground for more general versions of the Manin-Peyre conjecture, and we hope that the present rather complex example will initiate further research. \\ %For instance, in the case of singular del Pezzo surfaces it is well known that there is a  unique
%minimal canonical resolution.  In higher dimension one may not always expect that there is just one
%natural minimal resolution, and hence even the prediction of the number of log-factors and the Peyre-constant are not obvious any more. 

 %With this resolution at hand, the %Batyrev-Tschinkel conjectures are of the same kind as  the original conjectures of Manin and Peyre for non-singular varieties. 
 
%We will verify in detail in Section \ref{maninpeyre} that the result of Theorem \ref{thm1} agrees with Peyre's ``empiric formula".\\ 

Although often quoted as an asymptotic relation, Manin's conjecture should  be considered as a statement concerning the analytic continuation of the corresponding height zeta function.  In our case, the height zeta function attached to the non-trivial part of the cubic defined in \eqref{1} is given by
\begin{displaymath}
  Z(s) = \sum_{ (\textbf{x}, \textbf{y}) \in W^{\circ}} H(\textbf{x}, \textbf{y})^{-s}, \quad \Re s > 1. 
\end{displaymath}
A routine partial summation coupled with Theorem \ref{thm1} yields an analytic continuation for $Z(s)$.  
\begin{thm}\label{kor} Let $\delta$ be as in Theorem \ref{thm1}. The height zeta function $Z(s)$ has analytic continuation to a right half plane ${\rm Re \,} s > 1-\delta$ except for a pole at $s=1$ of order 5. In the region ${\rm Re\,} s > 1-\delta$, $| s-1| > 1/10$ one has the growth estimate $Z(s) \ll |s|$. 
\end{thm}

Finally we remark that the variety $W$ carries  an   algebraic structure. The open subset $W^{\circ}$ is  an abelian group if the product of the two points  $(x_1,x_2,x_3,y_1,y_2,y_3)$ and $ (x'_1,x'_2,x'_3,y'_1,y'_2,y'_3)$ is defined by 
\begin{displaymath}
  (x_1y_1' + x_1'y_1, x_2y_2' + x_2'y_2, x_3y_3' + x_3'y_3, y_1y_2', y_2y_2', y_3y_3').
\end{displaymath}
Somewhat more conceptually,  let $H$ be the algebraic  group of all $2 \times 2$  matrices of the form $\left(\begin{matrix}b & a\\ 0& b\end{matrix}\right)$ with $b$ invertible. Let $\Psi : H^3 \rightarrow \Bbb{G}_a$ be the homomorphism 
\begin{equation}\label{hom}
\left(\left(\begin{matrix}b_1 & a_1\\ 0& b_1\end{matrix}\right), \left(\begin{matrix}b_2 & a_2\\ 0& b_2\end{matrix}\right),   \left(\begin{matrix}b_3 & a_3\\ 0& b_3\end{matrix}\right) \right) \mapsto \frac{a_1}{b_1} + \frac{a_2}{b_2} + \frac{a_3}{b_3},
\end{equation}
and define $G = {\rm ker}(\Psi)/\Bbb{G}_m$ where $\Bbb{G}_m$ is embedded   into $H^3$ via 
$
 b \mapsto \left(\left(\begin{matrix}b & 0\\ 0& b\end{matrix}\right),  \left(\begin{matrix}b &0 \\0 & b\end{matrix}\right),  \left(\begin{matrix}b & 0\\ 0 & b\end{matrix}\right) \right). 
 $
Then 
\begin{equation}\label{group}
   W^{\circ} \cong G \cong \Bbb{G}_a \times \Bbb{G}_a \times \Bbb{G}_m \times \Bbb{G}_m
 \end{equation}  
    as groups, and  there is a natural open immersion $G \rightarrow W$ and a natural $G$-action on $W$. Hence we identify $W$ with the equivariant compactification of $G$, the product of two additive and two multiplicative groups.  In various cases the group structure can be employed to prove an asymptotic formula for the number of points of bounded height using adelic Fourier analysis. This has been carried out for instance for toric varieties in \cite{BT1} (with a different proof in \cite{Sa1} and \cite{dlB4}) and equivariant compactifications of additive and certain other groups (e.g.\ \cite{CLT,  STBT}) including  a non-commutative example   in \cite{TT}. None of these cases cover the situation of a mixed additive and multiplicative abelian group, and the present case seems to be the first example in the literature for a group of the type \eqref{group}. 
     It is possible, however, that an extension of the Fourier analytic techniques could also produce a result similar to the one announced in Theorem \ref{thm1}, and 
     it would be interesting to compare the two approaches.  
%Cases where the ratio between the number of variables and the
%degree is nearly $2$  are rare, and Theorem \ref{thm1} seems to be the first time in the literature that a cubic fourfold has been analyzed.  The authors are not aware of any such
%result where this ratio is actually $2$ as in Theorem \ref{thm1}, save for
%the classical theory of quarternary quadratic forms (cf.\ \cite[Theorem 7]{HB}). \\

%We believe that  Theorems \ref{thm1} -- \ref{kor} together with the computation of the universal torsor in Theorem \ref{thm7}  and the verification of Peyre's empiric formula in Theorem \ref{thmpey} and  provide 

\subsection{The methods} 
The proof of Theorem \ref{thm1} draws from a wide range of methods. The main
argument that we now describe uses three very different tools,
namely elementary lattice point considerations, analytic counting by
multiple Mellin integrals, and an Euler product identity for certain multiple
Dirichlet series.

The initial step is not new. Rather than counting integral points on
the cubic \rf{1} directly, we pass in   Section  \ref{torsor}) to a descent variety, a frequently used technique in this context.  %Torsors
%have been used in this context by Salberger \cite{S} and Peyre
%\cite{P}, among others, and also by Heath-Brown \cite{HB}, the latter
%in a style similar to our use of it. 
%This turns
%out to be straightforward . 
After a succession of divisibility
considerations, one ends up with a bilinear equation. The underlying
lattice structure then provides a complete parametrization of the
cubic \rf{1}, see Lemma \ref{lem2}. These elementary  arguments can be interpreted 
in terms of equations for the universal torsor of the variety; we carry this out explicitly in Section \ref{universaltorsor} (cf.\ the companion Lemma \ref{lemtorsor}) in order to provide further insight into the genesis of the leading term in the asymptotic formula \eqref{pey*}.

The parametrization now in hand, the count for $N(P)$ has a new
interpretation as a 10-dimensional lattice point problem with a
strangely shaped boundary. An enveloping argument along the lines of \cite{BB} would   produce the correct order of magnitude 
%that may be described
%as elementary but not necessarily as simple, then provides the
%preliminary estimates
$
  N(P ) \asymp P(\log P)^4. 
$ 
Alternatively, one may approach the lattice point problem as a heavily
convoluted divisor sum, and tackle it by Mellin transform techniques.
Due to the
complicated boundary conditions, we will require
multidimensional Mellin inversion formulae. 
This part of the argument seems new in this context and 
should have a wider range of applications to the Manin-Peyre conjecture and
cognate problems in diophantine analysis;  we shall mention some in
Section \ref{appl} below. It would take us too far afield at
this point to comment on the finer structures of our techniques. 
We content ourselves here with the remark  that the
method uses a delicate regularization process of a priori divergent integrals and eventually 
produces
an asymptotic formula for a counting function that mimics $N(P )$, but
has a smooth weight attached to the variables, see Section \ref{sec81}.  %The main term in this
%formula is of the shape
%$P^3Q_0(\log P)$ with {\it some} polynomial $Q_0$. In principal, it
%would be possible to control the degree of $Q_0$, but only by  tedious
%computation.  Instead, we first remove the weights and obtain an
%asymptotic formula for $N(P)$, as in \rf{as}. Then \rf{A2} shows that
%the degree of $Q$ must be 4. For a more detailed description of this
%part as well the relation of our work with an important result
%of de la Breteche \cite{B-comp}, we refer the reader to section 4.
The proof of Theorem \ref{thm1} is then completed by removing the weights and computing the main term   through the
analytic machinery as a certain residue. This approach is inspired by work of de la Bret\`eche \cite{dlB}, but our situation is rather more involved. Some of the additional complications are discussed  in the final remark of Chapter \ref{weights7}.    If we were only interested in the leading term of the asym\-ptotic formula, we could completely dispense with the rather long and technical Chapter \ref{remove} (as well as Chapter \ref{upper} and   Lemma \ref{lem12}) and use a standard Tauberian argument instead. However, the full asymptotic formula and the analytic continuation of the height zeta function in Theorem \ref{kor} require the more complex argument. 

The Euler product formula
\rf{pey*} arises from an identity for a family of multiple Dirichlet series 
that may be of some independent interest. It can be used 
successfully
in many cases where the   analytic method is
applicable. Therefore we highlight this auxiliary result as Theorem \ref{thm3}
below, and reserve Chapter \ref{Diri} for its discussion and demonstration.  \\

The proof of Theorem \ref{thmres} proceeds in two steps. We first consider the closure $W' \subset  \Bbb{P}^5 \times \Bbb{P}^2$   of the graph of the projection $W\setminus \Pi \rightarrow \Bbb{P}^2$  from the plane $\Pi$ given by $y_1=y_2=y_3=0$ and show that $\text{pr}_1 : \Bbb{P}^5 \times \Bbb{P}^2 \rightarrow \Bbb{P}^5$ restricts to a crepant birational morphism $ g: W'\rightarrow W$. Then, after a base extension of $\text{pr}_2 : W' \rightarrow \Bbb{P}^2$, we obtain a $\Bbb{P}^2$-bundle $\lambda : X \rightarrow B$  over a non-singular del Pezzo surface $B$ of degree 6, where $X$ is crepant over $W$. Since this desingularisation is $G$-equivariant,  we may  compute Peyre's alpha invariant in Lemma \ref{alpha} by means of a result from \cite{TT}.  

%Our method is fairly general and   in particular does not use the fact that $W$ is an equivariant compactification of a group.
%\textbf{say a few lines something about the proof of Theorem 2 and the computation of Peyre's constant if appropriate}\\

%The proof of Theorem 2 follows the same pattern, but only up to the
%elementary part. It would be interesting to make the analysis work for
%the bihomogenous case as well, and we hope to return to this problem
%in due course.

\subsection{Further applications}\label{appl}
 When $y_1 y_2 y_3 \neq 0$, one may rewrite
\rf{1} as
$$
\f{x_1}{y_1} + \f{x_2}{y_2} + \f{x_3}{y_3} =0,
$$
and there is then a natural generalization to more than three
summands. When $n>3$, %and denominators are cleared in 
the solutions of the equation
\be{C1}
\f{x_1}{y_1} + \f{x_2}{y_2} + \ldots + \f{x_n}{y_n} =0
\ee
are the zeros of  a form of degree $n$ in $2n$ variables. A further
development of our techniques yields 
results for this form that are comparable to
Theorem \ref{thm1}. In fact, parts of the arguments are carried out for arbitrary $n$. 

From a more arithmetic point of view, one may also count fractional zero sums of bounded height,
that is, solutions of \rf{C1} with $(x_j; y_j)=1$ and $|x_j|, |y_j|\le
P$. The extra coprimality conditions decrease the power of the
logarithm that appears in the asymptotic formula. This phenomenon has
been observed for other forms as well, see Fouvry \cite{F} for a
discussion in the case $x_0^3 = x_1 x_2 x_3$.\\ %We hope to return to this point at some later occasion. 

In a different direction, we note that the cubic form on the left of \rf{1} is a linear form in ${\bf
  x}$, and a quadratic form in ${\bf y}$. Therefore 
\rf{1} also defines a singular bi-projective cubic threefold  $\tilde{W} \subset \PP^2 \times \PP^2$, and again  one may count its rational points with respect to the
anticanonical height and compare the result with the predictions by Manin and Peyre. %Although our normalisation is somewhat
%different from what is common practice in the area, 
This amounts to
analyzing the number $\tilde{N}(P )$ of non-trivial solutions to \rf{1} that satisfy
the conditions $(x_1; x_2; x_3) = (y_1; y_2; y_3) =1$, $x_1x_2x_3y_1y_2y_3 \not= 0$ and the size
constraints
\begin{equation}\label{newnormalization}
1\le |x_i^2 y_j| \le P \qquad (1\le i, j\le 3).
\end{equation}
Note that $\tilde{W}$ cannot be written as a compactification of a group in a natural way. 
In \cite{BB} we were able to determine the order of magnitude of $\tilde{N}(P )$: with the normalization \eqref{newnormalization} one has
$$
P(\log P)^4 \ll \tilde{N}(P ) \ll P(\log P)^4.
$$
We take the opportunity to relate these estimates to the standard predictions  and show that the order of magnitude agrees with the expected one. The singular locus of $\tilde{W}$ is given by the three points $x_i = x_{i+1} = y_i = y_{i+1} = 0$, $i \in \{1, 2, 3\}$, with indices understood modulo 3. Let $\tilde{X} \subset \Bbb{P}^2 \times \Bbb{P}^2 \times \Bbb{P}^2$ be the tri-projective variety with tri-homogeneous coordinates $(\textbf{x}; \textbf{y}; \textbf{z}) = (x_1, x_2, x_3; y_1, y_2, y_3; z_1, z_2, z_3)$ defined by
\begin{equation}\label{tri2}
  x_1z_1 + x_2z_2 + x_3z_3 = 0,
\end{equation}
\begin{equation}\label{tri3}
  y_1z_1=y_2z_2 = y_3z_3.
\end{equation}
Similarly as in Theorem \ref{thmres} we prove
\begin{thm}\label{thmthree}
The restriction to $\tilde{X}$ of projection $\Bbb{P}^2 \times \Bbb{P}^2 \times \Bbb{P}^2 \rightarrow \Bbb{P}^2 \times \Bbb{P}^2$ onto the first two factors is a crepant resolution of $\tilde{W}$, and one has $\text{\rm rk } \text{\rm Pic}(\tilde{X}) = 5$. 
\end{thm}

%Similarly as in Theorem \ref{thmres} we construct a crepant resolution of $\tilde{W}$ and compute the Tamagawa numbers. On the basis of these calculations we make the following

%\textbf{Conjecture.} One has
%\begin{displaymath}
%  \tilde{V}^{\ast}(P ) = (1+o(1)) ??? P^3(\log P^4), \quad P \rightarrow \infty. 
%\end{displaymath}

\subsection{Leitfaden}  We start in Chapter \ref{Diri} with the graph theoretic proof of Theorem \ref{thm3} which will be used to compute explicitly the Euler product in \eqref{pey*}. Chapter \ref{secres} features the proof of Theorem \ref{thmres} (along with Theorem \ref{thmthree}). We also take the opportunity to compute Peyre's alpha invariant  at the end of this chapter. In Chapter \ref{descent} we pass to the universal torsor (see in particular Theorem \ref{thm7}). Then we are prepared to discuss the ``empirical  formula"  suggested by Peyre in Chapter \ref{maninpeyre}. The main ingredient for this is Theorem \ref{thmpey}. The remaining chapters are devoted to the analytic proof of Theorem \ref{thm1}. Chapter \ref{upper} provides some preliminary upper  bounds for $N(P )$ and related quantities that will be needed later. Chapter \ref{weights} is of technical nature and introduces certain smooth weight functions along with properties of their (multi-dimensional) Mellin transforms; the proofs can safely be skipped at a first reading.  Chapter \ref{analytic} is the heart of the proof. By Mellin transform and contour shifts the asymptotic formula is reduced to the calculation of certain residues. Chapter \ref{remove} is again rather technical, and its only purpose is to remove the smooth weights in order to get an asymptotic formula for a count in a box. Finally, the leading coefficient is computed in Chapter \ref{peyre}, completing the proof of Theorem \ref{thm1}.

\section{Some combinatorial identities}\label{Diri}
\subsection{Multiple Dirichlet series with coprimality constraints} %The
%main goal of this chapter is to establish a certain identity for a
%multiple Dirichlet series, Theorem \ref{thm3} below, that plays a pivotal role
%in the computation of the constant \eqref{pey*} in Theorem 1. 

For a natural number $r$ let $G = (V, E)$ be any graph on the set of vertices $V = \{1, \ldots, r\}$. 
Then, whenever $s_1, \ldots, s_r$ are complex numbers with $\Re s_j
>1$, the series
\be{A1}
D_G(s_1, \ldots, s_r) = 
\multsum{n_1,\ldots,n_r = 1}{(n_k; n_l) =1 \text{ for } (k,l)\in E}
 n_1^{-s_1}n_2^{-s_2} \cdots n_r^{-s_r}
\ee
is absolutely convergent.    %For
%${\bf s} \in\CC^r$ we put
%\be{A2}
%s(U)= \sum_{j\in\ver (U)} s_j.
%\ee
Note that when $r=1$, then $E$ is necessarily empty, and \rf{A1}
reduces to the definition of Riemann's zeta function
$\zeta(s_1)$. Likewise, for any $r\in\NN$, one finds from \rf{A1} that
$$
D_{(V, \emptyset)}(s_1, \ldots, s_r) = \zeta (s_1) \cdots \zeta(s_r),
$$
so that there is an analytic continuation to $\CC^r$ except for
singularities at $s_j=1$. The following result describes the
situation for any graph $G=(V,E)$. 
For a subset $U\subset E$  we define its vertex set $\ver U \subset V$
 as the
set of all vertices that are adjacent to at least one edge in $U$. 

\begin{thm}\label{thm3} Let ${\bf s}\in\CC^r$ with  ${\rm Re\,} s_j >1$ for
  $1\le j\le r$. Then 
\be{A3}
\zeta(s_1)^{-1} \cdots \zeta(s_r)^{-1} D_G({\bf s}) = \prod_p
\sum_{U\subset E}(-1)^{|U|} \prod_{j \in \ver U } p^{-s_j}. %p^{-s(U)}.
\ee
The Euler product on the right hand side of \rf{A3} converges absolutely in the
region ${\rm Re \,} s_j > \f{1}{2}$ $(1\le j\le r)$, and constitutes an
analytic continuation of the function 
$$\Xi_G({\bf s}) =
\zeta(s_1)^{-1} \cdots \zeta(s_r)^{-1} D_G({\bf s})$$ to this set.
\end{thm}

The simplest example not yet considered is when $r=2$,
$E=\{(1,2)\}$. Here the Euler factors in \rf{A3} are $1-p^{-s_1
  -s_2}$, and the principal conclusion of Theorem \ref{thm3} reduces to the
familiar identity
$$
\multsum{n_1, n_2=1}{(n_1; n_2) =1}^\inf n_1^{-s_1} n_2^{-s_2} =
\f{\zeta(s_1) \zeta(s_2)}{\zeta(s_1+s_2)}.
$$

We also note that it will suffice to establish \rf{A3}, because for
any $\emptyset \not= U\subset E$  the vertex set has at least
two elements, and consequently, one finds that whenever $\Re s_j
>\f{1}{2}$, then 
$$
\sum_{j \in \ver U} \Re s_j > 1.
$$
Hence the product in \rf{A3} indeed
converges absolutely in the indicated region.\\

Theorem \ref{thm3} can be generalized in various ways. For instance, if $\alpha_j$ ($1\leq j \leq r$) are arbitrary completely multiplicative functions with corresponding Dirichlet series $L_j(s) = \sum_{n=1}^{\infty} \alpha_j(n) n^{-s}$, 
then the identity
\be{general}
L_1(s_1)^{-1} \cdots L_r(s_r)^{-1} D_G({\bf s}) = \prod_p
\sum_{U\subset E}(-1)^{|U|} \prod_{j \in \ver U} \alpha_j( p) p^{-s_j} %p^{-s(U)}.
\ee
holds in the region of absolute convergence. The proof of \eqref{general} is the same. Martin \cite[Proposition A.4]{Mar} has proved a related result for more general multiplicative functions, but only for a complete graph $G$. \\

In applications, it is desirable to compute the number $\Xi_G
(1,\ldots, 1)$ explicitly. By \rf{A3}, one has
$$
\Xi_G (1,\ldots, 1) = \prod_p \sum_{U\subset E} (-1)^{|U|}
p^{-|\ver U|} = \prod_p \sum_{k=0}^r b_k p^{-k}
$$
where
\be{A5}
b_k = \multsum{U\subset E}{\ver U =k} (-1)^{|U|}.
\ee
Note that one has $b_0 = 1, b_1=0, b_2 = -| E|$ and
\be{A6}
\sum_{k=0}^r b_k = \sum_{U\subset E}(-1)^{|U|} =0.
\ee

In any concrete example, the numbers $b_k$ can be computed via
\rf{A5}. In this paper, the 
case of interest is when $r=6$ and $E$ is the set of pairs $\{(1,2), (1,3),
(2,3), (4,5), (5,6),$ $(4,6), (1,4), (2,5), (3,6)\}$.
\begin{equation}\label{graph}
\xymatrix@R=5pt@C=20pt{
1 \ar @{-} [rrrr]
  \ar @{-} [dd]
  \ar @{-} [rd]  &                 & &   & 4 \ar @{-} [dd] \ar @{-} [ld] \\
                 & 3 \ar @{-} [rr] & & 6 &   \\
2 \ar @{-} [rrrr]
  \ar @{-} [ru]  &                 & &   & 5 \ar @{-} [lu]
}
\end{equation}
Here one finds that $b_3 = 16$, $b_4= -9$, $b_5=0$ and $b_6=1$. This gives
\be{A7}
\Xi_{ G} (1,1,1,1,1,1)= \prod_p\Big(1-\f{9}{p^2} + \f{16}{p^3}
-\f{9}{p^4} + \f{1}{p^6}\Big)=
\prod_p\Big(1-\f{1}{p}\Big)^4\Big(1+\f{4}{p}+\f{1}{p^2}\Big).
\ee
We briefly indicate how one may compute $b_3$. By \rf{A5}, one first
has to determine all $U\subset E$ with $|\ver U|=3$. Such a set $U$
is either a complete subgraph on three vertices, or it consists of two
edges with one vertex in common. There are exactly two complete
subgraphs on three vertices in ${G}$, so this class contributes
$-2$ to $b_3$. In order to count the other class of sets $U$, one
first chooses one of the six vertices to determine the vertex that the
two edges should have in common. Since each vertex has three adjacent
edges, one can then make three choices of two edges. Hence, there are 18
pairs of edges with a vertex in common, and \rf{A5} yields $b_3 = 18-2
= 16$, as required. It is equally straightforward but more elaborate
to compute $b_4$ and $b_5$. The coefficient $b_6$ can then be determined from \eqref{A6}. We leave the details to the reader.

\subsection{Graphs and power series} In this section we reduce the
proof of Theorem \ref{thm3} to an identity for a power series associated with
the graph $G$. To define this series, let $\delta: \NN_0\to \{0,1\}$
be defined by 
$\delta(0)=1$, $\delta(n) =0$  for all $ n\ge 1$, 
and put
\be{A8}
\Delta_G (n_1, \ldots, n_r) = \prod_{(k,l)\in E} \del(n_k n_l).
\ee
The power series
\be{A9}
T_G(x_1, \ldots, x_r) = \sum_{{\bf n}=0}^\inf \Delta_G ({\bf n}){\bf
  x}^{\bf n}
\ee
converges in the disk $|x_j|<1$ $(1\le j\le r)$. It turns out that
$T_G$ is a rational function that becomes a polynomial when multiplied with
\be{A10}
\Pi_r (x_1, \ldots, x_r) = (1-x_1) (1-x_2) \ldots (1-x_r).
\ee
For any $U\subset E$ let
\be{A11}
{\bf x}_U = \prod_{j\in\ver U} x_j.
\ee

\begin{lem}\label{prop4}  One has
$$
\Pi_r ({\bf x}) T_G({\bf x}) = \sum_{U\subset E} (-1)^{|U|} {\bf
  x}_U.
$$
\end{lem}

In the next section, we establish Lemma \ref{prop4} by induction on $r$, but
now is the time to deduce Theorem \ref{thm3}. 
The indicator function on the conditions $(n_k;
n_l)=1$ for $(k,l)\in E$ is a multiplicative function on $(n_1,
\ldots, n_r)$. Hence, the Dirichlet series \rf{A1} has an Euler
product in its region of absolute convergence. By \rf{A1}, \rf{A8} and
\rf{A9}, this takes the shape
$$
D_G({\bf s}) = \prod_p T_G(p^{-s_1}, \ldots, p^{-s_r}).
$$
Hence, by Lemma \ref{prop4}, \rf{A10} and the Euler product for Riemann's zeta
function,
$$
D_G({\bf s}) = \zeta(s_1) \zeta(s_2) \cdots \zeta(s_r) \prod_p
\sum_{U\subset E} (-1)^{|U|} \prod_{j\in\ver U} p^{-s_j}.
$$
This confirms Theorem \ref{thm3}.

\subsection{An inductive strategy} It remains to prove
Lemma \ref{prop4}. First consider the case where $E=\emptyset$ is the empty
set. Then, by \rf{A8}, one has $\Delta_{(V,\emptyset)}({\bf n}) =1$ for
all ${\bf n}\in\NN^r$, whence
$$
T_{(V,\emptyset)}({\bf x}) = \sum_{{\bf n}=0}^\inf {\bf x}^{\bf n} =
(1-x_1)^{-1}\cdots (1-x_r)^{-1},
$$
as is claimed in Lemma \ref{prop4}.

When $r=1$, then $E=\emptyset$ is the only possibility. When $r=2$,
and $E$ is non-empty, then $E=\{(1,2)\}$. In this case,
\begin{eqnarray*}
T_G(x_1, x_2) = 
\sum_{n_1, n_2=0}^\inf \del (n_1 n_2) x_1^{n_1} x_2^{n_2}  
=
\sum_{n_1=0}^\inf x_1^{n_1} + \sum_{n_2=0}^\inf x_2^{n_2} -1 = \f{1}{1-x_1} + \f{1}{1-x_2} -1,
\end{eqnarray*}
which is equivalent to the conclusion of Lemma \ref{prop4}. This settles
Lemma \ref{prop4} when $r=1$ or $2$. We may now suppose that $r\ge 3$, and
that Lemma \ref{prop4} is already established for smaller values of
$r$. Moreover, we have already dealt with the case where $E$ is
empty. In the opposite situation, there is at least one edge in $E$,
and by renumbering the vertices, we may suppose that $(1,2)\in E$. We
then consider the graph $G' = (V, E')$ with $E'= E\bs\{(1,2)\}$ and note that \rf{A8}
implies the equation
$$
\Delta_G ({\bf n}) = \del (n_1 n_2) \Delta_{G'} ({\bf n}).
$$
It will now be convenient to write ${\bf n}=(n_1, n_2, {\bf m})$ with
${\bf m}=(n_3, \ldots, n_r)$, and likewise, ${\bf x}= (x_1, x_2,{\bf
  y})$. Then, since $\delta(n_1 n_2)=1$ holds if and only if $n_1
n_2=0$, one finds from \rf{A9} that
\begin{eqnarray}
T_G({\bf x}) &=& \sum_{n_2=0}^\inf \sum_{{\bf m}=0}^\inf \Delta_{G'}
(0,n_2, {\bf m}) x_2^{n_2} {\bf y}^{\bf m} \nonumber\\
&&+
\sum_{n_1=0}^\inf \sum_{{\bf m}=0}^\inf \Delta_{G'}(n_1, 0,{\bf m})
x_1^{n_1}{\bf y}^{\bf m}
-
\sum_{{\bf m}=0}^\inf \Delta_{G'} (0,0,{\bf m}) {\bf y}^{\bf m}.
\label{A12}
\end{eqnarray}

Let $G_1= (V \setminus \{1\}, E_1)$ be the graph that is obtained from $G$ by removing the vertex 1 and all edges
$(1,l)$ adjacent to $1$.  Similarly, let $G_2 = (V \setminus \{2\}, E_2)$ be the graph that is obtained from $G$ by
removing the vertex 2 and all edges adjacent to 2.  Finally, let $H= (V \setminus\{1, 2\}, E_{1, 2})$ be the graph
that is the graph $G$ with all edges adjacent to $1$ or $2$
removed. Then  
\rf{A8} implies that 
%\begin{eqnarray*}
$\Delta_{G'}(0,n_2, {\bf m}) = \Delta_{G_1}(n_2, {\bf m})$,  
$\Delta_{G'}(n_1,0,{\bf m}) = \Delta_{G_2}(n_1, {\bf m})$, and  
$\Delta_{G'}(0,0,{\bf m})= \Delta_H({\bf m})$.
%\end{eqnarray*}
By \rf{A9} and \rf{A12}, we infer that
$$
T_G({\bf x}) = T_{G_1}(x_2, {\bf y}) + T_{G_2}(x_1, {\bf y}) -
T_H({\bf y}).
$$
We may now apply the induction hypothesis three times on the right
hand side.  For any subgraph $G^{\ast} = (V^{\ast}, E^{\ast})$ of $G$ we put
\be{A13}
S_{G^{\ast}}({\bf x}) = \sum_{U\subset E^{\ast}} (-1)^{|U|} {\bf x}_U,
\ee
and then we have by \rf{A10} that
$$
\Pi_r ({\bf x}) T_G({\bf x}) = (1-x_1) S_{G_1}({\bf x}) + (1-x_2)
S_{G_2}({\bf x}) -(1-x_1)(1-x_2) S_H({\bf y})
$$
where it is worth remarking that $S_{G_1}$ is a function of $x_2, {\bf
  y}$ only, and similarly for $S_{G_2}$. Now let $Q_j = S_{G_j} -
S_H$. Then the previous identity may be rewritten as
\be{A14}
\Pi_r ({\bf x}) T_G({\bf x}) = Q_1 + Q_2 + S_H - x_1 Q_1 - x_2 Q_2 -
x_1 x_2 S_H.
\ee
By \rf{A13}, the induction will be complete if the right hand side of
\rf{A14} can be shown to equal $S_G$. Before proceeding in this
direction, we first derive a useful formula for $Q_1$. In fact, a set
$U\subset E_1 \setminus E_{1, 2}$ is characterized by the
condition that $2 \in\ver U$. Hence, by \rf{A13},
\be{A15}
Q_1 = S_{G_1} - S_H = \multsum{U\subset E_1}{2\in\ver U}(-1)^{|U|}
{\bf x}_U.
\ee
By symmetry,
\be{A16}
Q_2 = \multsum{U\subset E_2}{1\in\ver U} (-1)^{|U| } {\bf x}_U.
\ee

We now consider $S_G$ and split the sum into various subsums. The
subsets $U\subset E$ fall into exactly one of the following eight
classes:
$$\begin{array}{lll}
\mbox{(I)}  \,\, 1\notin \ver U, 2\notin \ver U,
 \quad\quad  
\mbox{(II)} \,\, 1\in \ver U,  2 \notin \ver U, \quad \quad 
\mbox{(III)}  \,\, 1\notin\ver U,  2\in \ver U.
\end{array}$$
Any remaining set $U\subset E$ will have $\{1,2\} \subset \ver U$. Any
set $U$ that contains the edge $(1,2)\in E$ is of this type, and in
this case we write $U= (1,2) \cup U'$ with $U'= U\bs \{(1,2)\}$, and
sort these sets into the classes
$$
\begin{array}{lll}
\mbox{(IV)}&  U=\{(1,2)\} \cup U'; & 1,2 \notin \ver U',
\\
\mbox{(V)} &U=\{(1,2)\} \cup U'; & 1\in\ver U', \, 2\notin \ver U',
\\
\mbox{(VI)} &U=\{(1,2)\} \cup U';  & 1 \notin \ver U', \, 2\in\ver U',
\\
\mbox{(VII)} & U=\{(1,2)\} \cup U';& 1,2 \in\ver U'.
\end{array}
$$
Any set $U\subset E$ that has not been considered must have $1,2 \in
\ver U$, but $(1,2)\notin U$, and this is the characteristics of sets
in the final class (VIII). Accordingly, we sort the sets $U\subset E$ 
into these classes and obtain the decomposition
$$
S_G({\bf x}) = S^{\rm (I)} + S^{\rm (II)} + \ldots + S^{\rm (VIII)}
$$
where $S^{(\dagger)}$ is the subsum over sets $U$ in class
$(\dagger)$.\\

It remains to identify the various sums $S^{(\dagger)}$ with certain terms in
\rf{A14}. Note that a set $U\subset E$ with $1,2\notin \ver U$ is
actually a subset of $E_{1,2}$, and {\it vice versa}. By \rf{A13}, this
shows that $S^{\rm (I)} = S_H$, the third summand in \rf{A14}.

Next, consider a subset $U\subset E$ of class (II). The condition
that $2\notin \ver U$ is equivalent with $U\subset E_2$. Then the
further requirement that $1\in \ver U$ in conjunction with \rf{A16}
shows that $S^{\rm (II)}= Q_2$. By symmetry and \rf{A15}, we also have
$S^{\rm (III)}= Q_1$. This corresponds to the first two summands in
\rf{A14}.
 
Now consider a set $U= \{(1,2)\} \cup U'$ of class (V). Then, by
\rf{A11}, we have
$$
(-1)^{|U|} {\bf x}_U = -x_2(-1)^{|U'|} {\bf x}_{U'},
$$
and $U'$ is a set of class (II). By summing over $U'$, the results
of the previous paragraph yields $S^{\rm (V)}=-x_2 S^{\rm (II)}= -x_2Q_2$. By
symmetry, one also has $S^{\rm (VI)}=-x_1Q_1$. The same argument also
shows that $S^{\rm (IV)}= -x_1 x_2 S_H$. This identifies the fourth, fifth
and last summand on the right hand side of \rf{A14} as $S^{\rm (VI)},
S^{\rm (V)}$ and $S^{\rm (IV)}$, respectively.

Finally, consider a set $U= \{(1,2)\}\cup U'$ of class (VII). Then
$U'$ is of class (VIII). Inversely, when $U'$ is of class (VIII),
the set $\{(1,2)\} \cup U'$ is of class (VII). Since
$$
(-1)^{|U|}{\bf x}_U = -(-1)^{|U'|}{\bf x}_{U'},
$$
it follows that $S^{\rm (VII)} = -S^{\rm (VIII)}$. On collecting together, we
have now shown that the right hand side of \rf{A14} equals $S_G$,
completing the induction.

\section{Resolution of singularities}\label{secres}

In this chapter we prove Theorem \ref{thmres}. In preparation for the verification of Peyre's empirical formula in Chapter \ref{maninpeyre} we will also compute Peyre's alpha invariant \cite[Def. 2.4]{P} for the variety $X$ in Lemma \ref{alpha} and show in Lemma \ref{almostfano} that $X$ is an ``almost Fano" variety in the sense of \cite[Def.\ 3.1]{P2}.\\

%, and we also construct a similar resolution for the bi-projective threefold $\tilde{W} \subset \Bbb{P}^2 \times \Bbb{P}^2$ described in Section \ref{appl}. 
It amounts to no extra effort to work in more generality and consider varieties and schemes over an arbitrary (fixed) base field $k$. Let $n \geq 2$ be an arbitrary integer and let $W = W_n \subset \Bbb{P}^{2n-1}$  be the normal projective hypersurface with homogeneous coordinates $(x_1, x_2,\ldots , x_n, y_1, y_2,\ldots ,y_n)$ defined by the equation 
\begin{equation}\label{res1}
    x_1y_2y_3 \cdots  y_{n-1}y_n+ x_2y_1y_3\cdots y_n + \cdots + x_ny_1y_2\cdots y_{n-1}=0.
\end{equation}    
As in the special case $n=3$, $k = \Bbb{Q}$ considered in the introduction, there is a natural $G$-action on $W$ by a commutative algebraic group $G$: if $H \subset GL_2$ denotes the subgroup of upper triangular matrices with equal diagonal elements and $\Psi : H^n \rightarrow \Bbb{G}_a$ is the obvious generalization of  \eqref{hom},  then setting $G = {\rm ker}\Psi/\Bbb{G}_m$ we have a natural open immersion 
\begin{displaymath}
  j : G \rightarrow W.
\end{displaymath}   
and a natural action $\alpha: G \times W \rightarrow W$. 

The variety $W$ is of multiplicity $n-1$ along the $(n-1)$-plane $\Pi$ given by $y_1=y_2=\ldots = y_n=0$. Let $\Gamma = {\rm Bl}_{\Pi}\Bbb{P}^{2n-1}$  be the blow-up of $\Bbb{P}^{2n-1}$ along $\Pi$. It is the subvariety of $\Bbb{P}^{2n-1} \times \Bbb{P}^{n-1}$  with bi-homogeneous coordinates $(x_1, x_2,\ldots , x_n, y_1, y_2,\ldots,y_n ; Y_1, Y_2, \ldots, Y_n)$ defined by the equations
\begin{equation}\label{res2}
  y_iY_j - y_j Y_i = 0, \quad 1 \leq i < j \leq n.
\end{equation}
The restriction $p: \Gamma \rightarrow \Bbb{P}^{2n-1}$ of  $\Bbb{P}^{2n-1}\times \Bbb{P}^{n-1} \rightarrow \Bbb{P}^{2n-1}$    to $\Gamma$ gives a birational morphism, which is an isomorphism outside the exceptional divisor $D=p^{-1}(\Gamma)$ of $p$.  The blow-up $W' = {\rm Bl}_{\Pi}W$ of   $W$ along $\Pi$ is the closure of $p^{-1}(W\setminus \Pi)$ in $\Gamma$ (see \cite[Prop.\ IV-21]{EH}). Hence $W'$ is the subvariety of  $\Gamma$ defined by the equation
\begin{equation}\label{res3}
  x_1 Y_2Y_3 \cdots Y_{n-1}Y_n + x_2 Y_1 Y_3 \cdots Y_n + \ldots + x_nY_1 Y_2 \cdots Y_{n-1} = 0.
\end{equation}
Let $g : W' \rightarrow W$ be the restriction of $p$ to $W$. 
The inverse image scheme $p^{-1}(W) = W \times_{\Bbb{P}^{2n-1}} \Gamma$ of $W$ is the subscheme of $\Bbb{P}^{2n-1} \times \Bbb{P}^{n-1}$ defined by  \eqref{res1}. From \eqref{res2} and \eqref{res3} we deduce that 
\begin{displaymath}
  (x_1y_2 y_3 \cdots y_{n-1}y_n + x_2 y_1 \cdots y_{n-1} y_n + \ldots + x_n y_1 y_2 \cdots y_{n-1})Y_i^{n-1} = 0
\end{displaymath}
for all $1 \leq i \leq n$. Hence \eqref{res1} holds on $W'$, and it is therefore a subscheme of $p^{-1}(W)$. As $p^{-1}(W)$ and $W'$ are locally principal subschemes of $\Gamma$, we may also regard them as Cartier divisors (see \cite[p.\ 145]{Ha}). By our initial remark that $W$ is of multiplicity $n-1$ along  $\Pi$, we conclude that  $p^{-1}(W) = W'+(n-1)D$  (see \cite[Cor. 4.2.2  and Section 4.3]{ Fu}).\\  

We recall that for a normal variety $V$ with open immersion $j : U \rightarrow V$ of the non-singular locus $U$ of $V$ the canonical sheaf $\omega_V$ of $V$ is defined by $\omega_V = j_{\ast}(\Lambda^{\dim V}\Omega_U)$ where as usual $\Omega_U$ denotes the cotangent sheaf (cf.\ \cite{Re1b} and \cite[pp.\ 127-128, 180-182]{Ha} for the notation). 

\begin{lem} The canonical sheaves $\omega_W$ and $\omega_{W'}$ are invertible, and there is a canonical isomorphism $g^{\ast}\omega_W = \omega_{W'}$ of $O_{W'}$-modules. 
\end{lem}

\emph{Proof.} The inclusions $i: W \rightarrow  \Bbb{P}^{2n-1}$ and $j : W'\rightarrow \Gamma$  embed $W$ and $W'$ as locally principal subschemes of non-singular varieties. Hence by the adjunction formula, there are canonical isomorphisms $i^{\ast}\omega_{\Bbb{P}^{2n-1}}(W) =  \omega_W$  and  $j^{\ast}\omega_{\Gamma}(W') = \omega_{W'}$  induced by the Poincar\'e residue (cf.\ Ex.\ 25 in \cite[Chapter 3]{Re2}). This shows the invertibility of  $\omega_W$ and $\omega_{W'}$    as well as the equality
\begin{displaymath}
  g^{\ast} \omega_W = g^{\ast} i^{\ast}(\omega_{\Bbb{P}^{2n-1}}(W)) = j^{\ast} p^{\ast}(\omega_{\Bbb{P}^{2n-1}}(W)) = j^{\ast}(p^{\ast}(\omega_{\Bbb{P}^{2n-1}}) \otimes_{O_{\Gamma}} O_{\Gamma}(p^{-1}W)). 
\end{displaymath}
By \cite[p.\ 608]{GH} and the obvious fact $\dim \Bbb{P}^{2n-1} - \dim\Pi - 1 = n-1$ we have the canonical isomorphism
\begin{displaymath}
  p^{\ast} \omega_{\Bbb{P}^{2n-1}} = \omega(-(n-1)D) = \omega_{\Gamma} \otimes_{O_{\Gamma}} O_{\Gamma}(-(n-1)D). 
\end{displaymath}
Since $O_{\Gamma}(p^{-1}W) \otimes_{O_{\Gamma}} O_{\Gamma}(-(n-1)D) = O_{\Gamma}(p^{-1}W - (n-1)D))$ (cf.\ \cite[p.\ 144]{Ha}), we finally obtain canonical isomorphisms
\begin{displaymath}
  g^{\ast}\omega_W = j^{\ast}(\omega_{\Gamma} \otimes_{O_{\Gamma}} O_{\Gamma}(p^{-1}W - (n-1)D)) = j^{\ast}(\omega_{\Gamma} \otimes_{O_{\Gamma}} O_{\Gamma}(W')) = \omega_{W'}. 
\end{displaymath}
This completes the proof.\\

Now let $X \subset \Bbb{P}^{2n-1} \times \Bbb{P}^{n-1} \times \Bbb{P}^{n-1}$ be the tri-projective variety with tri-homogeneous coordinates $(x_1, \ldots, x_n, y_1, \ldots, y_n; Y_1, \ldots, Y_n; Z_1, \ldots, Z_n)$ defined by the equations
\begin{equation}\label{res4}
  x_1Z_1 + \ldots + x_nZ_n = 0,
\end{equation}
\begin{equation}\label{res5}
  y_iY_j - y_jY_i = 0, \quad 1 \leq i < j \leq n,
\end{equation}
\begin{equation}\label{res6}
   Y_1Z_1 = \ldots = Y_nZ_n.
\end{equation}
Then \eqref{res6} and \eqref{res4} imply
\begin{displaymath}
\begin{split}
  &(x_1Y_2Y_3 \cdots Y_n + x_2 Y_1 Y_3 \cdots Y_n + \ldots + x_n Y_1 \cdots Y_{n-1}) Z_i \\
  =  & (x_1Z_1 + \ldots + x_nZ_n) Y_1 \cdots Y_{i-1} Y_{i+1} \cdots Y_n = 0.
 \end{split} 
\end{displaymath}
In particular, \eqref{res3} holds on $X$, and the projection $\Bbb{P}^{2n-1} \times \Bbb{P}^{n-1} \times \Bbb{P}^{n-1} \rightarrow \Bbb{P}^{2n-1} \times \Bbb{P}^{n-1}$ onto the first two factors restricts to a morphism $h : X \rightarrow W'$. Let $B \subset \Bbb{P}^{n-1} \times \Bbb{P}^{n-1}$  be the bi-projective variety with bi-homogeneous coordinates $ (Y_1, \ldots , Y_n; Z_1,\ldots, Z_n)$ defined by \eqref{res6}. Then we have the following lemma.

\begin{lem}\label{lemmares2} {\rm (i)} The variety $B$ is a projective toric variety of dimension $n-1$, which is locally a complete intersection. Its singular locus is the union of all closed subsets defined by equations
\begin{displaymath}
     Y_i= Y_j=  Z_i=Z_j=0, \quad 1\leq i < j \leq n.
\end{displaymath}     
In particular, if $n\leq 3$, then $B$ is non-singular.\\
{\rm (ii)} The variety $X$ is a $\Bbb{P}^{n-1}$-bundle over $B$. 
\end{lem}

\emph{Proof.} (i) Let $B^{\circ}$ be  the open subset of $B$ where none of the $Y_i$ vanishes. By \eqref{res6} this is equivalent to the condition that  none of the  $Z_i$  vanishes, and we may regard $B^{\circ}$ as the underlying variety of the $(n-1)$-dimensional algebraic torus $T$  obtained as the quotient group of the diagonal embedding of  $\Bbb{G}_m$  in $\Bbb{G}_m^n$ .  There is a $\Bbb{G}_m^n$-action on $B^{\circ}$ given by
\begin{displaymath}
  (t_1, \ldots, t_n) \cdot (Y_1, \ldots, Y_n;Z_1, \ldots, Z_n) = (t_1Y_1, \ldots, t_nY_n; Z_1/t_1, \ldots, Z_n/t_n),
\end{displaymath}
and this action factorizes to a $T$-action $\rho : T\times B \rightarrow B$ such that its restriction $T\times  B^{\circ}\rightarrow  B^{\circ}$ is the group law on $T$.  Hence $B$ is a toric variety of dimension $n-1$. %It follows by induction on $n$ that  it is non-singular outside $B^{\circ}$.

We cover $B$ by open subsets $B_{k, l}$ ($1 \leq k \not= l \leq n$), defined by the conditions  $Y_k Z_l \not= 0$.  For $i \not=k$ and $j \not=l$  let $t_i=Y_i/Y_k$  and $u_j=Y_j/Y_l$. Using \eqref{res6}, we may then   eliminate $t_l$  and $u_k$ and identify $B_{k,l}$ with a subvariety of $\Bbb{A}^{2n-4}$  with affine coordinates $t_i , u_j$ for  $i, j \in \{1,2,\ldots, n\}\setminus \{k,l\}$, defined by    $n-3$ equations.  For instance, if  $k=n-1$, $l=n$, then $B_{k,l} \subset \Bbb{A}^{2n-4}$   is the complete intersection given by the equations 
\begin{displaymath}
  t_1u_1- t_2u_2 = \ldots = t_{n-3}u_{n-3}-t_{n-2}u_{n-2}=0 
\end{displaymath} 
   with singular locus given by the union of all closed subsets where $t_i= t_j= u_i= u_j=0$ for $  1\leq i< j \leq n-2$. Permuting indices, the same holds for general $k, l$.  \\

(ii) The projection $\Bbb{P}^{2n-1} \times \Bbb{P}^{n-1} \times \Bbb{P}^{n-1} \rightarrow  \Bbb{P}^{n-1} \times \Bbb{P}^{n-1} $  restricts to a surjective morphism from $X$ to $B$ where the fibre at a point $(Y_1,\ldots, Y_n; Z_1,\ldots, Z_n) \in B$ is the $(n-1)$-plane  in $\Bbb{P}^{2n-1}$ defined by \eqref{res4} and \eqref{res5}.\\

We are now prepared to prove the following more general version of Theorem \ref{thmres}: 
\begin{thm} The  projection $\Bbb{P}^{2n-1} \times \Bbb{P}^{n-1} \times \Bbb{P}^{n-1} \rightarrow  \Bbb{P}^{2n-1} $  restricts to a proper $G$-equivariant morphism $f: X\rightarrow W$ from a normal variety $X$ with $f^{\ast} \omega_W \cong \omega_X$. If $n=3$, then $f$  is a crepant  resolution of $W$. 
\end{thm}

\emph{Proof.} Let $h : X \rightarrow W'$ and $g : W' \rightarrow W$ be as above. Then $f = g h$. As we have already shown that $g$ is birational with $g^{\ast}\omega_W = \omega_{W'}$, it remains to show that $h$ is a crepant resolution of $W'$. It follows from Lemma \ref{lemmares2} that $X$ is non-singular. Let $U_0 \subset W'$ be the open subset with at most one $Y_j = 0$, and for $1 \leq i \leq n$ let $U_i \subset W'$ be the open subset where $x_i \not = 0$ and $Y_j = 0$ for at most one $j \not= i$. We may then identify local sections $s_i : U_i \rightarrow X$ of $h$ as follows. For a point $P = (x_1, \ldots x_n, y_1, \ldots, y_n; Y_1, \ldots, Y_n) \in U_i$ we define the coordinates $(Z_1, \ldots, Z_n)$ of $s_i(P )$ by
\begin{displaymath}
  (Y_2 Y_3 \cdots Y_n, Y_1Y_3 \cdots Y_n, \ldots, Y_1Y_2 \cdots Y_{n-1})  
\end{displaymath}
if $i = 0$, and for $1 \leq i \leq n$ by
\begin{displaymath}
\begin{split}
 &  Z_i = - \sum_{k \not = i} x_k \frac{Y_1 Y_2 \cdots Y_n}{Y_k Y_i}, \quad\quad\quad Z_j = x_i\frac{Y_1 Y_2 \cdots Y_n}{Y_j Y_i}, \quad j \not = i.
\end{split}  
\end{displaymath}
It is now easy to see that the morphisms glue to a section $s  : U= U_0 \cup U_1 \ldots \cup  U_n \rightarrow X$  
of $h: X \rightarrow W'$ and that $h$ maps $h^{-1}(U)$  isomorphically onto $U$. This shows that $h$ is birational 
and that the restrictions of  $h^{\ast}\omega_{W'}$  and $\omega_X$   to $h^{-1}(U)$ are isomorphic. Next we show   that the restriction from ${\rm Pic}(X)$ to ${\rm Pic}(h^{-1}(U))$ is bijective. It follows from Lemma \ref{lemmares2}  that $X$ is locally a complete intersection that is non-singular when $n=3$, and has singular locus $X_{\text{sing}}$ disjoint to $h^{-1}(U)$  and of codimension $\geq 4$ if $n\geq 4$.  The restriction from ${\rm Pic}(X)$ to ${\rm Pic}(X\setminus  X_{\text{sing}})$ is therefore bijective by a theorem of Grothendieck \cite[exp.\ XI, \S3]{SGA}. Now let $Z = X \setminus h^{-1}(U)$, and  for $1 \leq i < j < k \leq n$ let  $Z_{i, j, k} \subset Z$ be the subset where $Y_i=Y_j=Y_k = 0$ and for $1 \leq i < j \leq n$ let $Z_{i, j}\subset Z$ be the subset where $x_i = x_j = Y_i = Y_j = 0$. Then $Z$ is the union of all $Z_{i, j, k}$ and all $Z_{i, j}$, and we see that $Z = X \setminus h^{-1}(U)$ is of codimension at least 2 in $X$. Hence the restriction from ${\rm Pic}(X\setminus  X_{\text{sing}})$ to ${\rm Pic}(h^{-1}(U))$ is also bijective   (cf.\ \cite[Chapter II, Prop 6.5b and Cor. 6.16]{Ha}).

It remains to show that the resolution is $G$-invariant. The set $V=j(G)$ is the open subset of $W$ defined by $y_1y_2\cdots y_n \not=0$.  Hence $f^{-1}(V)$  is mapped isomorphically onto $V$ under $f$, and the inverse map is given by  $(x_1,\ldots, x_n, y_1, \ldots, y_n) \mapsto 
%(x_1,\ldots, x_n, y_1, \ldots, y_n; Y_1, \ldots, Y_n; Z_1, \ldots, Z_n)=
(x_1,\ldots, x_n, y_1, \ldots, y_n; y_1, \ldots, y_n; 1/y_1, \ldots , 1/y_n)$.  We may thus embed $G$ as an open subset of $X$, and there is  a natural $G$-action $\beta : G \times X \rightarrow X$, given by
\begin{displaymath}
\begin{split}
 & \left(\left(\begin{matrix}b_1 & a_1\\ 0& b_1\end{matrix}\right), \ldots,  \left(\begin{matrix}b_n & a_n\\ 0& b_n\end{matrix}\right) \right)\cdot (x_1, \ldots, x_n, y_1, \ldots, y_n; Y_1, \ldots, Y_n; Z_1, \ldots, Z_n)\\
  &=  ((b_1x_1 + a_1y_1, \ldots, b_nx_n + a_ny_n, b_1y_1, \ldots, b_ny_n; b_1y_1, \ldots, b_ny_n; z_1/b_1, \ldots, z_n/b_n). 
 \end{split} 
  \end{displaymath}
The restriction of $\beta$ to $G \times f^{-1}j(G)$  reduces to the group law on $G$,  and it is easy to see that  there is a commutative square 
  $$\begin{CD}
G \times X @>\beta>> X\\
@V{(id, f)}VV       @VVfV\\
G \times W @>\alpha>> W
\end{CD}$$
This completes the proof of the theorem.\\

The proof of Theorem \ref{thmthree} proceeds along similar lines: We recall the definition of $\tilde{W}$ and $\tilde{X}$ in Section \ref{appl}. Let $\tilde{B} \subset \Bbb{P}^2 \times \Bbb{P}^2$ be the subvariety with bi-homogeneous coordinates $(y_1, y_2, y_3; z_1, z_2, z_3)$ defined by \eqref{tri3}. It is the blow-up of $\Bbb{P}^2$ at the three points $(1, 0, 0)$, $(0, 1, 0)$, $(0, 0, 1)$ and hence a non-singular del Pezzo surface of degree $6$ with $\text{Pic}(\tilde{B}) \cong \Bbb{Z}^4$. Moreover, as the map $(\textbf{x}; \textbf{y}; \textbf{z})$ makes $\tilde{X}$ to a $\Bbb{P}^1$-bundle over $\tilde{B}$ (cf.\ \eqref{tri2}), we conclude that $\tilde{X}$ is non-singular and $\text{Pic}(\tilde{X}) \cong \Bbb{Z}^5$. We now consider the restriction to $\tilde{X}$ of the projection $\Bbb{P}^2 \times \Bbb{P}^2 \times\Bbb{P}^2\rightarrow \Bbb{P}^2 \times \Bbb{P}^2$ which sends $(\textbf{x}; \textbf{y}; \textbf{z})$ to $(\textbf{x}; \textbf{y})$. As
\begin{displaymath}
  (x_1z_1 + x_2z_2 + x_3z_3)y_{i+1}y_{i+2} = (x_1y_2y_3 + x_2y_1y_3 + x_3 y_1y_2) z_i
\end{displaymath}
for $1 \leq i \leq 3$ and with indices modulo 3, we obtain a morphism $\tilde{f} : \tilde{X} \rightarrow \tilde{W}$. If we let $\tilde{X}_0 \subset \tilde{X}$ and $\tilde{W}_0 \subset \tilde{W}$ be the open subsets where $(y_1y_2, y_1y_3, y_2y_3) \not = 0$, the $\tilde{f}$ restricts to an isomorphism $\tilde{f}_0 :  \tilde{X}_0 \rightarrow \tilde{W}_0$ with an inverse map $\tilde{g}_0 :\tilde{W}_0 \rightarrow \tilde{X}_0$ which sends $(\textbf{x}; \textbf{y})$ to $(\textbf{x}; \textbf{y}; \textbf{z})$ with $z_i = y_{i+1}y_{i+2}$ (again taking indices modulo 3). Hence $\tilde{f} :  \tilde{X} \rightarrow \tilde{W}$ is a desingularization of $\tilde{W}$. Note that $\tilde{W}$ is locally a complete intersection and hence a Gorenstein variety with invertible canonical sheaf $\omega_{\tilde{W}}$. As $\tilde{f}_0$ is an isomorphism, the canonical homomorphism $f^{\ast}  \omega_{\tilde{W}} \rightarrow \omega_{\tilde{X}}$ of $O_{\tilde{X}}$-modules restricts to an isomorphism of $O_{\tilde{X}_0}$-modules. Since $\tilde{X}$ is smooth and $\text{\rm codim}(\tilde{X} \setminus \tilde{X}_0) \geq 2$, the isomorphism class of an invertible $O_{\tilde{X}}$-module is uniquely determined by its restriction to $\tilde{X}_0$ (see \cite[Chap.\ II, 6.5, 6.11, 6.15]{Ha}). So $f^{\ast} \omega_{\tilde{W}} \cong \omega_{\tilde{X}}$ as $O_{\tilde{X}}$-modules, which means that $\tilde{f} : \tilde{X} \rightarrow \tilde{W}$ is crepant.\\

We return to our original set-up and proceed to   compute Peyre's alpha invariant \cite[Def 2.4, p.\ 120]{P} of the non-singular fourfold $X$ in the case $n=3$.   It is convenient to  
identify $G$ with the open subvariety $f^{-1}(j(G))\subset X$ where $%
y_{1}y_{2}y_{3}\neq 0$ and make use of the fact that $X$ is an equivariant
compactification of $G$. 

We introduce the following notation: let  $D_{0}$ be the subvariety of $X$ defined by $%
y_{1}=y_{2}=y_{3}=0$. If $1\leq i\leq 3$ and $\left\{ j,k\right\}
=\left\{ 1,2,3\right\} \backslash \left\{ i\right\}$,  we define $D_{i}$ $\subset
X$   by $Y_{j}=Y_{k}=0$ and $D_{i+3}\subset X$ by $Z_{j}=Z_{k}=0$%
.  Then $D_{0}$ is a $\Bbb{P}^{1}$-bundle over $B$ contained in $X\backslash G$,
while $D_{1},\ldots ,D_{6}$ are the inverse images of the exceptional curves on
the toric del Pezzo surface $B$. It follows  that $X\backslash G=D_{0}\cup
D_{1}\cup ...\cup D_{6}$ and we shall write $\text{Div}_{X\backslash
G}X=\sum_{i=0}^{6}\Bbb{Z}D_{i}$ for the free abelian group of divisors
with support in $X\backslash G$. If $D$ is a divisor  on $X$, we
will write $\left[ D\right] $ for its class in $\text{Pic}( X)$ and $C_{\text{eff}%
}(X)\subset  \text{Pic}(X)\otimes _{\Bbb{Z}}\Bbb{R}$ for the
pseudo-effective cone spanned by the classes of the effective
divisors.

\begin{lem}\label{lemx} {\rm (i)} The canonical homomorphism from  ${\rm Div}
_{X\backslash G}X$ to ${\rm Pic}(X)$ is surjective and  its kernel is 
the subgroup generated by $D_{2}-D_{1}+D_{4}-D_{5}$ and $
D_{3}-D_{1}+D_{4}-D_{6}$. In particular,
\begin{equation}\label{picard}
  \text{\rm rk } \text{\rm Pic}(X) = 5.
\end{equation}
{\rm (ii)} Any element in $C_{\text{eff}}(X)$ is equal to 
$\sum_{i=0}^{6}\lambda _{i}\left[ D_{i}\right] $ for some non-negative 
real numbers $\lambda _{0},...,\lambda _{7}$.\\
{\rm (iii)} If $1\leq i\leq 3$ and $\{ j,k\} = \{ 1,2,3\}
\backslash \{ i \} ,$ then $3(D_{0}+D_{j}+D_{k}+D_{i+3})$ is 
an anticanonical divisor. 
\end{lem} 

\emph{Proof.} (i) Let $\text{Hom}(G, \Bbb{G}_{m})$ be the character
group of $G$. By \cite[Prop.\ 1.1]{TT} there is  a natural exact
sequence 
$$0\rightarrow \text{Hom}(G, \Bbb{G}_{m}) \rightarrow  \text{Div}_{X\backslash
G}X\rightarrow  \text{Pic} (X)\rightarrow  0$$
where the map from Hom($G$, $\Bbb{G}_{m}$) is the usual divisor map of
rational functions. Hence, as Hom($G$, $\Bbb{G}_{m}$) is a free abelian
group generated by $y_{1}/y_{2}=Y_{1}/Y_{2}$ and $y_{1}/y_{3}=Y_{1}/Y_{3}$,
it suffices to note that div($Y_{1}/Y_{2})=D_{2}-D_{1}+D_{4}-D_{5}$ and div$%
(Y_{1}/Y_{3})=D_{3}-D_{1}+D_{4}-D_{6}$ to get the desired assertion.\\
(ii)  This is a special case of \cite[Prop.\ 1.1(3)]{TT}.\\
(iii) As $f:X\rightarrow W$ is crepant and $\omega _{W}\cong  O_{W}(-3)$%
, it suffices to show that the principal closed subscheme of $X$ defined by $%
y_{i}$ gives rise to the divisor $D_{0}+D_{j}+D_{k}+D_{i+3}$ (cf.\ \cite[II.6.17.1]{Ha}). To see this, we first note that $y_{i}$ has multiplicity 1
along $D_{0}$. It is therefore enough to prove that the closed subscheme of $X$
defined by $Y_{i}$ has divisor $D_{j}+D_{k}+D_{i+3}$, which is easy to check
by computing the divisor of $Y_{i}=0$ on $B$. \\

Now let $C_{\text{eff}}(X)^{\vee }\subset \text{Hom}(\text{Pic} (X)\otimes \Bbb{R}, 
\Bbb{R})$ be the dual cone of all linear maps $\Lambda : \text{Pic} (X)\otimes  
\Bbb{R}\rightarrow \Bbb{R}$ such that $\Lambda (\left[ D\right]
)\geq 0$ for every effective divisor $D$ on $X$. Moreover, let $l : \text{Hom}(\text{Pic}(X)\otimes 
\Bbb{R}, \Bbb{R})\rightarrow \Bbb{R}$ be the linear map which
sends $\Lambda $ to $\Lambda (\left[ -K_{X}\right] )$. We then endow $\text{Hom}(\text{Pic} 
(X)\otimes \Bbb{R}, \Bbb{R})$ with the Lebesque measure $ds$
normalized such that $L=\text{Hom}(\text{Pic}(X), \Bbb{Z})$ has covolume 1, and $%
H_{X}=l^{-1}(1)$ with the measure $ds/d(l-1)$. If $z_{1},\ldots ,z_{r}$ are
coordinates for $\text{Hom}(\text{Pic}(X)\otimes \Bbb{R},  \Bbb{R}) =\Bbb{R}^{r}$
with respect to a $\Bbb{Z}$-basis of $L$ and $%
l(z_{1},\ldots ,z_{r})=\alpha _{1}z_{1}+\ldots +\alpha _{r}z_{r}$, then $%
ds/d(l-1)=dz_{1}...\widehat{dz_{i}}...dz_{r}/\left\vert \alpha
_{i}\right\vert $ whenever $\alpha _{i}\neq 0$. After these preparations, we may now define $\alpha (X)
$ as 
\begin{equation}\label{defalpha}\alpha (X)=\int\limits_{C_{\text{eff}}(X)^{\vee
}\cap H_{X}}\frac{ds}{d(l-1)}.\end{equation}

The following result evaluates the alpha invariant explicitly.

\begin{lem}\label{alpha} One has
 $$\alpha (X)=\frac{1}{2^{4}3^{5}}=\frac{1}{3888}.$$
\end{lem}

\emph{Proof.}  By part (i) of the preceding lemma, the classes of $%
D_{0},D_{1},D_{2},D_{3},D_{4}$ form a $\Bbb{Z}$-basis of $\text{Pic} (X)$%
. Let $e_{0},e_{1},e_{2},e_{3},e_{4}$ be the dual $\Bbb{Z}$-%
basis of $L$ with $e_{i}([ D_{j}] )=\delta _{ij}$ and $%
(z_{0},z_{1},z_{2},z_{3},z_{4})$ be the coordinates of $\text{Hom}(\text{Pic}(X)\otimes 
\Bbb{R}, \Bbb{R})$, with respect to this basis. Then, by part (ii)
and (iii) of Lemma \ref{lemx}, we have that $C_{\text{eff}}(X)^{\vee }$ is the
subset of $\Bbb{R}_{\geq 0}^{5}$  defined by the inequalites $%
z_{2}+z_{4}-z_{1}\geq 0$ and $z_{3}+z_{4}-z_{1}\geq 0$, and that $H_{X}$ is
the hyperplane in $\Bbb{R}^{5}$ defined by $3z_{0}+3z_{2}+3z_{3}+3z_{4}=1$. Hence $0\leq z_{0}\leq \frac{1}{3}$ on $H_{X}$ and 

$$ \alpha (X)=\int_0^{1/3}\left(
\int \int \int \int_{\Delta }\frac{dz_{1}dz_{2}dz_{3}dz_{4}}{%
d(l_{0}-(1-3z_{0}))}\right) dz_{0}$$
for $l_{0}(z_{1},...,z_{4})=3z_{2}+3z_{3}+3z_{4}$ and the subset $\Delta
\subset  \Bbb{R}_{\geq 0}^{4}$ defined by $z_{2}+z_{4}-z_{1}\geq 0$ and 
$z_{3}+z_{4}-z_{1}\geq 0$. To compute the inner integral, we   substitute $\tilde{z}_{i} = z_{i}/(1-3z_{0})$ for $1\leq i\leq
4.$ Then $$l_{0}(\tilde{z}_{1} ,\ldots ,\tilde{z}_{4})-1=\frac{%
l_{0}(z_{1},\ldots ,z_{4})-(1-3z_{0})}{1-3z_{0}}$$ so that 
%$\int \int \int \int_{\Delta }\frac{dz_{1}dz_{2}dz_{3}dz_{4}}{%
%d(l_{0}-(1-3z_{0}))}=(1-3z_{0})^{3}\int_{\Delta }\frac{%
%dz_{1}dz_{2}dz_{3}dz_{4}}{d(l_{0}(z_{1},...,z_{4})-1)}$ \ \ and\bigskip 
$$\alpha (X)=\int_{0}^{1/3}(1-3z_{0})^{3}dz_{0}\int \int \int \int_{\Delta }%
\frac{dz_{1}dz_{2}dz_{3}dz_{4}}{d(l_{0}(z_{1},\ldots,z_{4})-1)}.$$
To compute the  last integral, we use the equation $l_{0}(z_{1},...,z_{4})-1=0$
on $H_{1}$ to eliminate $z_{4}$. In this way we see that the  multiple  integral over $\Delta$ %$\int \int \int
%\int_{\Delta }\frac{dz_{1}dz_{2}dz_{3}dz_{4}}{d(l_{0}(z_{1},...,z_{4})-1)}=%
equals $\frac{1}{3}\text{Vol}( \Pi) $ for the subset$\ \Pi \subset \Bbb{R}%
_{\geq 0}^{3}$ of the first octant defined by $3z_{1}+3z_{2}\leq 1$, $%
3z_{1}+3z_{3}\leq 1$ and $3z_{2}+3z_{3}\leq 1$. As Vol$\left( \Pi \right) =%
\frac{1}{108}=\frac{1}{2^{2}3^{3}}$ and $\int_{0}^{1/3}(1-3z_{0})^{3}dz_{0}=%
\frac{1}{2^{2}3}$, we obtain $\alpha (X)=\frac{1}{2^{4}3^{5}}=\frac{%
1}{3888}$. \\

We conclude this chapter by  showing that $X$ satisfies the three conditions in
\cite[Def. 3.1]{P2}   for being an ``almost Fano" variety.  

\begin{lem}\label{almostfano}  Let $X\subset  \Bbb{P}%
^{5}\times \Bbb{P}^{2}\times \Bbb{P}^{2}$ be the fourfold %
defined by \eqref{res4}--\eqref{res6}  over some field $k$. Then the 
following holds:\\
{\rm (i)} $H^{1}(X,O_{X})=H^{2}(X,O_{X})=0.$\\
{\rm (ii)} The geometric Picard group $ \text{Pic}(\overline{k}\times X)$ is %
torsion-free.\\
{\rm (iii)} The anticanonical class is in the interior of $C_{%
\text{eff}}(X)$.
\end{lem}

\emph{Proof.}  (i) Use Lemma \ref{lemmares2}(ii) and the fact that $%
H^{1}(B,O_{B})=H^{2}(B,O_{B})=0$.\\
(ii) By Lemma \ref{lemx}(i) we have that $\text{Pic}(K\times X)= \Bbb{Z}^5$ for any
field $K\supset k$.\\
(iii) This follows from Lemma \ref{lemx}(ii) and  (iii).\\

\emph{Remark:} The variety $\tilde{X}$ featured in the proof of Theorem \ref{thmthree} is also easily seen to be ``almost Fano": as $\tilde{f} : \tilde{X} \rightarrow \tilde{W}$ is a crepant resolution and $\omega_{\tilde{W}}^{-1}$ is ample, $\omega_{\tilde{X}}^{-1}$ is
a big $O_{\tilde{X}}$-module \cite[Def.\ 2.2.1]{La} and hence $[ÐK_{\tilde{X}}]$  in the interior of $C_{\text{eff}}(\tilde{X})$ by \cite[Th.\ 2.2.25]{La}. This proves (iii), while (i) follows as in the previous lemma and (ii) from the fact that $\tilde{X}$ is a $\Bbb{P}^1$-bundle over $\tilde{B}$.

\section{The descent variety}\label{descent}

In this chapter, we show that the cubic can be parametrized. We start with simple divisibility considerations that resemble the argument  in  \cite{BB}. In the following section we then show that the descent variety so obtained is   the universal torsor. 

\subsection{An elementary argument}\label{torsor}
Let $\cal W$ denote the set of integer solutions to \eqref{1} and \eqref{nonzero} with no further coprimality conditions. % $({\mathbf x}, {\mathbf y})$
%of the conditions
%$$  x_1y_2y_3 + x_2y_3y_1 + x_3y_1y_2 = 0, \quad y_1y_2y_3 \neq 0. $$
As a first preparatory step we
will link this with 
the bilinear equation
\be{2.1}
u_1 v_1 + u_2 v_2 + u_3 v_3 =0.
\ee
For $({\bf x}, {\bf y})\in\cal W$ 
%with $y_1y_2y_3\neq 0$. Then 
we put
$$ u= (y_1;y_2;y_3),\; u_1=(y_2/u;y_3/u),\; u_2=(y_3/u;y_1/u),\;
u_3=(y_1/u;y_2/u) $$
and observe that
\be{2.2}
(u_1;u_2)=(u_2;u_3)=(u_3;u_1)=1.
\ee
Hence we can write
\be{2.3} y_1=uu_2u_3w_1, \quad y_2= uu_1u_3w_2, \quad y_3 = uu_1u_2w_3
\ee
with integers $w_j\neq 0$, and the equation \rf{1} now reads
\be{2.4} u_1x_1w_2w_3 + u_2x_2w_1w_3 + u_3x_3 w_1w_2 = 0. 
\ee
By construction, the coprimality conditions
\be{2.5}
(u_j;w_j)=(w_1;w_2)=(w_2;w_3)=(w_3;w_1)=1 \quad (1\le j \le 3) 
\ee
hold in addition to \rf{2.2}. By \rf{2.4}, we see that
$ w_1 \mid u_1x_1w_2w_3 $, and \rf{2.5} then implies that $w_1\mid x_1$.
By symmetry, it also follows that $w_2\mid x_2$, $w_3\mid x_3$, and we write
\be{2.6} x_j=w_jv_j \ee
with $v_j\in\ZZ$. In this notation, \rf{2.4} reduces to \rf{2.1}. These
transformations can be reversed: if natural numbers $u,u_1,u_2,u_3$
and
integers $v_j,w_j$ are given, then the numbers $x_j,y_j$ defined by
\rf{2.6} and \rf{2.3} satisfy \rf{1}. In particular, this proves the
following.

\begin{lem}\label{lem1} Let $\mathcal A$ denote the set of all $10$-tuples
  $u,u_1,u_2,u_3, v_1,v_2,v_3,w_1,w_2,w_3$ with $u,u_j\in\NN$,
  $w_j\in\ZZ\setminus \{0\}$, $v_j\in\ZZ$ that satisfy \rf{2.1},
  \rf{2.2} and \rf{2.5}. Then the map ${\mathcal A} \to \ZZ^6$ defined
  by \rf{2.3} and \rf{2.6} is a bijection between $\mathcal A$ and $\cal W$.
\end{lem}

%It follows that $V(P )$ equals the number of 10-tuples $(d,{\mathbf d}, {\mathbf z},{\mathbf a})\in{\mathcal A}$ with \be{2.7} dd_1d_2|z_3|\le P, \quad dd_2d_3|z_1|\le P, \quad dd_1d_3|z_2|\le P,\quad |a_jz_j|\le P\quad (1\le j\le 3). \eeThese conditions are independent of the signs of the $z_j$. Hence, if $U(P)$ is the number of elements in $\mathcal A$ with ${\mathbf z}\in\NN^3$ and \rf{2.7}, then\be{2.8} V(P)=8U(P). \ee This is our basic reinterpretation of $V(P)$. We are reduced to counting solutions of the equation \rf{2.1} with convolution type weights induced by \rf{2.7} and coprimality constraints \rf{2.2} and \rf{2.5}.\\

\def\bfd{{\mathbf u}}\def\bfa{{\mathbf v}}

For the next step, consider $\bfd\in{\Bbb N}^3$ as fixed, and 
study the set ${\cal L}(\bfd)$ of solutions
$\bfa\in {\Bbb Z}^3$ of \rf{2.1} as a lattice. 
For any integers $r_1, r_2, r_3$, the numbers
\be{2.9}
v_1 = {u_2 r_3} - {u_3 r_2},\quad 
v_2 = {u_3 r_1} - {u_1 r_3},\quad
v_3 = {u_1 r_2} - {u_2 r_1}
\ee
are a solution of \rf{2.1}. Fix a complete 
set  $\cal S$ of residues modulo $u_1$, and consider \rf{2.9}  as a map
$$  {\cal S}  \times {\Bbb Z}\times{\Bbb Z}
\to {\cal L}(\bfd),\quad {\mathbf r} \mapsto \bfa. $$
If \rf{2.2} holds, this map is actually a bijection. 
This fact is certainly well known, but we 
include the simple proof for completeness: 
Suppose that ${\mathbf r}, {\mathbf r}'$
with $r_1,r'_1\in\cal S$
map to the same $\bfa \in {\cal L}(\bfd)$. Then by \rf{2.9} for $v_3$,
one finds that $u_2r_1 \equiv u_2r'_1 \bmod u_1$, and hence that $r_1=r'_1$.
By \rf{2.9} again, it is now immediate that  ${\mathbf r}= {\mathbf r}'$, as
required to show that the map is injective. To show that the map is surjective,
let $\bfa\in{\cal L}(\bfd)$. By \rf{2.2}, there are integers $a,b$ with
$v_1=u_2a-u_3b$. Then, for any $k\in\Bbb Z$, one has
$$ v_1= u_2(a+ku_3) - u_3(b+ku_2) . $$
Similarly,  there are integers $r_1, r_3$ with 
$ v_2 = u_3r_1 - u_1r_3 $. Injecting these expressions into \rf{2.1}, we
deduce that $u_1u_2(a-r_3) \equiv u_1v_1 + u_2v_2 \equiv 0 \bmod u_3$, and 
hence we
may chose $k$ such that $r_3= a + ku_3$. With this choice, we put $r_2
= b+ku_2$. Then, by construction, the first two equations
in \rf{2.9} hold. The third equation must then also hold, because
$\mathbf r$ maps to the solution of \rf{2.1} with given values $v_1,v_2$.    
This shows that any solution of \rf{2.1} can be
written as in \rf{2.9}, for some ${\mathbf r}\in{\Bbb Z}^3$ . 
For any $j\in{\Bbb Z}$, the transformation
$
(r_1, r_2, r_3) \mapsto (r_1 +j{u_1}, \, r_2 +
j{u_2}, \, r_3 + j{u_3})
$
leaves \rf{2.9} invariant. Hence, an appropriate choice of $j$
guarantees that $r_1\in \cal S$, as required. This last invariance property
also shows that whenever $r_1$, $r'_1$ are natural numbers with
$r_1\equiv r'_1\bmod u_1$, then the sets 
${\cal R}(r_1)=\{(r_1,r_2,r_3): r_2,r_3\in{\Bbb Z}\}$ and 
${\cal R}(r'_1)$ are mapped to the same image.

We may now use \rf{2.9} within the conclusion of Lemma \ref{lem1}.
This yields the following.

\begin{lem}\label{lem2} Let ${\mathcal S}(q)$ denote a
  complete set of residues modulo $q$. Let $\mathcal B$ denote the set
of all $10$-tuples   $u,u_1,u_2,u_3, w_1,w_2,w_3,r_1,r_2,r_3$ with $u,u_j\in\NN$,
  $w_j\in\ZZ\setminus \{0\}$, $r_1\in {\cal S}(d_1)$, $r_2\in\ZZ$, 
$r_3\in\ZZ$.  
that satisfy
  \rf{2.2} and \rf{2.5}. Then the map ${\mathcal B} \to \ZZ^6$ defined
  by \rf{2.3} and
$$ 
x_1 = w_1({u_2 r_3} - {u_3 r_2}),\quad 
x_2 = w_2({u_3 r_1} - {u_1 r_3}),\quad
x_3 = w_3({u_1 r_2} - {u_2 r_1})
$$
 is a bijection between $\mathcal B$ and $\cal W$. 
\end{lem}

It will be relevant later 
to know that products $r_i u_j w_k$ with $\{i,j,k\} = \{1,2,3\}$ are not much 
larger than the original variables $x_j$, $y_j$. 
The following lemma makes this precise.

\begin{lem}\label{lem2A} Let $({\mathbf x,\mathbf y})\in\cal W$  with
$ |x_j| \le P$, $|y_j|\le P$ for $1\le j\le 3$. Suppose that
${\mathcal S}(u_1) \subset [1,2u_1]$  and 
 \begin{equation}\label{assump}
   |w_1u_2u_3| \leq 2\min(|w_2u_1u_3|, |w_3u_1u_2|).
\end{equation}   
     Then one has
\begin{equation}\label{2.10}
|r_1u_2 w_3|\le 2P, \quad |r_1 u_3 w_2|\le 2P,\quad
| r_2 u_1 w_3| \le 3P, \quad
|r_3 u_1 w_2| \le 3P,\
\end{equation}
\begin{equation}\label{2.10a}
|r_2u_3 w_1| \le 7P, \quad |r_3u_2w_1| \le 7P.
\end{equation}
\end{lem}

{\em Proof}. By Lemma \ref{lem2} and \rf{2.3}, the conditions 
$|x_j|,|y_j| \le P$ may be rewritten as the six constraints
\be{lem2a}
uu_1 u_2 w_3 \le P, \quad u u_2 u_3 w_1 \le P, \quad u u_1 u_3 w_2
\le P,
\ee
and
\be{lem2b}
| {u_2 r_3} - {u_3 r_2}| \le 
  \f{P}{|w_1|}, \quad
| u_3 r_1 - {u_1r_3}| \le  \f{P}{|w_2|}, \quad
| {u_1 r_2} - u_2 r_1| \le  \f{P}{|w_3|}. 
\ee
Using the bound $1\le  r_1 \le 2u_1$, one notes that \eqref{lem2a} 
implies the first two bounds in 
\eqref{2.10}. These together with \eqref{lem2b} imply the rest of \eqref{2.10}. 
With \eqref{2.10} in hand, one applies  \eqref{assump} and bounds the
minimum by the geometric mean to conclude that
\begin{displaymath}
  \min(|r_3u_2 w_1|,  |r_2u_3w_1|) \leq \sqrt{|r_2r_3u_2u_3w_1^2|} 
\leq  \frac{3P|w_1 u_2u_3|}{\sqrt{|w_2u_3u_1 w_3u_2u_1|}} \leq 6P. 
\end{displaymath}
Appealing to \eqref{lem2b} once again, we derive \eqref{2.10a}. %that
 %$\max(|r_2d_3 z_1|,  |r_3d_2z_1|)  \leq 7P$, as required. 

 \subsection{The universal torsor}\label{universaltorsor} In Lemma \ref{lem1} we proved a useful parametrization of the cubic \eqref{1} in an elementary \emph{ad hoc} fashion. In this section we take a very different route, and use much more sophisticated tools, to compute the universal torsor of the variety \eqref{1} by applying the general theory of Colliot-Th\'el\`ene and Sansuc \cite{CS}. The main result of this section  is Theorem \ref{thm7}. As a corollary we obtain a new proof of Lemma \ref{lem1} which is contained in the equivalent companion Lemma \ref{lemtorsor} below.\\

Let $K$ be a perfect field with algebraic closure $\bar{K}$,  $\mathfrak{g} = \text{Gal}(\bar{K}/K)$ and $\bar{V} = \bar{K} \times_K V$ for a variety $V$ over $K$. We recall that an $X$-torsor over a $K$-variety $X$ under a $K$-torus $T$ is a principal homogeneous space $%\phi_{\mathcal{R}} : 
\mathcal{T} %\overset{\phi}{\long\
\rightarrow X$ under $T$ (see \cite[Ch. III, \S 4]{Mi}). The isomorphism classes $[\mathcal{T}]$ of $X$-torsors under $T$ are in bijection with elements of $H^1_{\text{et}}(X, T)$ (cf.\ \cite[Section 1.2]{CS}), and if $X$ is a smooth, geometrically integral $K$-variety with $H^0_{\text{et}}(X, \Bbb{G}_m) = K^{\ast}$, then there is a natural exact sequence \cite[2.0.2]{CS}
\begin{equation}\label{tor1}
 0 \longrightarrow H^1_{\text{et}}(K, T) \longrightarrow H^1_{\text{et}}(X, T) \overset{\chi}{\longrightarrow} \text{Hom}_{\mathfrak{g}}(\widehat{T}, \text{Pic} (\bar{X}))
\end{equation}
where $\chi([\mathcal{T}]) \in \text{Hom}_{\mathfrak{g}}(\widehat{T}, \text{Pic}(\bar{X}))$ sends a character $\Psi : \bar{T} \rightarrow \Bbb{G}_{m, \bar{K}}$ to the $\bar{X}$-torsor $\overline{\mathcal{T}} \times^{\bar{T}} \Bbb{G}_{m, \bar{K}}$ under $\Bbb{G}_{m, \bar{K}}$, which one obtains from $\overline{\mathcal{T}}$ and $\Psi$ by changing the structure group of the torsor from $\bar{T}$ to $\Bbb{G}_{m, \bar{K}}$. The image $\chi([\mathcal{T}]) \in  \text{Hom}_{\mathfrak{g}}(\widehat{T}, \text{Pic} (\bar{X}))$ is called the \emph{type} of the torsor. If the $K$-torus $T$ is split, then $H^1_{\text{et}}(K, T) = 0$ by Hilbert's Theorem 90 (cf.\ \cite[III.4.9]{Mi}). As all $K$-tori in this paper are split, we shall therefore (by a slight abuse of language) refer to \emph{the} $X$-torsor of a certain type.

Now suppose that $\text{Pic}(\bar{X})$ is finitely generated and torsion free and that $T$ is the dual torus with character group $\widehat{T} = \text{Pic}(\bar{X})$ (here we identify canonically isomorphic $\mathfrak{g}$-modules). Then an $X$-torsor $\mathcal{T}$ under $T$ is said to be \emph{universal} if $ \chi([\mathcal{T}])  : \widehat{T} \rightarrow  \text{Pic} (\bar{X})$ is the identity map. It is known  \cite[2.2.9]{CS} that a universal torsor exists whenever $X(K) \not= \emptyset$. \\

We specialize now to the situation relevant in our case. Let  $K = \Bbb{Q}$ and $X \subset \Bbb{P}^5 \times \Bbb{P}^2 \times \Bbb{P}^2$ be the fourfold given by  \eqref{res4} -- \eqref{res6}. It is a hypersurface in the fivefold $\Xi \subset \Bbb{P}^5 \times \Bbb{P}^2 \times \Bbb{P}^2$ defined by \eqref{res5} and \eqref{res6}. The projection $ \Bbb{P}^5 \times \Bbb{P}^2 \times \Bbb{P}^2 \rightarrow  \Bbb{P}^2 \times \Bbb{P}^2$ restricts to morphisms $\lambda : X \rightarrow B$  and $\gamma : \Xi \rightarrow B$, which makes $X$ a $\Bbb{P}^2$-bundle and $\Xi$ a $\Bbb{P}^3$-bundle over the surface $B \subset \Bbb{P}^2 \times \Bbb{P}^2$ defined by \eqref{res6}. 

As a first step, we will describe the universal torsor over $\Xi$, which is  a (split) smooth projective toric variety: the torus  is the open subset $U \subset \Xi$ where all coordinates are different from zero and the $U$-action $U \times \Xi \rightarrow \Xi$ is given by coordinate-wise multiplication of the two 12-tuples representing points in $U$ and $\Xi$. It was shown in \cite[Prop.\ 8.5]{Sa1} that the universal torsor $ \mathcal{T} $ of a split smooth projective toric variety $\Xi$ coincides   with the toric morphism from the open toric subvariety $\Bbb{A}^n \setminus F$  of $\Bbb{A}^n$ described by Cox in \cite{Co}. Here the $n$ affine coordinates $t_{\rho}$ of $\mathcal{T} \subset \Bbb{A}^n$ are indexed by the one-dimensional cones (or edges) of the fan $\Delta$ of $\Xi$ (see \cite{Fu2}) and $F \subset \Bbb{A}^n$ is the closed subset defined by the monomials $t^{\sigma} = \prod_{\rho \not\in \sigma(1)} t_{\rho}$   for the maximal cones $\sigma$ of $\Delta$.

\begin{lem}\label{thm6} Let $\Omega \subset \Bbb{A}^{10}$ be the open subvariety with coordinates $(\xi_0, \xi_1, \xi_2, \xi_3, u_1, u_2, u_3, w_1, w_2, w_3)$ defined by
\begin{equation}\label{416}
\begin{split}
& u_iu_kw_jw_k \not= 0 \text{ for at least one triple }\{i, j, k\} = \{1, 2, 3\},\\
& (\xi_0, \xi_1, \xi_2, \xi_3) \not= (0, 0, 0, 0).
\end{split}
\end{equation}
Let $\phi : \Omega \rightarrow \Xi$ be the morphism which sends $(\xi_0, \xi_1, \xi_2, \xi_3, u_1, u_2, u_3, w_1, w_2, w_3)$ to
\begin{equation}\label{418}
\begin{split}
(x_1, x_2, x_3, y_1, y_2, y_3) & = (\xi_1, \xi_2, \xi_3, \xi_0u_2u_3w_1, \xi_0u_1u_3w_2, \xi_0u_1u_2w_3),\\
(Y_1, Y_2, Y_3; Z_1, Z_2, Z_3) & = (u_2u_3w_1, u_1u_3w_2, u_1u_2w_3; u_1w_2w_3, u_2w_1w_3, u_3w_1w_2).
\end{split}
\end{equation}
Then $\phi : \Omega \rightarrow \Xi$ is the underlying morphism of a universal torsor over $\Xi$.
\end{lem}
 
\emph{Proof.} For a cone $\tau$ of $\Delta$ under the action of the torus $U$ of $\Xi$ let $V(\tau)$ be the   closure of the orbit $O_{\tau}$, see \cite[Section 3.1]{Fu2}. There is a bijection between edges $\rho \in \Delta$ and irreducible components $D_{\rho} = V(\rho)$ of $\Xi \setminus U$, and there is also a bijection between maximal cones $\sigma \in \Delta$ and fixed points $P_{\sigma} = V(\sigma) = O_{\sigma} \in \Xi\setminus U$ under the action of $U$. Moreover, we have $\rho \in \sigma(1)$ if and only if $P_{\sigma} \in D_{\rho}$. There are ten irreducible components  $D_{\rho}$ of $\Xi \backslash U$. For $i = 1, 2, 3$ we let $\xi_i$ correspond to the prime divisor where $x_i = 0$ and $\xi_0$ to the prime divisor where $y_1=y_2=y_3 = 0$. For a triple $\{i, j, k\} = \{1,2, 3\}$ we let $u_j$ be the coordinate corresponding to the prime divisor where $Y_i=Y_k = Z_j = 0$ and $w_k$ the coordinate corresponding to the prime divisor where $Y_k=Z_i=Z_j=0$. There is then a natural embedding of the universal $\Xi$-torsor $\mathcal{T}$ in the affine space $\Bbb{A}^{10}$ with coordinates $(\xi_0, \xi_1, \xi_2, \xi_3, u_1, u_2, u_3, w_1, w_2, w_3)$. 

A point on $\Xi \subset \Bbb{P}^5 \times \Bbb{P}^2 \times \Bbb{P}^2$ is a fixed point $P_{\sigma}$ under $U$ if and only if its image has exactly one non-zero coordinate under each of the projections $\text{pr}_1 : \Xi \rightarrow \Bbb{P}^5$, $\text{pr}_2 : \Xi \rightarrow \Bbb{P}^2$ and $\text{pr}_3 : \Xi \rightarrow \Bbb{P}^2$. Such a point  either satisfies $y_jY_jZ_i(P_{\sigma}) \not= 0$ or $z_lY_jZ_i(P_{\sigma}) \not= 0$ for $i \not= j$ and $1 \leq l \leq 3$. In the first case,  $P_{\sigma}$ does  not lie on the divisors corresponding to $\xi_0, u_i, u_k, w_j, w_k$ and in the second case, $P_{\sigma}$ does not lie on the divisors corresponding to $\xi_l, u_i, u_k, w_j, w_k$. The exceptional set $F \subset \Bbb{A}^{10}$ of Cox is thus given by the monomials $\xi_lu_iu_kw_jw_k$ where $(i, j, k, l)$ runs over all quadruples with $\{i, j, k\} = \{1, 2, 3\}$ and $0 \leq l \leq 3$. Hence $\mathcal{T} = \Bbb{A}^{10} \setminus F$ is just the open subset $\Omega \subset \Bbb{A}^{10}$ defined in \eqref{416}. 

The structure morphism $\phi  : \mathcal{T} \rightarrow \Xi$ of the universal torsor is given in terms of fans in \cite{Co}  and \cite[Prop.\ 8.5]{Sa1}; it follows from this or from the general local description of torsors \cite[Section 2.3]{CS} that the restriction of $\phi $ to $\Bbb{G}_{m, \Bbb{Q}}^{10}$ is the homomorphism of tori $\Bbb{G}_{m, \Bbb{Q}}^{10}\rightarrow U$ dual to the divisor map $\Bbb{Q}[U]^{\ast}/\Bbb{Q}^{\ast} \rightarrow \text{Div}_{\Xi \setminus U}(\Xi)$, where the latter denotes  the free group of divisors in $\Xi$ with support in $\Xi \setminus U$. A set of generators of
$\Bbb{Q}[U]^{\ast}/\Bbb{Q}^{\ast}$ is given by $x_1/y_3$, $x_2/y_3$, $x_3/y_3$, $y_1/y_3$ and $y_2/y_3$. Computing the divisors of this set,   we can determine  the restriction of $\phi$ to the open subset of $\Omega$ where all coordinates are different from zero. We conclude that it  is given by \eqref{418}  on this Zariski dense
subset and hence everywhere on $\Omega$.\\

Having completed the proof of Lemma \ref{thm6}, we proceed to relate this result to the universal torsor over $X$.  By the Leray spectral sequence $H_{\text{et}}^p(\bar{B}, R^q \bar{\lambda}_{\ast}\Bbb{G}_m) \rightarrow H^{p+q}_{\text{et}}(\bar{X}, \Bbb{G}_m)$ (see \cite[(4.5)]{Gr1}) applied to the corresponding morphisms $\bar{\lambda} : \bar{X} \rightarrow \bar{B}$, $\bar{\gamma} : \bar{\Xi} \rightarrow \bar{B}$ over $\overline{\Bbb{Q}}$, there is a  commutative diagram  with exact rows  of trivial $\mathfrak{g}$-modules
% \begin{equation}\label{413}
$$ \begin{CD}
0 @>>> \text{Pic}(\bar{B}) @>{\bar{\gamma}^{\ast}}>> \text{Pic}(\bar{\Xi}) @>>> \Bbb{Z} @>>> 0\\
@.  @V{\text{id}}VV       @VVV @VV{\text{id}}V\\
 0 @>>> \text{Pic}(\bar{B}) @>{\bar{\lambda}^{\ast}}>> \text{Pic}(\bar{X}) @>>> \Bbb{Z} @>>> 0
\end{CD}$$
%\end{equation}
  where $\bar{\gamma}^{\ast}$ and $\bar{\lambda}^{\ast}$ are the contravariant functorial maps from $H^1_{\text{et}}(\bar{B}, \Bbb{G}_m) = \text{Pic}(\bar{B})$. The restriction $\text{Pic}(\bar{\Xi}) \rightarrow \text{Pic}(\bar{X})$ is thus an isomorphism, and there is a dual sequence of $\Bbb{Q}$-tori
$$ 1 \rightarrow \Bbb{G}_m \rightarrow T \rightarrow S \rightarrow 1$$
where the character groups of $T$ and $S$ are given by $\widehat{T} = \text{Pic}(\bar{\Xi}) = \text{Pic}(\bar{X})$ and $\widehat{S}= \text{Pic}(\bar{B})$.   From the functoriality of \eqref{tor1} under $X \rightarrow \Xi$ we conclude that the universal $\Xi$-torsor under $T$ restricts to the universal $X$-torsor under $T$.  We now restrict the map $\phi$ defined in Lemma \ref{thm6} to the closed subset $\phi^{-1}(X) \subset \Omega$ defined by \eqref{res4}. If we apply \eqref{418}, then \eqref{res4} takes the form
\begin{equation}\label{420}
  \xi_1u_1w_2w_3 + \xi_2u_2w_1w_3 + \xi_3u_3w_1w_2 = 0.
\end{equation}
We may now define three regular functions $v_1, v_2, v_3$ on $\phi^{-1}(X)$ as follows. For notational simplicity we agree that all indices are understood modulo 3. On the principal open subset where $u_iw_{i+1}w_{i+2} \not= 0$, we let 
\begin{displaymath}
  v_i = - \frac{u_{i+1}w_{i+2}\xi_{i+1} + u_{i+2} w_{i+1}\xi_{i+2}}{u_iw_{i+1}w_{i+2}}, \quad v_{i+1} = \frac{\xi_{i+1}}{w_{i+1}}, \quad v_{i+2} = \frac{\xi_{i+2}}{w_{i+2}}.
\end{displaymath}
If in addition $w_i \not= 0$, then $v_i = \xi_i/w_i$ by \eqref{420}, such that $v_1, v_2, v_3$ are well-defined. If we write $u$ instead of $\xi_0$, we obtain the following main result of this section.

\begin{thm}\label{thm7} Let $O \subset  \Bbb{A}^{10}$ be the  subvariety with coordinates $(u, v_1, v_2, v_3, u_1, u_2, u_3, w_1, w_2, w_3)$ defined by \eqref{2.1} and 
%\begin{equation}\label{421}
 % u_1v_1+u_2v_2+u_3v_3 = 0,
%\end{equation}
\begin{equation}\label{422}
  u_iu_kw_jw_k \not= 0 \text{ for at least one triple }\{i, j, k\} = \{1, 2, 3\},
\end{equation}
\begin{equation}\label{423}
 (u, v_1, v_2, v_3) \not= (0, 0, 0, 0).
 \end{equation}
Let $\phi_O : O \rightarrow X$ be the morphism which sends $(u, v_1, v_2, v_3, u_1, u_2, u_3, w_1, w_2, w_3)$ to
\begin{equation}\label{424}
\begin{split}
(x_1, x_2, x_3, y_1, y_2, y_3) & = (v_1w_1, v_2w_2, v_3w_3, uu_2u_3w_1, uu_1u_3w_2, uu_1u_2w_3),\\
(Y_1, Y_2, Y_3; Z_1, Z_2, Z_3) & = (u_2u_3w_1, u_1u_3w_2, u_1u_2w_3; u_1w_2w_3, u_2w_1w_3, u_3w_1w_2).
\end{split}
\end{equation}
Then $\phi_O : O \rightarrow X$ is the underlying $X$-scheme of a universal torsor over $X$.
\end{thm}

Indeed, it follows easily from the definition of the $v_i$  that \eqref{420} and   \eqref{2.1} as well as the second condition in \eqref{416} and  \eqref{423} are equivalent if \eqref{422} holds.\\

Finally we turn our attention to integral points. Let $\underline{X} \subset \Bbb{P}^5_{\Bbb{Z}} \times \Bbb{P}^2_{\Bbb{Z}} \times \Bbb{P}^2_{\Bbb{Z}}$ be defined by \eqref{res4} -- \eqref{res6} and let $\underline{\Xi} \subset \Bbb{P}^5_{\Bbb{Z}} \times \Bbb{P}^2_{\Bbb{Z}} \times \Bbb{P}^2_{\Bbb{Z}}$ be defined by \eqref{res5} and \eqref{res6}. We may extend the Cox morphism $\phi : \Omega \rightarrow \Xi$ from Lemma \ref{thm6} to a morphism $\underline{\phi} : \underline{\Omega} \rightarrow \underline{\Xi}$ between toric schemes, since the Cox morphism is derived from a morphism of fans (see \cite[pp.\ 22-23]{Fu2}). By repeating the arguments in the proof of Lemma \ref{thm6} over $\Bbb{Z}$, one obtains an open subscheme $\underline{\Omega}$ of $\Bbb{A}_{\Bbb{Z}}^{10}$ with coordinates $(\xi_0, \xi_1, \xi_2, \xi_3, u_1, u_2, u_3, w_1, w_2, w_3)$ defined by \eqref{416}. The morphism  $\underline{\phi} : \underline{\Omega} \rightarrow \underline{\Xi}$ defined by \eqref{418} is  the underlying morphism of a   torsor $\underline{\phi}_{\underline{\mathcal{T}}} : \underline{\mathcal{T}} \rightarrow \underline{\Xi}$  under a split $\Bbb{Z}$-torus $\underline{T} \cong \Bbb{G}_{m, \Bbb{Z}}^5$ with $H^1_{\text{et}}(\Bbb{Z}, \underline{T}) = 1$ (cf.\ \cite[III.4.9]{Mi}). The $\Bbb{Z}$-torsor obtained by base extension of $\underline{\phi}_{\underline{\mathcal{T}}} :  \underline{\mathcal{T}} \rightarrow \underline{\Xi}$ to an integral point is therefore always trivial. Hence there is a bijection between $\underline{T}(\Bbb{Z})$-orbits of integral points on $\underline{\Omega}$ and integral points on $\underline{\Xi}$. 

If we restrict $\underline{\phi}$ to the closed subset $\underline{\phi}^{-1}(\underline{X})$ of $\underline{\Omega}$ defined by \eqref{424}, we may again introduce new coordinates such that $\underline{O} = \underline{\phi}^{-1}(\underline{X})$ is the (locally closed) subscheme of $\Bbb{A}_{\Bbb{Z}}^{10}$ defined by \eqref{2.1}, \eqref{422}, \eqref{423}, and $\underline{\phi}_{\underline{O}} : \underline{O} \rightarrow \underline{X}$ is given by \eqref{424}.

We are now ready to state and prove the following equivalent version of Lemma \ref{lem1}.

\begin{lem}\label{lemtorsor} Let $\mathcal{A}_0$ denote the set of $10$-tuples $(u, v_1, v_2, v_3, u_1, u_2, u_3, w_1, w_2, w_3)$ with $v_j \in \Bbb{Z}$, $u, u_j \in \Bbb{N}$ and $w_j \in \Bbb{Z} \setminus \{0\}$ satisfying \eqref{2.1} as well as the coprimality conditions
\begin{equation}\label{427}
(u_1u_2w_1w_3; u_1u_2w_2w_3; u_1u_3w_1w_2; u_1u_3w_2w_3; u_2u_3w_1w_2; u_2u_3w_1w_3) = 1,
\end{equation}
\begin{equation}\label{428}
(u; v_1w_1; v_2w_2; v_3w_3) = 1.
\end{equation}
Then the map $\mathcal{A}_0 \rightarrow \Bbb{Z}^6$ defined by \eqref{2.3} and \eqref{2.6} gives a bijection between $\mathcal{A}_0$ and the set of \emph{primitive} integral solutions to \eqref{1} and \eqref{nonzero}.
\end{lem}

Note that \eqref{427} is equivalent to  \eqref{2.2} and \eqref{2.5}, in which case $(u_2u_3w_1; u_1u_3w_2; u_1u_2w_3) = 1$, so that \eqref{428} is equivalent to $(x_1; x_2; x_3; y_1;y_2;y_3) = 1$.\\

\emph{Proof.} Let $\underline{W} \subset \Bbb{P}^5_{\Bbb{Z}}$ be the subscheme defined by \eqref{1} and $\underline{f} : \underline{X} \rightarrow \underline{W}$ be the extension of the resolution $f : X \rightarrow W$ induced by the projection $ \Bbb{P}^5_{\Bbb{Z}} \times  \Bbb{P}^2_{\Bbb{Z}} \times  \Bbb{P}^2_{\Bbb{Z}} \rightarrow  \Bbb{P}^5_{\Bbb{Z}}$. Then there are natural bijections $\underline{X}(\Bbb{Z}) = X(\Bbb{Q})$, $\underline{W}(\Bbb{Z}) = W(\Bbb{Q})$ and $X^{\circ}(\Bbb{Q}) = W^{\circ}(\Bbb{Q})$ for the open subsets $X^{\circ} \subset X$, $W^{\circ} \subset W$ where $y= y_1 y_2 y_3 \not= 0$. If we let $\underline{X}(\Bbb{Z})^{\circ} \subset \underline{X}(\Bbb{Z})$ correspond to $X^{\circ}(\Bbb{Q}) \subset X(\Bbb{Q})$ and $\underline{W}(\Bbb{Z})^{\circ} \subset \underline{W}(\Bbb{Z})$ correspond to $W^{\circ}(\Bbb{Q}) \subset W(\Bbb{Q})$, then we get a bijection $\underline{X}(\Bbb{Z})^{\circ} = \underline{W}(\Bbb{Z})^{\circ}$. Next, let $O^{\circ}\subset O$ be the open subset where $y_1y_2y_3 = u^3(u_1u_2u_3)^2 w_1w_2w_3 \not= 0$ and let $\underline{O}(\Bbb{Z})^{\circ} \subset \underline{O}(\Bbb{Z})$ correspond to $O^{\circ}(\Bbb{Q}) \subset O(\Bbb{Q})$ under the bijection $\underline{O}(\Bbb{Z}) = O(\Bbb{Q})$. As $\phi^{-1}(W^{\circ}) = O^{\circ}$, we obtain a bijection between the $\underline{T}(\Bbb{Z})$-orbits in $\underline{O}(\Bbb{Z})^{\circ}$ and the points in $\underline{X}(\Bbb{Z})^{\circ}$. We observe that an integral $10$-tuple $(u, v_1, v_2, v_3, u_1, u_2, u_3, w_1, w_2, w_3)$ belongs to   $\underline{O}(\Bbb{Z})$ if and only if \eqref{2.1}, \eqref{422} and \eqref{423}  hold for all reductions modulo $p$. Hence it is in $\underline{O}(\Bbb{Z})$ if and only if \eqref{2.1}, \eqref{427} and \eqref{428} hold.  There are $2^{\dim T} = 32$ integral points in each $\underline{T}(\Bbb{Z})$-orbit in $\underline{O}(\Bbb{Z})$ with coordinates only differing by signs. For orbits in $\underline{O}(\Bbb{Z})^{\circ}$, the four $u$-coordinates do not vanish; there are exactly two integral points in each such $\underline{T}(\Bbb{Z})$-orbit with $u > 0$ and all $u_j > 0$,  and these two points have the same $v_j$-coordinates and opposite non-zero $w_j$-coordinates.  Summarizing the above discussion, we have shown that there is a bijection between the set of such pairs in $\mathcal{A}_0$ and $\underline{X}(\Bbb{Z})^{\circ} = \underline{W}(\Bbb{Z})^{\circ}$ where the latter set may be identified with the pairs $\pm(x_1, x_2, x_3, y_1, y_2, y_3)$ of primitive sextuples of integers satisfying \eqref{1} and \eqref{nonzero}.  The map from $\underline{O}(\Bbb{Z})^{\circ}$ to $\underline{W}(\Bbb{Z})^{\circ}$ comes from $\underline{f} \circ \underline{\phi}$ and is thus given by \eqref{2.3} and \eqref{2.6}. As the signs of $(x_1, x_2, x_3, y_1, y_2, y_3) $ are opposite for the two tuples in $\mathcal{A}_0$ in the same $\underline{T}(\Bbb{Z})$-orbit, we have established the desired bijection.

\section{Peyre's conjecture}\label{maninpeyre}

The aim of this chapter is to formulate Peyre's conjecture on the asymptotic
behaviour of $N(P)$. As the fourfold $W$ $\subset \Bbb{P}^{5}
$ defined by \eqref{1} is singular, we cannot refer to the original conjectures
of Manin \cite{FMT} and Peyre \cite{P} for Fano varieties. But as $%
f:X\rightarrow W$ restricts to an isomorphism from $X^{\circ }$ to $W^{\circ
}$, we obtain that
\begin{equation*}\label{51} 
  N(P)= |\left\{ x\in X^{\circ }(\Bbb{Q}%
):(H\circ f)(x)\leq P\right\}|
\end{equation*}
where the height function $H\circ f:X(\Bbb{Q}) \rightarrow \Bbb{N}$
is anticanonical, as $f$ is crepant. Since $X$ is an ``almost
Fano" variety (see Lemma \ref{almostfano}), we may refer to the ``formule empirique" in \cite[5.1]{P2}   for anticanonical counting functions on such varieties. This
formula predicts that 
\begin{equation}\label{52}
 N(P)\sim \Theta _{H}(X)P(\log P)^{\text{rk Pic}(X)-1}
\end{equation}
where $\Theta _{H}(X)=\alpha (X)\tau _{H}(X)$ for a suitable adelic Tamagawa
volume $\tau _{H}(X)$ of $X(\mathbf{A})$ and $\alpha(X)$ as in \eqref{defalpha}. Peyre demands for his
formula to hold that the cohomological Brauer group $Br'(\overline{X})=H_{\text{et}}^{2}(\overline{X},\Bbb{G}_{m})$ vanishes which is 
 true for our rational fourfold since $Br'(\overline{X})$ is a birational invariant \cite[Thm 7.1]{Gr1}   and $Br'
(\Bbb{P}_{\overline{\Bbb{Q}}}^n)=0$. There is also one hypothesis on 
$C_{\text{eff}}(X)$ in \cite[3.3]{P2}  which is satisfied thanks to Lemma
\ref{lemx}(ii). Finally, Peyre  assumes that there are no weakly accumulating subsets
(see \cite[(3.1)]{P2}) on $X^{\circ }$ for \eqref{52} to be valid. This will  
follow from our asymptotic formula, but is expected here because of the
group structure on $G=X^{\circ }$. 

The main goal of this chapter  is to define and  compute $\tau _{H}(X)$ for our
fourfold $X$, which together with the previous results on $\text{Pic}(X)$ in \eqref{picard} and $%
\alpha (X)$ in Lemma \ref{alpha}  gives an explicit conjecture for the asymptotic formula of $N(P)$.

For $i = 4, 5, 6$ it will be convenient to set $x_{i}=y_{i-3}$  and $F(x_{1},\ldots ,x_{6})=x_{1}y_{2}y_{3}+x_{2}y_{1}y_{3}+x_{3}y_{1}y_{2}$. For $1 \leq i \leq 6$ we write $\Bbb{P}_{(i)}^5 \subset \Bbb{P}^5$, $\Xi_{(i)} \subset \Xi$ (with $\Xi$ as in the previous chapter), $W_{(i)} \subset W$ and $X_{(i)} \subset X$ for the principal open subsets where $x_i \not= 0$, and we introduce 
 the
affine coordinates $x_{j}^{(i)}=x_{j}/x_{i}$, $j\neq i$ for $\Bbb{P}^5_{(i)} = \Bbb{A}^5$ and $\Xi_{(i)} \subset \Bbb{P}^5_{(i)} \times \Bbb{P}^2 \times \Bbb{P}^2$. Then $W_{(i)}$ is the affine hypersurface in $\Bbb{P}^5_{(i)}$ defined by   $F_{i}(x_{1}^{(i)},\ldots ,\widehat{%
x_{i}^{(i)}},\ldots ,x_{6}^{(i)})=F(x_{1}^{(i)},\ldots,1,\ldots ,x_{6}^{(i)})$.

To define $\tau _{H}(X)$ we need another description of the
height function $H\circ f:X(\Bbb{Q})\rightarrow \Bbb{N}$ in terms of
an adelic metric on the anticanonical sheaf $\omega _{X}^{-1}$. This adelic
metric will be defined by means of  global sections of $\omega
_{X}^{-1}=f^{\ast }(\omega _{W}^{-1})$, which are inverse images of global  sections
on $\omega _{W}^{-1}$. If $s\in \Gamma
(U,L)$ is a local section of an $O_{W}$-module $L$, we shall write $f^{\ast }(s)$ for the local section $%
f^{-1}(s)\otimes _{f^{-1}O_{W}}1\in \Gamma (f^{-1}(U),f^{\ast }(L))$ of   $f^{\ast }(L)= f^{-1}(L)\otimes _{f^{-1}O_{W}}O_{X}$.

The global sections of $\omega _{W}^{-1}$ and $\omega_X^{-1}$ that we shall  use are dual to certain 4-forms on $W$ and $X$. These 4-forms are  given by Poincar\'{e} residues of rational 5-forms on $\Bbb{P}^5$ and $\Xi$. To control the rational 5-forms on $\Bbb{P}^5$ and $\Xi$, we need the following   lemma from the theory of toric varieties \cite[p.\ 86]{Fu2}.

\begin{lem}\label{lemtoric} Let  $V$ be a non-singular $n$-dimensional 
toric   variety with torus U and $\omega _{V}(\sum
_{k=1}^{r}D_{k})$ be the sheaf of $n$-forms on $V$ with at %
most simple poles along all irreducible components $D_{1},\ldots, %
D_{r}$ of $\delta V=V \setminus U$. Then there is a global section %
$s_{V}\in \Gamma (V,\omega _{V}(\sum _{k=1}^{r}D_{k}))$ such that %
$s_V=\pm \frac{d\chi _{1}}{\chi _{1}}\wedge \ldots \wedge \frac{d\chi _{n}}{\chi
_{n}}$ on $U$ for any set of $n$ characters $\chi
_{i}:U\rightarrow \Bbb{G}_{m}$ that form a basis of $M= {\rm Hom}%
(U, \Bbb{G}_{m})$. The section $s_V$ generates the $O_V$-module $\omega_V(\sum_{k=1}^r D_k)$. 
\end{lem}

We now apply this lemma to the torus $U$ given by the cokernel of the
diagonal inclusion of $\Bbb{G}_{m}$ in $\Bbb{G}_{m}^{6}$, and to the
toric fivefolds $V=\Bbb{P}^{5}$ and $V=\Xi $. If we let $p_{1}: \Xi
\rightarrow \Bbb{P}^{5}$ be the restriction of  $\Bbb{P}%
^{5}\times \Bbb{P}^{2}\times \Bbb{P}^{2}\rightarrow \Bbb{P}^{5}$ to 
$\Xi $, then $p_{1}$ gives an isomorphism between the open subsets $U_{\Xi
}\subset \Xi $ and $U_{\Bbb{P}^{5}}\subset \Bbb{P}^{5}$, where $%
x_{i}\neq 0$ for all $1\leq i\leq 6$. If we make the obvious identifications of
these two open subsets with $U$, then the group law on $U$ extends to
actions of $U$ on $\Xi $ and $\Bbb{P}^{5}$ such that $p_{1}$ is $U$%
-equivariant.

In particular (see \cite[p.\ 41]{Re2}), if $V=\Bbb{P}^{5}$ and $H_{k}$, $1\leq
k\leq 6$, are the coordinate planes of $\Bbb{P}^{5}$, then there is a
global nowhere vanishing section $s_{\Bbb{P}^{5}}$ of $\omega _{\Bbb{P}%
^{5}}(\Sigma _{k=1}^{6}H_{k})$, such that the restriction $s^{(i)}$ of $s_{%
\Bbb{P}^{5}}$ to $\Bbb{P}_{(i)}^{5}$ is equal to
\begin{equation*}\label{53}
 s^{(i)}=(-1)^{i}\frac{dx_{1}^{(i)}}{x_{1}^{(i)}}\wedge
\ldots \wedge \widehat{\frac{dx_{i}^{(i)}}{x_{i}^{(i)}}}\wedge \ldots \wedge \frac{%
dx_{6}^{(i)}}{x_{6}^{(i)}}\in \Gamma \Bigl(\Bbb{P}_{(i)}^{5}, \omega _{%
\Bbb{P}^{5}}\bigl(\sum _{k=1}^{6}H_{k}\bigr)\Bigr).
\end{equation*}
On the other hand, if $V=\Xi $, then we conclude from  the proof of Lemma \ref{thm6} that  there are ten irreducible
components of $\Xi \setminus U$ corresponding to the ten edges of the fan
of $\Xi $ and to the ten coordinate hyperplane sections of the universal
torsor $\Omega \subset \Bbb{A}^{10}$. We let $D(\xi _{i})$, $0\leq i\leq 3
$, be the image in $\Xi $ of the subset of $\Omega $ defined by $\xi _{i}=0$, 
and we let $D(u_{j})$ and $D(w_{j}),1\leq j\leq 3$, be the prime divisors on $\Xi $
defined in the same way. Then by Lemma \ref{lemtoric} there is a global section $%
s_{\Xi }$ on $\omega _{\Xi }(E)$ for $E=\sum _{j=0}^{3}D(\xi _{j})+\sum
_{j=1}^{3}(D(u_{j})+D(w_{j}))$ such that the restriction of $s_{\Xi }$ to
the open subset $\Xi _{(i)}$ of $\Xi $ where $x_{i}\neq 0$ is given
by
\begin{equation*}\label{54}
  s^{(i)}=(-1)^{i}\frac{dx_{1}^{(i)}}{x_{1}^{(i)}}\wedge
\ldots \wedge \widehat{\frac{dx_{i}^{(i)}}{x_{i}^{(i)}}}\wedge \ldots \wedge \frac{%
dx_{6}^{(i)}}{x_{6}^{(i)}}\in \Gamma (\Xi _{(i)}, \omega _{\Xi }(E)).\end{equation*}
Now let
\begin{displaymath}
\begin{split}
& \omega _{i}=\frac{x_{1}x_{2}x_{3}x_{4}x_{5}x_{6}%
}{x_{i}^{3}F}s_{\Bbb{P}^{5}}\in \Gamma \bigl(\Bbb{P}^{5},\omega _{%
\Bbb{P}^{5}}(W+3H_{i})\bigr),\\
& \varpi _{i}=\frac{x_{1}x_{2}x_{3}x_{4}x_{5}x_{6}%
}{x_{i}^{3}F}s_{\Xi }\in \Gamma\bigl (\Xi, \omega _{\Xi
}(X+3p_{1}^{\ast }H_{i})\bigr).
\end{split}
\end{displaymath}
Then on the open subsets where $x_{i}\neq 0$, we have for $1\leq i\leq 6$
that 
\begin{equation}\label{55}
\begin{split}
& \omega _{i}=\frac{(-1)^{i}}{F_{i}}dx_{1}^{(i)}\wedge \ldots \wedge 
\widehat{dx_{i}^{(i)}}\wedge \ldots\wedge dx_{5}^{(i)}\in \Gamma \bigl(\Bbb{P}%
_{(i)}^{5},\omega _{\Bbb{P}^{5}}(W)\bigr),\\
 & \varpi _{i}=\frac{(-1)^{i}}{F_{i}}dx_{1}^{(i)}\wedge
\ldots \wedge \widehat{dx_{i}^{(i)}}\wedge \ldots \wedge dx_{5}^{(i)}\in \Gamma\bigl (%
\Xi_{(i)},\omega _{\Bbb{P}^{5}}(X)\bigr).
\end{split}
\end{equation}

We now consider Poincar\'e residues of these forms. The Poincar\'e residue map 
is usually given as a homomorphism $\Omega _{V}^{n}(W)\rightarrow i_{\ast }\Omega
_{_{W}}^{n-1}$ for the inclusion map $i$ of a non-singular hypersurface $%
W\subset V$ in an $n$-dimensional non-singular variety (cf.\ \cite[p.\ 89]{Re2}, for
example). More generally,  one can also use Poincar\'{e} residues to define local
sections on the canonical sheaf $\omega _{W}$ of an arbitrary normal
hypersurface (cf.\ \cite{We}) as one still gets regular $(n-1)$-forms on the
non-singular locus $W_{ns}$ of $W$ and since $\omega _{W}=j_{\ast }\Omega
_{_{W_{ns}}}^{n-1}$ for the open embedding $j:$ $W_{ns}\rightarrow W$. For  our singular hypersurface $i:W\rightarrow 
\Bbb{P}^{5}$ we  therefore have a unique homomorphism $\text{Res} :\omega _{\Bbb{P}%
^{5}}(W)\rightarrow i_{\ast }\omega _{W}$ of $O_{\Bbb{P}^{5}}$-modules,
which sends 
 $\omega _{i} \in \Gamma (\Bbb{P}%
_{(i)}^{5}, \omega_{\Bbb{P}^5}(W))$  to the section $\text{Res}(\omega _{i})\in \Gamma (\Bbb{P}_{(i)}^{5},i_{\ast
}\omega _{W})=\Gamma (W_{(i)},\omega _{W})$ which at the open subset of $%
W_{(i)}$ where $\partial F_{i}/\partial x_{j}^{(i)}\neq 0$ is given
by
\begin{equation}\label{56}
%\begin{split}
 \text{Res}(\omega _{i})= \begin{cases} \frac{(-1)^{i+k}}{\partial F_{i}/\partial
x_{k}^{(i)}}dx_{1}^{(i)}\wedge \ldots \wedge \widehat{dx_{i}^{(i)}}\wedge
\ldots \wedge \widehat{dx_{k}^{(i)}}\wedge \ldots \wedge dx_{6}^{(i)}&  \text{if} \quad i<k,\\
% \text{Res}(\omega _{i})=
\frac{(-1)^{i+k-1}}{\partial
F_{i}/\partial x_{k}^{(i)}}dx_{1}^{(i)}\wedge \ldots \wedge \widehat{dx_{k}^{(i)}%
}\wedge \ldots \wedge \widehat{dx_{i}^{(i)}}\wedge \ldots \wedge dx_{6}^{(i)} & \text{if} \quad  k<i.\end{cases}
%end{split}
\end{equation}
Similarly, for the inclusion $\iota :X\subset \Xi$, we note that $%
X_{(i)}\subset \Xi_{(i)}$ is defined by $F_{i}$ on the open subset
of $\Xi_{(i)}$ where two of $y_{1},$ $y_{2},$ $y_{3}$ are
different from zero. As $\varpi _{i}\in \Gamma (\Xi%
_{(i)}^{5},\omega _{\Xi}(X))$ is given by the same $5$-form as in
\eqref{55}, we obtain similarly that
\begin{equation}\label{57}
%\begin{split}
  \text{Res}(\varpi _{i})= \begin{cases}\frac{(-1)^{i+k}}{\partial F_{i}/\partial
x_{k}^{(i)}}dx_{1}^{(i)}\wedge \ldots \wedge \widehat{dx_{i}^{(i)}}\wedge
\ldots \wedge \widehat{dx_{k}^{(i)}}\wedge \ldots \wedge dx_{6}^{(i)} & \text{if} \quad 
 i<k,\\
 %\text{Res}(\varpi _{i})=
 \frac{(-1)^{i+k-1}}{\partial
F_{i}/\partial x_{k}^{(i)}}dx_{1}^{(i)}\wedge \ldots \wedge \widehat{dx_{k}^{(i)}%
}\wedge \ldots \wedge \widehat{dx_{i}^{(i)}}\wedge \ldots \wedge dx_{6}^{(i)} & \text{if} 
\quad k<i\end{cases}
%\end{split}
\end{equation}
on the open subset of $X_{(i)}$ where $(y_{1}y_{2},y_{1}y_{3},y_{2}y_{3})%
\neq (0,0,0)$. 

\begin{lem}\label{extend} {\rm (i)}  The section  $\text{\rm Res}(\omega _{i})$ extends uniquely to a %
global nowhere vanishing section of $\omega _{W}(3(H_{i}\cap W))$.\\
{\rm (ii)} The section $\text{\rm Res}(\varpi _{i})$ extends uniquely to a global nowhere %
vanishing section of $\omega _{X}(3f^{\ast }(H_{i}\cap W))$.\\
{\rm (iii)} The section $f^{\ast }\left( \text{\rm Res}(\omega _{i})\right) \in \Gamma (X,f^{\ast }\omega
_{W}(3(H_{i}\cap W)))$ is sent to $\text{\rm Res}(\varpi _{i})$ under the %
natural homomorphism from $f^{\ast }\omega _{W}(3(W\cap H_{i}))$ to $%
\omega _{X}(3f^{\ast }(H_{i}\cap W))$.
\end{lem}

\emph{Proof.} To prove (i) and (ii), we make use of the adjunction
formula (see \cite[pp.\ 146-147]{GH}). In this way  we obtain isomorphisms $i^{\ast
}\omega _{\Bbb{P}^{5}}(W)\rightarrow \omega _{W}$ and $\iota ^{\ast
}\omega _{\Xi }(X)\rightarrow \omega _{X}$ adjoint to the Poincar\'{e}
residue maps. These maps induce in  turn isomorphisms $i^{\ast }\omega _{%
\Bbb{P}^{5}}(W+3H_{i})\rightarrow \omega _{W}(3(H_{i}\cap W))$ and $\iota
^{\ast }\omega _{\Xi }(X+3p_{1}^{\ast }H_{i})\rightarrow \omega
_{X}(3f^{\ast }(H_{i}\cap W))$, which send $i^{\ast }\omega _{i}$ to Res$%
(\omega _{i})$ and $\iota ^{\ast }\varpi _{i}$ to Res$(\varpi _{i})$. Hence
it suffices for the proof of (i) and (ii) to show that $\omega _{i}=\frac{%
x_{1}x_{2}x_{3}x_{4}x_{5}x_{6}}{x_{i}^{3}F}s_{\Bbb{P}^{5}}$ generates the 
$O_{\Bbb{P}^{5}}$-module $\omega _{\Bbb{P}^{5}}(W+3H_{i})$ and that $%
\varpi _{i}=\frac{x_{1}x_{2}x_{3}x_{4}x_{5}x_{6}}{x_{i}^{3}F}s_{\Xi }$
generates the $O_{\Xi }$-module $\omega _{\Xi }(X+3p_{1}^{\ast }H_{i})$. This follows from the last assertion of Lemma \ref{lemtoric}. Part (iii) follows from \eqref{56} and \eqref{57}. \\

Thanks to Lemma \ref{extend}, we may now define global nowhere vanishing sections $%
\tau _{i}= \text{Res}(\omega _{i})^{-1}$ of $\omega _{W}^{-1}(-3(H_{i}\cap W))$
and $\sigma _{i}= \text{Res}(\varpi _{i})^{-1}$of $\omega _{X}^{-1}(-3f^{\ast
}(H_{i}\cap W))$. We will regard them as anticanonical global sections and
use the following result to define $v$-adic norms and measures.

\begin{lem}\label{lemy} {\rm (i)} The section $\tau _{i}\in \Gamma (W,\omega _{W}^{-1})$ does not %
vanish anywhere on $W_{(i)}$.\\
{\rm (ii)} The section $\sigma _{i}\in \Gamma (X,\omega _{X}^{-1})$ does not vanish anywhere %
on $X_{(i)}$.\\
{\rm (iii)} The section $f^{\ast }\tau _{i}$ is mapped to $\sigma _{i}$ under the %
canonical isomorphism from $f^{\ast }\omega _{W}^{-1}$ to $\omega
_{X}^{-1}$.
\end{lem}

\emph{Proof.} This is an immediate consequence of the previous lemma
since $\omega _{W}^{-1}(-3(H_{i}\cap W))=\omega _{W}^{-1}$ on $W_{(i)}$ and $%
\omega _{X}^{-1}(-3f^{\ast }(H_{i}\cap W))=\omega _{X}^{-1}$ on $X_{(i)}$.\\

 %We observe that $\omega _{\Bbb{P}^{5}}(W)$ is generated by $\omega
%_{i}$ on $\Bbb{A}_{(i)}^{5}$ and that $\omega _{W}$ is generated by $\text{Res}%
%(\omega _i)$ on $W_{(i)}$ and that the $O_{W}$-module homomorphism $%
%i^{\ast }\omega _{\Bbb{P}^{5}}(W)\rightarrow \omega _{W}$ adjoint to $\text{Res} : 
%\omega _{\Bbb{P}^{5}}(W)\rightarrow i_{\ast }\omega _{W}$ is an
%isomorphism. We recall the 

%As $\omega _{i}=\frac{x_{1}x_{2}x_{3}x_{4}x_{5}x_{6}}{x_{i}^{3}F}s^{(i)}$
%outside $W\cup H_{i}$, we may thus extend $\omega _{i}$ to a nowhere
%vanishing global section in $\Gamma \bigl(\Bbb{P}^{5},\omega _{\Bbb{P}%
%^{5}}(W+3H_{i})\bigr)$, and $\text{Res}(\omega _{i})$ to a nowhere vanishing global
%section in $\Gamma (W,\omega _{W}(3(W\cap H_{i}))$. Now write $\tau _{i}\in
%\Gamma (W,\omega _{W}^{-1}(-3(W\cap H_{i}))\subset \Gamma (W,\omega
%_{W}^{-1})$ for the inverse of this extension of $\text{Res}(\omega _{i})$. Then,

In the following we shall use the standard absolute values $|.|_v : \Bbb{Q}_v \rightarrow [0, \infty)$ for the places $v$ of $\Bbb{Q}$ (including the archimedean place). As $\tau _{i}$ vanishes nowhere on $W_{(i)}$, we obtain for each place $v$   a $v$-adic norm on $%
\omega _{W}^{-1}$   by letting
\begin{equation*}\label{58}
  \left\Vert \tau (w_{v})\right\Vert _{v}=\min_j\Bigl|
\frac{\tau }{\tau _{j}}(w_{v})\Bigr|_{v}=\min_j\left\vert (\tau 
\text{Res}(\omega _{j}))(w_{v})\right\vert _{v}\end{equation*}
for a local section $\tau $ of $\omega _{W}^{-1}$ defined at $w_{v}\in W(%
\Bbb{Q}_{v})$ and where the minimum is taken over all $j\in \{1,\ldots ,6\}$ such that $\tau
_{j}(w_{v})\neq 0$.  
%Here $\left\vert 0\right\vert _{v}=0$ and $\left\vert
%p^{r}u\right\vert _{v}=p^{-r}$ for a unit $u$ in $\mathbf{Z}_{p}$ for the $p$%
%-adic place $v$, while $\left\vert ..\right\vert _{v}$ is the ordinary
%absolute valute for the real place.\bigskip
This definition is the same as in \cite[pp.\ 107-108]{P}, although it is called a $v$-adic metric there. For more on $v$-adic norms
on invertible sheaves, see also \cite[Chapter 1]{Sa1}. 

As $\sigma _{i}$ vanishes nowhere on $X_{(i)}$, we obtain in the same way a $%
v$-adic norm on $\omega _{X}^{-1}$ for each place $v$ of $\Bbb{Q}$ by
letting\
\begin{equation}\label{59}
 \left\Vert \sigma (x_{v})\right\Vert
_{v}=\min_j\Bigl\vert \frac{\sigma }{\sigma _{j}}(x_{v})\Bigr\vert
_{v}=\min_j \left\vert (\sigma \text{Res}(\varpi _{j}))(x_{v})\right\vert
_{v} 
\end{equation}
for a local section $\sigma $ of $\omega _{X}^{-1}$ defined at $x_{v}\in X(%
\Bbb{Q}_{v})$ and where again the minimum is taken over all $j\in \{1,\ldots,6\}$ with $\sigma
_{j}(x_{v})\neq 0$.

\begin{lem}  {\rm (i)} Let $w\in W(\Bbb{Q})$ and let $\tau $ be a local 
section of $\omega _{W}^{-1}$ with $\tau (w)\neq 0$. Then 
\begin{equation}\label{510}
 H(w)=\prod\limits_{\text{\rm all }v} 
\left\Vert \tau (w)\right\Vert _{v}^{-1}.
\end{equation} 
{\rm (ii)} Let $x\in X(\Bbb{Q})$ and let $\sigma $ be a local 
section of $\omega _{X}^{-1}$ with $\sigma (x)\neq 0$. Then 
\begin{equation}\label{511}
 H(f(x))=\prod\limits_{\text{\rm all }v} 
\left\Vert \sigma (x)\right\Vert _{v}^{-1}.
\end{equation} 
\end{lem}

\emph{Proof.} (i) As $\prod_v \vert \alpha \vert _{v}=1$ for $\alpha \in \Bbb{Q}%
^{\ast }$ it suffices to show \eqref{510} for one such local section $\tau $. So let $C\in 
\Bbb{Q}[x_{1},x_{2},x_{3},x_{4},x_{5},x_{6}]$ be a cubic form with $%
C(w)\neq 0$ and $\tau =\frac{C}{x_{j}^{3}}\tau _{j}$ for $j$ with $%
x_{j}(w)\neq 0$. Then
$$ \left\Vert
\tau (w)\right\Vert _{v}^{-1}= \max_{1\leq j\leq 6}\left\vert \frac{\tau
_{j}}{\tau }(w)\right\vert _{v}=\max_{1\leq j\leq 6}\Bigl|\frac{%
x_{j}^{3}}{C}(w)\Bigr|_{v}$$
which immediately gives the desired formula for $H(w)$.\\
To prove (ii), we use   the canonical isomorphism $f^{\ast }(\omega
_{W}^{-1}) = \omega _{X}^{-1}$ and choose $\sigma$ to be the image of $f^{\ast}(\tau)$ for some local section $\tau$ of $\omega_W^{-1}$ where $\tau(w) \not= 0$ for $w = f(x)$. It follows from Lemma \ref{lemy}(iii) that $\|\sigma(x) \|_v = \| \tau(w) \|_v$ for each $v$, so that \eqref{511} follows from \eqref{510}. \\

We now apply Peyre's definition \cite[(2.2.1)]{P} of a measure $\mu
_{v}$ on $X(\Bbb{Q}_{v})$ associated to a $v$-adic norm on $\omega
_{X}^{-1}$. Let $\left\vert \text{Res}(\varpi _{i})\right\vert _{v}$ be the $%
v$-adic density on $X_{(i)}(\Bbb{Q}_{v})$ of the volume form $\text{Res}(\varpi
_{i})$ on $X_{(i)}$. Then for our particular $v$-adic norm $\left\Vert
{.}\right\Vert _{v}$ defined in \eqref{59}, we get the measure where
\begin{equation}\label{512}
   \mu _{v}(N_{v})=\int_{N_{v}}\frac{\left\vert \text{Res}%
(\varpi _{i})\right\vert _{v}}{\max_{1\leq j\leq 6}\left\vert \sigma _{j}%
\text{Res}(\varpi _{i})\right\vert _{v}}=\int_{N_{v}}\frac{\left\vert \text{%
Res}(\varpi _{i})\right\vert _{v}}{\max_{1\leq j\leq 6}\left\vert
(x_{j}/x_{i})^{3}\right\vert _{v}}
\end{equation}
for a Borel subset $N_{v}$ of $X_{(i)}(\Bbb{Q}_{v})$. 

To get a more explicit description of $\mu _{v}$, let us write $%
t_{j}=x_{j}^{(6)}=x_{j}/x_{6}$. Then, by \eqref{55} and \eqref{57}, we have
that
$$ \omega _{6}=\frac{1}{F_{6}}dt_{1}\wedge
dt_{2}\wedge dt_{3}\wedge dt_{4}\wedge dt_{5}, \quad \text{ Res}(\omega _{6})=\frac{(-1)^{k-1}}{\partial
F_{6}/\partial t_{k}}dt_{1}\wedge \ldots\wedge \widehat{dt_{k}}\wedge
\ldots dt_{5} $$
for any $k=1,2,3,4,5$ and $%
F_{6}(t_{1},t_{2},t_{3},t_{4},t_{5})=t_{1}t_{5}+t_{2}t_{4}+t_{3}t_{4}t_{5}$.
For instance, choosing  $k=3$, we obtain 
\begin{equation}\label{513} \mu _{v}(N_{v})=\int_{N_{v}}\frac{dt_{1}dt_{2}dt_{4}dt_{5}}{%
\left\vert t_{4}t_{5}\right\vert _{v}\max (\left\vert t_{1}\right\vert
_{v}^{3},\left\vert t_{2}\right\vert _{v}^{3}, \vert \frac{t_{1}}{t_{4}}+%
\frac{t_{2}}{t_{5}} \vert _{v}^{3},\left\vert t_{4}\right\vert
_{v}^{3},\left\vert t_{5}\right\vert _{v}^{3},1)}
\end{equation}
for any Borel subset $N_{v}$ of $\bigcap_{3\leq i\leq 6}X_{(i)}(%
\Bbb{Q}_{v})$. Here and elsewhere we assume that the underlying Haar
measure on $\Bbb{Q}_{v}$ is the usual Lebesgue mesure if $\Bbb{Q}_{v}=%
\Bbb{R}$, and that it is normalized by $\int_{\Bbb{Z}_{p}}dx=1$
if $v$ is $p$-adic. %From now on we write $\mu_{\infty}$ for $v = \infty$ und
%For the real place we will in the sequel write $\mu
%_{\infty }$ instead of $\mu _{v}$ and for the $p$-adic place we will write $%
%\mu _{p}$ instead of $\mu _{v}$.

Now let $L_{p}(s, \text{Pic}(\overline{X}))= \det(1-p^{-s} \text{Fr}_{p}\mid \text{Pic}%
( X_{\overline{\Bbb{F}}_{p}} ) \otimes \Bbb{Q})^{-1}$ for
a prime $p$. As Pic $ ( X_{\overline{\Bbb{F}}_{p}} ) =%
\Bbb{Z}^{5}$ with trivial Galois action, we get that  
$$L(s, \text{Pic}(\overline{X}))=\prod\limits_{\text{all }p}L_{p}(s, \text{Pic}(%
\overline{X}))=\prod\limits_{\text{all }p}(1-p^{-s})^{-1}=\zeta
(s)^{5}$$
for $s\in \Bbb{C}$ with $\Re s >1$. In particular, $\lim_{s\rightarrow
1}(s-1)^{5}L(s, \text{Pic}(\overline{X}))=1$ and $L_{p}(1, \text{Pic}(\overline{X}))^{-1}=(%
\frac{p-1}{p})^{5}$. For our particular fourfold $%
X$, Peyre's Tamagawa measure $\mu _{H}$ on $X(\mathbf{A} )=X(\Bbb{R})\times \prod_p X(\Bbb{Q}%
_{p}) $ (see \cite[Def.\ 4.6]{P2}) is therefore  given by $\mu _{H}=\mu _{\infty }\times \prod_p(\frac{%
p-1}{p})^{5}\mu _{p}$, and it is shown in \cite{P2} that this gives a well-defined measure on $X(\mathbf{A})$. As $X(\Bbb{Q})$ is dense
in $X(\mathbf{A})$, we thus have  
\begin{equation}\label{514}
  \tau _{H}(X)= \mu _{H}(X(\mathbf{A} ))=\mu
_{\infty }(X(\Bbb{R}))\prod\limits_{\text{all }p}\left(\frac{p-1}{p}\right)^{5}\mu
_{p}(X(\Bbb{Q}_{p}))
\end{equation}
  by \cite[Def.\ 4.8]{P2}, and it remains to determine the local volumes $\mu_{\infty}(X(\Bbb{R}))$ and $\mu_p(X(\Bbb{Q}_p))$. \\
  
 To compute $\mu _{\infty }(X(\Bbb{R}))$, let $N_{\infty }(\Bbb{R}%
)=\bigcap_{3\leq i\leq 6}X_{(i)}(\Bbb{R})$. Then, $\mu
_{\infty }(X(\Bbb{R}))=\mu _{\infty }(N_{\infty }(\Bbb{R}))$ by Sard's
theorem. Hence, by \eqref{513}, we obtain  
\begin{equation*}\label{515}
 \mu _{\infty }(X(\Bbb{R}))= \int_{-\infty }^{\infty
}\int_{-\infty }^{\infty }\int_{-\infty }^{\infty }\int_{-\infty }^{\infty }%
\frac{dt_{1}dt_{2}dt_{4}dt_{5}}{\left\vert t_{4}t_{5}\right\vert \max
(\left\vert t_{1}\right\vert ^{3},\left\vert t_{2}\right\vert
^{3}, \vert \frac{t_{1}}{t_{4}}+\frac{t_{2}}{t_{5}} \vert
^{3},\left\vert t_{4}\right\vert ^{3},\left\vert t_{5}\right\vert ^{3},1)}.
\end{equation*}
It is a long, but elementary and straightforward calculation to check that 
     \begin{equation}\label{4int}
 \mu_{\infty}(X(\Bbb{R})) = 
    12(\pi ^{2} +24\log 2-3).
 \end{equation}    \\

We proceed to compute $\mu_p(X(\Bbb{Q}_p))$. Let $\Omega \subset \Bbb{A}^{10}$ and $\varphi :\Omega
\rightarrow \Xi $ be as in Lemma \ref{thm6}. Then $O =\varphi ^{-1}(X)$, and the
structure morphism of the $X$-torsor $O$ is given by the restriction $\varphi
_{O}:O\rightarrow X$ of $\phi $ to $O$ (see Theorem \ref{thm7}). Hence $O\subset \Omega $ is the hypersurface defined by $\Phi =\xi
_{1}u_{1}w_{2}w_{3}+\xi _{2}u_{2}w_{1}w_{3}+\xi _{3}u_{3}w_{1}w_{2}$, cf.\ \eqref{420}. By Lemma \ref{lemtoric}, there is a rational 10-form $s_{\Omega}$ on the toric variety  $%
\Omega \subset \Bbb{A}^{10}$ defined by
$$ s_{\Omega }=\frac{d\xi _{0}}{\xi _{0}}\wedge \frac{d\xi
_{1}}{\xi _{1}}\wedge \frac{d\xi _{2}}{\xi _{2}}\wedge \frac{d\xi _{3}}{\xi
_{3}}\wedge \frac{du_{1}}{u_{1}}\wedge \frac{du_{2}}{u_{2}}\wedge \frac{%
du_{3}}{u_{3}}\wedge \frac{dw_{1}}{w_{1}}\wedge \frac{dw_{2}}{w_{2}}\wedge 
\frac{dw_{3}}{w_{3}}.$$ 
For $1 \leq i \leq 6$ we let $\varpi _{i}^{\Omega }=\frac{x_{1}x_{2}x_{3}x_{4}x_{5}x_{6}}{%
x_{i}^{3}F}s_{\Omega }$. Finally, we define
\begin{equation}\label{516}
   \varpi ^{\Omega }=\frac{1}{\Phi }d\xi _{0}\wedge
d\xi _{1}\wedge d\xi _{2}\wedge d\xi _{3}\wedge du_{1}\wedge du_{2}\wedge
du_{3}\wedge dw_{1}\wedge dw_{2}\wedge dw_{3}.
\end{equation}
Then, from $F=\xi _{0}^{2}u_{1}u_{2}u_{3}\Phi $, we conclude that $\varpi
_{i}^{\Omega }=\varpi ^{\Omega }/x_{i}^{3}$ on $\Omega _{(i)}=\varphi
^{-1}(\Xi _{(i)})$ for $1\leq i\leq 6$, where $x_{i}$ and $x_{i+3}=y_{i}$, $%
1\leq i\leq 3$, now denote the \emph{affine} coordinates given by the expressions
in \eqref{418}. 

The following construction works for an arbitrary place $v$, although in the present situation we are only interested in non-archimedean places. Any $\varpi _{i}^{\Omega }\in \Gamma (\Omega _{(i)}, \omega _{\Omega
}(O))$ has a Poincar\'{e} residue $\text{Res}(\varpi _{i}^{\Omega })\in \Gamma
(O_{(i)},\omega _{O})$ on $O_{(i)}=\varphi ^{-1}(X_{(i)})$. We may now, just
as in \eqref{512}, use the six local volume forms $\text{Res}(\varpi _{i}^{\Omega })$ to
construct a $v$-adic measure $m_{v}$ on $O(\Bbb{Q}_{v})$ by letting
\begin{equation*}\label{517}
 m_{v}(M_{v})=\int_{M_{v}}\frac{\left\vert \text{Res}%
(\varpi _{i}^{\Omega })\right\vert _{v}}{\max_{1\leq j\leq 6}\left\vert
(x_{j}/x_{i})^{3}\right\vert _{v}}
\end{equation*}
for a Borel subset $M_{v}$ of $O_{(i)}(\Bbb{Q}_{v})$. %We shall also write 
%$m_{p}$ if $v$ is $p$-adic and $m_{\infty }$ if $v$ is real. 
The connection between $m_p$ and $\mu_p$ for a prime $p$ will become clear in Lemma \ref{lemmeasures} below. 
To start with,  we consider the relative canonical sheaves $\omega
_{\Omega /\Xi }$, $\omega _{O/X}$ and apply the following result.

\begin{lem}
Let $E^{ \Omega }=\varphi ^{\ast }E\in \text{\rm Div}%
(\Omega) $ be the sum of the ten prime divisors of $\Omega $ 
defined by the ten coordinate hyperplanes of $\Bbb{A}^{10}$. %
Then there is a unique global nowhere vanishing section %
$s_{\Omega /\Xi }\in \Gamma (\Omega , \omega _{\Omega /\Xi })$ such that 
$s_{\Omega }=s_{\Omega /\Xi }\otimes \varphi ^{\ast }s_{\Xi }$ under the  
natural isomorphism $\omega _{\Omega }(E^{_{\Omega }})=\omega _{\Omega
/\Xi }\otimes \varphi ^{\ast }\omega _{\Xi }(E).$ Moreover, if we let %
$\iota _{O}:O\rightarrow \Omega $ be the inclusion map, $%
s_{O/X}\in \Gamma (O, \omega _{O/X})$ be the image of $\iota
_{O}^{\ast }s_{\Omega /\Xi }\in \Gamma (O, \iota _{O}^{\ast }\omega
_{\Omega /\Xi })$ under the functorial isomorphism from $\iota
_{O}^{\ast }\omega _{\Omega /\Xi }$ to $\omega _{O/X}$ and $%
s_{O/X}^{(i)}$ be the restriction of $s_{O/X}$ to $%
O_{(i)}$, then
\begin{equation}\label{518}
\text{\em Res}(\omega _{i}^{\Omega })=s_{O/X}^{(i)}\otimes
\varphi _{O}^{\ast }\text{\em Res}(\varpi _{i})
\end{equation}
for $1\leq i\leq 6$ under the canonical isomorphism $\omega _{O} %
=\omega _{O/X}\otimes \varphi _{O}^{\ast }\omega _{X}$.
 \end{lem}

\emph{Proof.} The isomorphism between $\omega _{\Omega }(E^{\Omega
})$ and $\omega _{\Omega /\Xi }\otimes \varphi ^{\ast }\omega _{\Xi }(E)$
is induced by the canonical isomorphism between $\omega _{\Omega }$ and $%
\omega _{\Omega /\Xi }\otimes \varphi ^{\ast }\omega _{\Xi }$, and the first
statement is obvious as $s_{\Omega }$ (resp.\ $\varphi ^{\ast }s_{\Xi }$) is
a global generator of $\omega _{\Omega }(E^{\Omega })$ (resp. $\varphi
^{\ast }\omega _{\Xi }(E)$). The second statement follows from a
functoriality property of Poincar\'{e} residues, which says that there is a
natural commutative diagram of isomorphisms of $O_{O}$-modules
\begin{displaymath}
\begin{CD}
\iota _{O}^{\ast }\omega _{\Omega
}(E^{\Omega })  @>>>  \iota _{O}^{\ast }(\omega _{\Omega /\Xi
})\otimes \iota _{O}^{\ast }\varphi ^{\ast }\omega _{\Xi }(E)\\
@VVV       @VVV\\
\omega_O @> >> \omega _{O/X}\otimes \varphi _{O}^{\ast }\omega
_{X}
\end{CD}
\end{displaymath}
 where the first vertical map is given by the adjoint Poincar\'{e} residue
map for $\iota _{O}:O$ $\rightarrow \Omega $ and the second vertical map
makes use of the isomorphism from $\iota _{O}^{\ast }\varphi ^{\ast }\omega
_{\Xi }(E)=\phi _{O}^{\ast }\iota ^{\ast }\omega _{\Xi }(E)$ to $\phi
_{O}^{\ast }\omega _{X}$ induced by the adjoint Poincar\'{e} residue map for 
$\iota :X$ $\rightarrow \Xi $.\\

We now apply this result to the $v$-adic analytic manifolds associated to $O$
and $X$ and refer to \cite[Ch.\ III]{Se}  for basic definitions and properties of
such manifolds and to \cite[Ch.\ 3]{Sa1} for the notion of torsors over $v$-adic
analytic manifolds. 

\begin{lem}\label{lemma17}  {\rm (i)}  The map $\phi _{O,v}: O(\Bbb{Q}%
_{v})\rightarrow X(\Bbb{Q}_{v})$ induced by $\phi _{O}$ is a %
submersion of $v$-adic analytic manifolds, which makes $O(%
\Bbb{Q}_{v})$  an analytic $X(\Bbb{Q}_{v})$-torsor under $%
T(\Bbb{Q}_{v})$.\\
{\rm (ii)} The relative volume form $s_{O/X}\in \Gamma (O, \omega
_{O/X})$ defines $v$-adic  measures on the fibres of $\phi
_{O,v}$ which yields a linear functional  $\Lambda
_{v}:C_{c}(O(\Bbb{Q}_{v}))\rightarrow C_{c}(X(\Bbb{Q}_{v}))$ when %
we  integrate along the fibres of $\phi _{O,v}$.\\
{\rm (iii)}  If $\beta _{v}\in C_{c}(O(\Bbb{Q}_{v}))$, then $\int_{O(%
\Bbb{Q}_{v})}\beta_v \, m_{v}=\int_{X(\Bbb{Q}_{v})}\Lambda _{v}(\beta
_{v})\mu _{v}$.
\end{lem}

 \emph{Proof.}  For (i), see \cite[pp.\ 126-127]{Sa1}. To obtain (ii) and (iii),
use \cite[Theorem 1.22]{Sa1} and \eqref{518}.\\

 We may now reinterpret the $p$-adic factor of Peyre's constant $\Theta
_{H}(X)$ as a $p$-adic density of the universal torsor over $X$.

\begin{lem}\label{lemmeasures}   Let $\underline{O}$ be the scheme defined %
in Section \ref{universaltorsor}. Then $ m_{p}(\underline{O}(\Bbb{Z}_{p}))=(\frac{p-1}{p%
})^{5}\mu _{p}(X(\Bbb{Q}_{p}))$ for any prime $p$. 
\end{lem}

\emph{Proof.}  We embed $\underline{O}(\Bbb{Z}_{p})$ as an open
subset of $\underline{O}(\Bbb{Q}_{p})= O(\Bbb{Q}_{p})$ and let $\chi
_{p}:O(\Bbb{Q}_{p}) \rightarrow \{0,1\}$ be the characteristic function
of $\underline{O}(\Bbb{Z}_{p})$. Then $\chi _{p}\in C_{c}(O(\Bbb{Q}%
_{p}))$ and $m_{p}(\underline{O}(\Bbb{Z}_{p}))=\int_{X(\Bbb{Q}%
_{p})}\Lambda _{p}(\chi _{p})\mu _{p}$ by the previous lemma. It is
therefore enough to show that $\Lambda _{p}(\chi _{p})\in C_{c}(X(\Bbb{Q}%
_{p}))$ has value $(\frac{p-1}{p})^{5}$ at all points of $\underline{X}(%
\Bbb{Z}_{p})=X(\Bbb{Q}_{p})$. But it is clear that the decomposition $%
s_{\Omega }=s_{\Omega /\Xi }\otimes \varphi ^{\ast }s_{\Xi }$ may be carried
out over $\Bbb{Z}$ such that $s_{\Omega /\Xi }$ extends to a $\underline{T%
}$-equivariant generator of $\omega _{\underline{\Omega }/\underline{\Xi }}$
and $s_{O/X}$ to a $\underline{T}$-equivariant generator $s_{\underline{O}/%
\underline{X}}$ of $\omega _{\underline{O}/\underline{X}}$. If $P$ is a $%
\Bbb{Z}_{p}$-point on $\underline{X}$ and $\underline{O}_{P}\rightarrow $ 
$P$ the base extension of $\underline{O}\rightarrow $ $\underline{X}$, then $%
s_{\underline{O}/\underline{X}}$ will therefore pull back to a $\underline{T}%
_{\Bbb{Z}_{p}}$-equivariant global section $s_{\underline{O}_{P}}$ on $%
\omega _{\underline{O}_{P}/\Bbb{Z}_{p}}$. As the torsor over $P$ is
trivial and $\underline{T}\cong\Bbb{G}_{m,\Bbb{Z}}^{5}$,
there are   affine coordinates $(t_{1},\ldots ,t_{5})$ for the affine $%
\Bbb{Z}_{p}$-scheme $\underline{O}_{P}$ such that $s_{\underline{O}_{P}}=%
\frac{dt_{1}}{t_{1}}\wedge \ldots \wedge \frac{dt_{5}}{t_{5}}$.
Hence $$\Lambda _{p}(\chi _{p})(P)=\int_{\underline{O}_{P}(\Bbb{Z}%
_{p})}\vert s_{\underline{O}_{P}}\vert =\prod\limits_{1\leq i\leq
5}\int_{\Bbb{Z}_{p}^{\ast }}\frac{dt_{i}}{t_i}=\left(\frac{p-1}{p}%
\right)^{5}, $$ and we are done.\\

To compute $m_{p}(\underline{O}(\Bbb{Z}_{p}))$, we give an alternative
definition of $m_{p}$. As $\omega ^{\Omega }\in  \Gamma (\Omega , \omega
_{\Omega }(O))$ by \eqref{516}, it has a residue form $\text{Res}(\omega ^{\Omega })\in 
 \Gamma (O,\omega _{O})$. If we again let $x_{i}$, $1\leq i\leq 6$, denote
the affine coordinates given by  the expressions in \eqref{418}, then $\omega
^{\Omega }$ restricts to $x_{i}^{3} \,\text{Res}(\omega _{i}^{\Omega })$ on $O_{(i)}$
such that
\begin{equation}\label{519}
 m_{p}(M_{p})=\int_{M_{p}}\frac{\vert \text{Res}%
(\omega ^{\Omega })\vert _{p}}{\max_{1\leq j\leq 6}\vert
x_{j}^{3}\vert _{p}}=\int_{M_{p}}\frac{d{\bm \xi} \, d\textbf{u} \, d\textbf{w} /d\Phi }{\max_{1\leq j\leq
6}\vert x_{j}^{3}\vert _{p}}
\end{equation}
for Borel subsets $M_{p}$ of $O(\Bbb{Q}_{p})$ and for the $p$-adic
density $\vert \text{Res}(\omega ^{\Omega })\vert _{p}$ of $\text{Res}%
(\omega ^{\Omega })$, which we also denote by $d{\bm \xi} \, d\textbf{u} \, d\textbf{w}/d\Phi $. 
%d\xi _{0}d\xi _{1}d\xi
%_{2}d\xi _{3}\, du_{1}du_{2}du_{3} \, dw_{1}dw_{2}dw_{3}/d\Phi $. 
Here we have written $d {\bm \xi} = d\xi _{0}d\xi _{1}d\xi
_{2}d\xi _{3}$, $d\textbf{u} =  du_{1}du_{2}du_{3}$ and $d\textbf{w} = dw_{1}dw_{2}dw_{3}$ for notational simplicity. 

\begin{lem}\label{lemmab} One has
$$ m_{p}(\underline{O}(\Bbb{Z}_{p}))= \frac{|\underline{%
O}(\Bbb{F}_{p})|}{p^{\dim O}} =\frac{(p-1)^{5}(p^{2}+p+1)(p^{2}+4p+1)%
}{p^{9}}.$$
\end{lem}

\emph{Proof.} For $P\in \underline{\Omega }(\Bbb{Z}_{p})$   we
may find some $j\in \{1,\ldots ,6\}$ such that $p$ does not divide $x_{j}(P)$
(cf.\ \eqref{418}), since  $\underline{\varphi }$ restricts to a morphism from $%
\underline{\Omega }_{\Bbb{F}_{p}}$ to $\underline{\Xi }_{\Bbb{F}%
_{p}}\subset \Bbb{P}_{\Bbb{F}_{p}}^{5}\times \Bbb{P}_{\Bbb{F}%
_{p}}^{2}\times \Bbb{P}_{\Bbb{F}_{p}}^{2}$ modulo $p$. We conclude $\max_{1\leq j\leq 6}\vert x_{j}^{3}\vert _{p}=1$ for $P\in 
\underline{O}(\Bbb{Z}_{p})\subset \underline{\Omega }(\Bbb{Z}_{p})$
and 
\begin{equation*}\label{520}
   m_{p}(\underline{O}(\Bbb{Z}_{p}))=\int_{\underline{O}(%
\Bbb{Z}_{p})}\frac{ d{\bm \xi} \, d\textbf{u} \, d\textbf{w}}{d\Phi }
\end{equation*}
by \eqref{519}.  It is easy to see that this measure coincides with the $p$-adic
model measure defined in \cite[2.9]{Sa1}. As $\underline{O}$ is smooth over $%
\Bbb{Z}$, we may thus apply \cite[Cor.\ 2.15]{Sa1}   and conclude that $%
m_{p}(\underline{O}(\Bbb{Z}_{p}))=|\underline{O}(\Bbb{F}_{p})|/p^{\dim
O}$. To determine $|\underline{O}(\Bbb{F}_{p})|$, we note that the $%
\underline{X}_{\Bbb{F}_{p}}$-torsor $\underline{O}_{\Bbb{F}_{p}}$
under $\underline{T}_{\Bbb{F}_{p}}$ is locally trivial, such that $|%
\underline{O}(\Bbb{F}_{p})|/p^{\dim O}=|\underline{T}(\Bbb{F}_{p})| \cdot |%
\underline{X}(\Bbb{F}_{p})|/p^{9}$. To finish we note that $|\underline{T}(%
\Bbb{F}_{p})|=(p-1)^{5}$ and $|\underline{X}(\Bbb{F}%
_{p})|=(p^{2}+p+1)(p^{2}+4p+1)$ as $\underline{X}_{\Bbb{F}_{p}}$ is a $%
\Bbb{P}^{2}$-bundle over a split del Pezzo $\Bbb{F}_{p}$-surface of
degree 6.\\

We remark on the side that   \cite[Cor.\ 2.15]{Sa1}   gives  $ m_{p}(\underline{O}(\Bbb{Z}_{p}))=|\underline{O}(\Bbb{Z}%
/p^{r})|/p^{r\dim O}$ more generally for all $r\geq 1$. We have now all the ingredients to give an explicit evaluation of Peyre's empirical formula for the counting function $N(P )$. A combination of \eqref{514}, \eqref{4int}, Lemma \ref{lemmeasures} and Lemma \ref{lemmab} yields

\begin{thm}\label{thmpey}
 Let $X \subset \Bbb{P} ^{5}\times 
\Bbb{P}^{2}\times \Bbb{P}^{2}$ be the %
fourfold defined by \eqref{res4} -- \eqref{res6} and let $\tau _{H}(X)$ be Peyre'%
s  adelic Tamagawa  volume of  $X(\mathbf{A})$ %
associated to the $v$-adic norms on $\omega _{X}^{-1}$ in \eqref{59}
. Then
\begin{equation*}\label{521}
  \tau _{H}(X)=12(\pi ^{2}+24\log 2-3)\prod\limits_{\text{\rm all }%
p}\left(1-\frac{1}{p}\right)^{5}\left( 1+\frac{5}{p}+\frac{6}{p^{2}}+\frac{5}{p^{3}}%
+\frac{1}{p^{4}}\right).  
\end{equation*}
\end{thm}

Combining this with \eqref{52}, Lemma \ref{alpha} and \eqref{picard}, we confirm that \eqref{pey*} agrees with Peyre's prediction.\\

\emph{Remark.} It is possible to interpret $\alpha (X)$ as a real
analogue of the $p$-adic convergence factor $|\underline{T}(\Bbb{F}%
_{p}|/p^{\dim T}$. Let $X^{\circ }(\Bbb{Q},P)=\left\{ x\in X^{\circ }(%
\mathbf{Q}):(H\circ f)(x)\leq P\right\} $. Then, if \eqref{55} holds, we have by
partial summation  and \eqref{defalpha} that
\begin{equation}\label{514a}
   \sum\limits_{x\in X^{\circ }(\mathbf{Q},P)}\frac{1}{%
(H\circ f)(x)}\sim \tau _{H}(X)\left( \frac{\alpha (X)}{\text{rk Pic}(X)}%
\right) (\log P)^{\text{rk Pic}(X)} = \tau_H(x)  \int_{\Delta (P)}ds
\end{equation}
%On the other hand, from \eqref{alpha} we deduce that $\frac{\alpha (X)}{\text{rk Pic}(X)}(\log P)^{\text{rk Pic}(X)}= \int_{\Delta (P)}ds$ 
where $\Delta (P)
$ is the set of all linear forms $\Lambda $ on $\text{Pic}(X)\otimes \Bbb{R}$ such
that $\Lambda ([-K_{X}])\leq \log P$ and $\Lambda \in C_{\text{eff}%
}(X)^{\vee }$. Now let $x\in X$ be the point where all six $x_{i}$%
-coordinates are equal. We may then identify $T$ with the fibre of the
torsor $\varphi _{O}:O\rightarrow X$ over $x$ such that the neutral element
of $T$ corresponds to the point in $\Bbb{A}^{10}$ with all coordinates
equal to $1$. Let $D(P)\subset T(\Bbb{R})$ be the subset
where $\min (u,u_{1},u_{2},u_{3},\left\vert w_{1}\right\vert ,\left\vert
w_{2}\right\vert ,\left\vert w_{3}\right\vert )\geq 1$ and where one and
hence all $\left\vert x_{i}^{3}\right\vert $ are at most $P$. Then, as $%
y_{1}=y_{2}=y_{3}$ on $T(\mathbf{R})$ all $w_{i}$   have the same sign.
If we let $D_{+}(P)\subset D(P)$ be the subset where all $w_{i}>0$ and $dt$
be the measure on $T(\Bbb{R})$ in Lemma \ref{lemma17} (see also Lemma \ref{lemmeasures}), we  obtain $\int_{D(P)}dt=2\int_{D_{+}(P)}ds$. Furthermore, by Lemma
\ref{lemx}(ii) we have that $\int_{D_{+}(P)}dt=\int_{\Delta (P)}ds$
such that  \begin{equation}\label{remark}\frac{\alpha (X)}{\text{rk Pic}(X)}(\log P)^{\text{rk Pic}(X)} 
= \frac{1}{2}\int_{D(P)}dt.\end{equation}
We may now give an heuristic derivation of the factor $\alpha (X(\Bbb{R}%
))\mu _{\infty }(X(\Bbb{R}))$ in Peyre's constant. Let $F(P)$ be the set
of all $\mathbf{r}=\left(
u,u_{1},u_{2},u_{3},v_{1},v_{2},v_{3},w_{1},w_{2},w_{3}\right) \in \Bbb{R}%
^{10}$ such that $$\min (u,u_{1},u_{2},u_{3},\left\vert w_{1}\right\vert
,\left\vert w_{2}\right\vert ,\left\vert w_{3}\right\vert) \geq 1,\quad  %
u_{1}v_{1}+u_{2}v_{2}+u_{3}v_{3}=0, \quad  \max \left\vert
x_{i}^{3}\right\vert \leq P$$ with $(x_{1},\ldots ,x_{6})$ as in \eqref{424}. Then,
by Lemma \ref{lemtorsor}  and \eqref{remark}, %the equality $\frac{1}{2}\int\limits_{D(P)}dt=\left( \frac{%
%\alpha (X)}{\text{rk Pic }X}\right) (\log P)^{\text{5}}$, 
we see that \eqref{514a} corresponds to the conjecture 
$$ \sum\limits_{\mathbf{r}\in F(P)\cap \Bbb{Z%
}^{10}}\frac{1}{\max \vert x_{i}^{3}\vert }\sim \mu _{\infty }(X(%
\Bbb{R}))\int\limits_{D(P)}dt.$$
To motivate this, let us approximate the sum on the left hand side by 
$$  \int_{F(P)}\frac{%
dudu_{1}du_{2}du_{3}dv_{1}dv_{2}dv_{3}dw_{1}dw_{2}dw_{3}}{\max_{1\leq j\leq
6}\vert x_{j}^{3}\vert d\left(
u_{1}v_{1}+u_{2}v_{2}+u_{3}v_{3}\right) }=m_{\infty }(F(P )).$$
Then by Lemma \ref{lemma17} we obtain $m_{\infty }(F(P))\sim  \mu _{\infty }(X(%
\Bbb{R}))\int_{D(P)}dt$ provided that 
the average contribution to $m_{\infty }(F(P))$ from the fibres of $\phi
_{O,\infty }: O(\Bbb{R})\rightarrow X(\Bbb{R})$ is the same as the contribution from the fibre over $x$.

%, we may define $v$-adic norms $\left\Vert . \right\Vert _{v}$ on 
%$\omega _{X}^{-1}$ from the previous norms $\left\Vert .\right\Vert _{W,v}$
%on $\omega _{W}^{-1}$ in \eqref{55} by letting 
%\begin{equation} 
 % \left\Vert f^{\ast }(\tau )(x_{v})\right\Vert
%_{v}=\left\Vert \tau (f(x_{v}))\right\Vert _{W,v} 
%\end{equation}
%for a local section $\tau $ of $\omega _{W}^{-1}$ where $\tau \circ f$ is
%defined at $w_{v}\in W(\Bbb{Q}_{v})$. It follows from \eqref{57} that
%\begin{equation} 
% (H\circ f)(x)=\prod\limits_{\text{all }v} \left\Vert
%\sigma (x)\right\Vert _{v}^{-1}
%\end{equation}
%for any $x\in X(\Bbb{Q})$ and any local section $\sigma $ of $\omega
%_{X}^{-1}$ where $\sigma (x)\neq 0$. \\

%To define $\tau _{H}(X),$ Peyre \cite[(2.2.1)]{P} uses the $v$-adic norms on $\omega _{X}^{-1}$
%appearing in the anticanoncial height function to defines a measure $ \mu
%_{v}$ on $X(\Bbb{Q}_{v})$ for each place $v$ of $\Bbb{Q}$ .\\

%Although one could use the definition of $\left\Vert .\right\Vert _{v}$ in \eqref{58} to define $\mu _{v}$,   it will be more useful to work with an
%slightly different definition of this $v$-adic norm based on the Poincar\'{e}
%residue map $\text{Res} : \omega _{\mathbf{\Xi }}(X)\rightarrow \iota _{\ast }\omega
%_{X}$ for the embedding $\iota :X\subset \mathbf{\Xi }$ into the toric
%subvariety $\Xi \subset \Bbb{P}^{5}\times \Bbb{P}^{2}\times \Bbb{P}%
%^{2}$ defined by \eqref{res5} and \eqref{res6}.  

\section{Preliminary upper bound estimates}\label{upper}

The rest of the paper  features analytic techniques, and it is convenient to introduce the following notation. Let $V(P )$ denote the number of integer sextuples $(\textbf{x}, \textbf{y})$ satisfying \eqref{1} and \eqref{nonzero} as well as the size condition $|x_j|\le P$, $|y_j|\le P$
$(1\le j\le 3)$. Since any rational point counted by $N(P)$ has exactly two representations
$(\mathbf x,\mathbf y) \in \Bbb{Z}^6$ with coprime coordinates, we conclude by one of M\"obius's
inversion formulae that
\begin{equation}\label{mobius}
   N(P) = \f12 \sum_{d=1}^\infty \mu(d) V(P^{1/3}/d).
\end{equation}  
The remainder of this paper is devoted to a proof of the
asymptotic relation
\be{asympV}
V(P)= P^3 Q_0(\log P) + O(P^{3-\tau})
\ee
in which $Q_0$ is a certain real polynomial of degree $4$ with leading 
coefficient 
\be{lead}
%\lim_{P\to\infty} \f{N(P)}{P^3(\log P)^4} = 
\frac{1}{2}(\pi^2 + 24\log 2 - 3)  \prod_p \Big(1-\f{9}{p^2}+\f{16}{p^3}-\f{9}{p^4}+\f{1}{p^6}\Big),
\ee
and $\tau$ is a suitable positive real number. Theorem \ref{thm1} follows easily from \eqref{mobius} once \eqref{asympV} and \eqref{lead}  are established. \\

The analytic counting 
procedures in the proof of Theorem \ref{thm1} will force us to implement a smooth approximation to the domain 
of counting. We start by deriving two upper bound estimates that help
 controlling the error in this transition. At the same time, this will
illustrate the use of Lemma \ref{lem1}. 
Our only additional tool is
the following simple estimate.
 
\begin{lem}\label{lem4} Let $A_1, A_2 \ge 1$, and let $u_1, u_2,u_3 \in\NN$
be coprime in pairs, 
  with $u_3\le A_2$. Then the number of solutions of $u_1 v_1 + u_2
  v_2 + u_3 v_3=0$ with $v_j\in\ZZ $ and $|v_1|\le A_1, |v_2|\le
  A_2 $ is $O(A_1 A_2 u_3^{-1})$.
  \end{lem}

To see this, note that we have to count solutions of $u_1 v_1 + u_2
v_2\equiv 0 \bmod u_3$. Choose $v_1$. Then one has to solve $u_2 v_2
\equiv c \bmod u_3$ for some $c\in\ZZ$. This has $O(1+A_2u_3^{-1})$ 
solutions in $v_2$.\\

Let ${\mathcal Z}\subset \{1,2,\ldots, [P]\}$ be a set
of $Z$ natural numbers, and let $V^*(P,{\mathcal Z})$ denote the number
of solutions of \rf{1} counted by $V(P)$ that satisfy $|y_j|\in{\mathcal
    Z}$ for at least one $j\in\{1,2,3\}$
Similarly, let $V_*(P,{\mathcal Z})$ denote the number
  of solutions of \rf{1} counted by $V(P)$ that satisfy $|x_j|\in{\mathcal
    Z}$ for at least one $j\in\{1,2,3\}$.

\begin{lem}\label{lem3} One has $V^*(P,{\mathcal Z})\ll
  P^{2+\eps} Z$ and  $V_*(P,{\mathcal Z})\ll
  P^{2+\eps} Z$.
%$$
%V^*(P,[1,Z])\ll P^2 Z (\log P)^3 (\log Z)^3.
%$$
\end{lem}

{\it Proof.} We begin with estimating $V^*(P,{\mathcal Z})$.
First observe that by symmetry in the
indices $1,2,3$, it suffices to estimate the number of solutions with
$|y_3|\in{\mathcal Z}$. %By  Lemma \ref{lem1}  and \rf{2.3}, it follows that $V^* (P,{\mathcal
%  Z})$ does not exceed $3$ times the number of $10$-tuples counted by
%$U(P)$ that satisfy the additional condition $z_3 d_1 d_2 d \in{\mathcal
 % Z}$. More explicitely,  by \rf{2.3}, \rf{2.9} and Lemma \ref{lem4},
 Hence, by Lemma \ref{lem1}, 
$$
V^*(P,{\mathcal Z})  \le 
3 \underset{\substack{ w_1 u_2 u_3 u \le P \\ w_2 u_1 u_3 u \le P\\ w_3 u_1 u_2 u
  \in{\mathcal Z} }}{\left.\sum\right.^{\ast}}
\trisum{v_1, v_2, v_3:}{|v_j w_j|\le P}{u_1 v_1 + u_2 v_2 + u_3 v_3=0}
1
$$
where $\sum^{\ast}$ indicates that the summation is subject to the coprimality conditions \eqref{2.2}. By Lemma \ref{lem4},  
$$V^*(P,{\mathcal Z}) \ll 
\trisum{w_1 u_2 u_3 u \le P}{w_2 u_1 u_3 u \le P}{w_3 u_1 u_2 u
  \in{\mathcal Z}}
\f{P^2}{w_1 w_2u_3} \ll P^2 (\log P)^3 \sum_{w_3 u_1 u_2 u\in{\mathcal Z}} 
  1.
$$
%where the very last sum is over all $z_3, d_1, d_2, d$ in $\NN$.
The
standard divisor estimate gives $O(P^\eps Z)$ for the last sum, as required. 

Now consider $V_*(P,{\mathcal Z})$. By symmetry, this quantity also does not
exceed 3 times the number of 
solutions counted by $V(P)$ with $x_1 \in {\mathcal Z}$. By Lemma \ref{lem1} and \eqref{2.9} we then see that   $V_*(P,{\mathcal Z})$ does
not exceed $3$ times the number of $10$-tuples $w_1, w_2, w_3$, $u,u_1,
u_2, u_3$, $v_1, v_2, v_3$ with $w_j, u_j, u \in\NN$, $v_j \in\ZZ$
satisfying
$$
w_1 u_2 u_3 u\le P, \, w_2 u_1 u_3 u\le P, \, w_3 u_1 u_2 u\le P, \,
(u_2; u_3)=1,
$$
and $u_1 v_1 + u_2 v_2 + u_3 v_3 =0$ with $|v_j w_j|\le P$ and $|v_1
w_1|\in{\mathcal Z}$. This leaves $ZP^\eps$ possibilities for $v_1, w_1$ 
by a divisor estimate. For given values of $u_1, u_2,u_3$, there are
$O(P/(w_2 u_3))$ possibilities for $v_2$ with $u_1 v_1 + u_2 v_2
\equiv 0 \bmod u_3$, and this fixes $v_3$ through the linear equation;
here we took advantage of the condition that $(u_2; u_3)=1$. It
follows that
$$
V_* (P,{\mathcal Z}) \ll \trisum{u_2 u_3 u\le P}{w_2 u_1 u_3 u \le P}{w_3
  u_1 u_2 u \le P}
\frac{ZP^{1+\eps}}{w_2 u_3} \ll ZP^{2+\eps},
$$
as required.\\

%Lemmas \ref{lem3} and \ref{lem6} are slightly more general than needed in our application. We will later apply them for expressions of the type 
%\begin{displaymath}
%  V^{\ast}(P(1+\Delta), [P, P(1+\Delta)] \cap \Bbb{Z}) \ll P^{3+\varepsilon} \Delta
%\end{displaymath}
%for some parameter $0 < \Delta < 1/50$. \\

Lemma \ref{lem3} shows that once one of the variables $x_j,y_j$ in \rf{1}
is restricted to a slim set, then there are few solutions.
However, more general
slim regions may well contain many integral points. 
An argument similar to the  above shows that there are $\gg P^3$ 
points on the subvariety defined by \eqref{1}, \eqref{2}, 
$y_1y_2y_3 \not=0$ and the additional condition $y_1 = y_2$. 
By symmetry, one may be tempted to reduce the evaluation of $V( P)$ to  
counting integral solutions in a  cone of the type $|y_1| \leq |y_2| \leq |y_3|$, but 
the error introduced from multiple counts of
 the subvarieties $y_1 = y_2$ and 
$y_2 = y_3$ is not negligible. This will cause extra difficulties later that 
we bypass with 
the introduction of a certain partition
of unity in \rf{defh}.

\section{Weights and integral kernels}\label{weights7}

% In the next sections, we will apply analytic methods to derive an asymptoticformula for $V(P)$. It will turn out to be convenient to count the solutions of \eqref{1} with $y_1y_2y_3 \not= 0$  in the box \eqref{2} with smooth weights, and then apply Mellin inversion to translate the quantity $V(P )$ into a multiple integral of certain Dirichlet series. Unfortunately a naive implementation of this programme would lead to highly divergent integrals, and the necessary regularization procedure requires considerable effort. 

In this section, we compile a number of technical results 
that will be needed in the analytical counting argument. The main topic 
is Mellin inversion for certain smooth and rough indicator functions.
Some  of the analysis is  routine. However, 
the two-dimensional Beta type kernels to be discussed in 
Lemmas \ref{lem10}, \ref{lem11} and \ref{lem12}
seem to be a new feature in the study of diophantine problems.

\subsection{Weight functions}\label{weights}

Counting problems are related to characteristic functions on appropriate 
regions, and the latter are discontinuous in a natural way. 
One obtains  smooth approximations by convolving them with a smooth 
approximate delta-distribution. The regions of relevance in this paper
are intervals, simplices, and a strip near the diagonal in the 
two-dimensional plane. 

The  smoothing will be 
controlled by two parameters $\delta$ and $\Delta$. We suppose from now on
that $\Delta \in (0, 1]$. All estimates will be uniform in $\Delta$.
In the end, we shall choose $\Delta$ as a small negative power of $P$.
The role of $\delta$ will be described in due course. 

The simplest smoothing is that of a sum over residue classes. Choose
a smooth function
$q : [0,\infty) \rightarrow [0, 1]$ with 
 $q=0$ on $[0, 1/4] \cup [7/4,\infty)$, and
$   q(x) + q(1+x) = 1$   
 for $x \in [0, 1]$.  Then $q^{(j)} \ll_j 1$ for all $j \in \Bbb{N}_0$.
Also, when $F:{\Bbb N} \to \CC$ is a function with period $D\in\NN$,
then    
\begin{equation}\label{trick}
\sum_{r=1}^D F(r) = \sum_{r=1}^\infty q(r/D) F(r).
\end{equation}

We proceed by smoothing the characteristic function of the interval $[0,1]$, 
denoted hereafter by $f_0$.
Let $\rho_{\Delta} : [0, \infty) \rightarrow [0,\infty)$ be a smooth non-negative function with
\begin{equation}\label{supprho}
  \text{supp}(\rho_{\Delta}) \subset (1, 1+\Delta)
\end{equation}
and 
\begin{equation}\label{volrho}
\int_0^{\infty} \rho_{\Delta}(x) \frac{dx}{x} = 1,
\end{equation}
and such that
\begin{equation}\label{diffrho}
  \rho_{\Delta}^{(j)} (x) \ll_j \Delta^{-1-j}
\end{equation}
holds for all $j \in \Bbb{N}_0$. For $x\in(0,\infty)$,  define
\begin{equation}\label{deff}
  f_{\Delta}(x) = \int_0^{\infty} \rho_{\Delta}(z) f_0\left(\frac{x}{z}\right) \frac{dz}{z} = \int_x^{\infty} \rho_{\Delta}(z) \frac{dz}{z}. 
\end{equation}
It follows from \eqref{supprho}, \eqref{volrho} and \eqref{diffrho} that 
\begin{equation}\label{suppfsimple}
   0\le f_{\Delta}(x)\le 1 \text{ for all } x\in[0,\infty), \quad  f_{\Delta} = 1 \text{ on } [0, 1], \quad \text{supp}(f_{\Delta}) \subset [0, 1+\Delta],  
 \end{equation}  
and that
\begin{equation}\label{difff}
f_{\Delta}^{(j)} \ll_j \Delta^{-j}
\end{equation}
holds for all $j \in \Bbb{N}_0$. We also note that
$ \text{supp}(f'_{\Delta}) \subset [1, 1+\Delta]$.
 Thus, $f_{\Delta}$ is indeed a smooth approximation 
to $f_0$.

For $n\in\NN$  
let 
$$
\mathcal{Q} =  {\mathcal Q}(n) = \{{\bf x}\in\RR^n: x_j \geq  0 \text{ for }  1\le j\le n\}
$$
denote the positive quadrant. 
Our next aim is to construct a certain smooth partition of unity of $\cal Q$.
 For $0 \le \delta < 1/10$, 
consider the  (infinite)  simplex
\begin{displaymath}
   \mathcal{T}_{\delta} = \{ \textbf{x} \in \mathcal{Q}: x_1 \leq (1+\delta)x_2 \leq \ldots \leq (1+\delta)^{n-1}x_n\}.
\end{displaymath}
   For $n \geq 2$,  define the function 
$h_{\delta} : \mathcal{Q} \rightarrow [0, 1]$ by
\begin{displaymath}
    h_{\delta}(\textbf{x}) = f_{\delta}(x_1/x_2) \cdot \ldots \cdot f_{\delta}(x_{n-1}/x_n)
\end{displaymath}
provided that $x_2 \cdots x_n \not= 0$, and put 
$h_{\delta}(\textbf{x}) = 0$ otherwise. 
It is easy to see that $h_{\delta}$ vanishes on  
$\mathcal{Q} \setminus \mathcal{T}_{\delta}$, and $h_{\delta}(\textbf{x}) = 1$ 
if $\textbf{x} \in \mathcal{T}_0$. 

The group $S_n$ acts on $\mathcal{Q}$ by permuting coordinates.  For $\pi \in S_n$  and $\textbf{x} \in \mathcal{Q}$ define 
\begin{equation}\label{defh}
  h_{\pi, \delta}(\textbf{x}) = \frac{h_{\delta}(\pi(\textbf{x}))}{ \sum_{\sigma \in S_n}h_{\delta}(\sigma(\textbf{x}))}. 
\end{equation}
 Note that the denominator is between 1 and $n!$ by construction of $h_{\delta}$. Clearly,  \begin{equation}\label{part}
   \sum_{\pi \in S_n} h_{\pi, \delta} = 1 \text{ on } \mathcal{Q} \quad \text{ and } \quad h_{\pi, \delta} = 0 \text{ on } \mathcal{Q} \setminus \pi(\mathcal{T}_{\delta}). 
 \end{equation}
The function $h_{\pi, 0}$ is simply the characteristic function on $\pi(\mathcal{T})$.

In our later work, we will use the functions $h_{\pi,\delta}$ only for $\delta =0$
and one specific positive value of $\delta$. It is therefore not necessary
to keep track of the dependence of implicit constants on $\delta$. Hence,
from now on, {\em implied constants may depend on} $\delta$.

For $\delta > 0$, the function  
$h_{\delta}$ is smooth on  $\mathcal{Q}$,  and it is  a simple exercise using 
\eqref{suppfsimple} to show that any fixed 
$(\nu_1, \ldots, \nu_n) \in \Bbb{N}_0^n$, the estimate
\begin{displaymath}
   \frac{\partial^{\nu_1}}{\partial x_1^{\nu_1}} \cdots  \frac{\partial^{\nu_n}}{\partial x_1^{\nu_n}} h_{\delta}(\textbf{x}) \ll \prod_{j=1}^{n} x_j^{-\nu_j} 
\end{displaymath}
holds uniformly in the range $ x_j > 0$ $(1\le j\le n)$. Consequently, 
for the same values of $\bf x$,
 \begin{equation}\label{diffh}
   \frac{\partial^{\nu_1}}{\partial x_1^{\nu_1}} \cdots  \frac{\partial^{\nu_n}}{\partial x_n^{\nu_n}} h_{\pi, \delta}(\textbf{x}) \ll \prod_{j=1}^{n} x_j^{-\nu_j}.
\end{equation}

Now we specialize to $n=3$. For $\pi \in S_3$ and 
$\textbf{x} \in \mathcal{Q}(3)$ we let
\begin{equation}\label{deffpi0delta}
  f_{\pi, 0, \delta}(\textbf{x}) = 
f_0(x_1)f_0(x_2)f_0(x_3) h_{\pi, \delta}(\textbf{x}).
\end{equation}
Then, $f_{\pi, 0, 0}$ is the characteristic function on the 
tetrahedron $0 \leq x_{\pi(1)} \leq x_{\pi( 2)} \leq x_{\pi(3)} \leq 1$.  
The corresponding smooth version is defined by
\begin{equation}\label{deffpi}
\begin{split}
  f_{\pi, \Delta, \delta}(\textbf{x}) &  =   \int_{\mathcal{Q}(3)} \rho_{\Delta}(z_1) \rho_{\Delta}(z_2) \rho_{\Delta}(z_3) f_{\pi, 0, \delta}\left(\frac{x_1}{z_1}, \frac{x_2}{z_2}, \frac{x_3}{z_3}\right) \frac{dz_1}{z_1} \frac{dz_2}{z_2} \frac{dz_3}{z_3}\\
 & = \int_{x_3}^{\infty} \int_{x_2}^{\infty} \int_{x_1}^{\infty} \rho_{\Delta}(z_1) \rho_{\Delta}(z_2) \rho_{\Delta}(z_3) h_{\pi, \delta}\left(\frac{x_1}{z_1}, \frac{x_2}{z_2}, \frac{x_3}{z_3}\right) \frac{dz_1}{z_1} \frac{dz_2}{z_2} \frac{dz_3}{z_3}. 
  \end{split}
\end{equation}
It is immediate from \eqref{part}, \eqref{deffpi0delta} and \eqref{deff} that
\begin{equation}\label{sym}
  \sum_{\pi \in S_3} f_{\pi, \Delta, \delta}(\textbf{x}) = f_{\Delta}(x_1) f_{\Delta}(x_2) f_{\Delta}(x_3).
\end{equation}
Note that the right hand side (and hence the left hand side) of \rf{sym} is 
independent of $\delta$. By \eqref{defh} and the last expression in \eqref{deffpi} we see that
\begin{equation}\label{suppf}
  \text{supp}(f_{\pi, \Delta, \delta}) \subset \{\textbf{x} \in \mathcal{Q}(3) \mid x_{\pi(1)} \leq \gamma x_{\pi(2)} \leq\gamma^2 x_{\pi(3)} \leq \gamma^3\}.
\end{equation}
where $\gamma = (1+\delta)(1+\Delta)$. 
Moreover, for $\delta>0$, the function $f_{\pi, \Delta, \delta}(\textbf{x})$ is 
smooth for each $\pi \in S_3$  and satisfies the crude bound 
\begin{equation}\label{difffpi}
   \frac{\partial^{\nu_1}}{\partial x_1^{\nu_1}}   \frac{\partial^{\nu_2}}{\partial x_2^{\nu_2}}   \frac{\partial^{\nu_3}}{\partial x_3^{\nu_3}}   f_{\pi, \Delta, \delta}(\textbf{x}) \ll \Delta^{-(\nu_1+\nu_2+\nu_3)}
\end{equation}   
for any fixed $(\nu_1, \nu_2, \nu_3) \in \Bbb{N}_0^3$, as can be seen from \eqref{diffh} and the last expression in \eqref{deffpi}. 
%\begin{displaymath}
%  f_{\pi, \Delta, \delta}(\textbf{x}) = \int_{x_3}^{\infty} \int_{x_2}^{\infty} \int_{x_1}^{\infty} \rho_{\Delta}(z_1) \rho_{\Delta}(z_2) \rho_{\Delta}(z_3) h_{\pi, \delta}\left(\frac{x_1}{z_1}, \frac{x_2}{z_2}, \frac{x_3}{z_3}\right) \frac{dz_1}{z_1} \frac{dz_2}{z_2} \frac{dz_3}{z_3}
%\end{displaymath}
%which in turn is inferred from \rf{deffpi0delta} and \rf{deffpi}.

Finally let $k^+_0, k^{-}_0 : \mathcal{Q}(2) \rightarrow [0, 1]$ be the characteristic functions on the sets 
\begin{equation}\label{defkpm}
  \{\textbf{x} \in \mathcal{Q}(2) : x_1 + x_2 \leq 1\}, \quad  \text{resp.}  \quad \{\textbf{x} \in \mathcal{Q}(2) : x_1 \leq 10, x_2 \leq 10, |x_1 - x_2| \leq 1\}. 
 \end{equation}
 Note that the region $|x_1 - x_2 | \leq 1$ has infinite intersection with $\mathcal{Q}(2)$, therefore we need an additional truncation. Define the smooth functions
\begin{equation}\label{defk}
  k^{\pm}_{\Delta}(x_1, x_2) = \int_{\mathcal{Q}(2)}  \rho_{\Delta}(z_1) \rho_{\Delta}(z_2) k^{\pm}_0\left(\frac{x_1}{z_1}, \frac{x_2}{z_2}\right) \frac{dz_1}{z_1} \frac{dz_2}{z_2}. 
\end{equation}
As before one then finds that for any fixed $\nu_1,\nu_2$ one has
\begin{equation}\label{diffk}
   \frac{\partial^{\nu_1}}{\partial x_1^{\nu_1}}   \frac{\partial^{\nu_2}}{\partial x_2^{\nu_2}}      k^{\pm}_{ \Delta}(\textbf{x}) \ll \Delta^{-(\nu_1+\nu_2)}.
\end{equation}   
It is also clear that
\begin{equation}\label{suppk+}
  k^{+}_\Delta(\textbf{x}) = 1 \text{ if } x_1 + x_2 \leq 1, \quad \text{supp}(k^+_\Delta) \subset \{\textbf{x} \in \mathcal{Q}(2) : x_1 + x_2 \leq 1 + \Delta\}
\end{equation}
and
\begin{equation}\label{suppk-}
\begin{split}
 & k^{-}_\Delta(\textbf{x}) = 1 \quad \text{ if } \quad |x_1 - x_2| \leq 1 - 10\Delta,  \,\, x_1, x_2 \leq 10,\\
  & \text{supp}(k^-_\Delta) \subset \{\textbf{x} \in \mathcal{Q}(2) : |x_1 - x_2| \leq 1 + 10 \Delta, \,\, x_1, x_2 \leq 10(1+\Delta)\}.
  \end{split}
\end{equation}
  
\subsection{Mellin inversion}   We begin with a short summary of 
well-known facts concerning Mellin transforms and the related
inversion theorem in a multidimensional set-up. For ${\bf s}\in\CC^n$ and ${\bf x}\in{\mathcal Q}$ write ${\bf x}^{\bf s}=
x_1^{s_1} x_2^{s_2} \cdots x_n^{s_n}$ in the interest of brevity, and
put ${\bf 1} = (1,1,\ldots,1)$.  A function $g:{\mathcal Q}\to \CC$ is {\it
  piecewise continuous} if it is continuous everywhere on ${\mathcal Q}$
except for a compact part of ${\mathcal Q}$ that is contained in the union
of finitely many $(n-1)$-dimensional submanifolds of $\RR^n$.  
Whenever $g:{\mathcal
  Q}\to\CC$ is a piecewise continuous, compactly supported and bounded function and
${\bf s}\in\CC^n$ with $\Re(s_j)>0$ $(1\le j\le n)$, then the integral
$$
\widehat{g}( {\bf s}) = \int_{\mathcal Q} g({\bf x}){\bf x}^{{\bf s}-{\bf
    1}} d{\bf x}
$$
defines a holomorphic function $\widehat{g}$. 
  If in addition $g$ is continuous, Mellin's inversion formula asserts
that for any ${\bf c}\in\RR^n$ with $c_j>0$ $(1\le j\le n)$ one has
\be{mellininv}
g({\bf x}) = \Big(\f{1}{2\pi i}\Big)^n \int_{({\bf c})} \widehat{g}(\textbf{s})
{\bf x}^{-{\bf s}} d{\bf s}.
\ee
Here and later, $\int_{({\bf c})}$ denotes $n$-fold integration over
the lines $s_j=c_j + it_j, t_j\in\RR$. \\

We consider the Mellin transforms of the weight functions that we 
defined in the previous section.  
%\begin{displaymath}
%\begin{split}
%  F_{\Delta}(s) & : =   \int_0^{\infty} f_{\Delta}(x) x^{s-1} dx, \quad Q(s) := \int_0^{\infty} q(x) x^{s-1} dx,\\
%   F_{\pi, \Delta}(\textbf{s}) & := \int_{\mathcal{Q}(3)} f_{\pi, \Delta}(\textbf{x}) \textbf{x}^{\textbf{s}-\textbf{1}} d\textbf{x},   \quad K^{\pm}_{\Delta}(\textbf{s})  :=  \int_{\mathcal{Q}(2)} k_{\Delta}^{\pm}(\textbf{x})\textbf{x}^{\textbf{s}-\textbf{1}} d\textbf{x}
 % \end{split}
%\end{displaymath}
%for $\Delta \geq 0 $ and $R_{\Delta}(s) = \int_0^{\infty} \rho_{\Delta}(x) x^{s-1} dx$ for $\Delta > 0$. 
The following estimates for the smooth versions are almost immediate.

\begin{lem}\label{lem10} Let $j_1, j_2, j_3 \in \Bbb{N}$.    \\
{\rm (i)} The function $\widehat{q}$ can be extended to an entire function satisfying $\widehat{q}(s) \ll_{j_1} (1+|s|)^{-j_1}$. \\
{\rm(ii)} The function $\widehat{f}_{\Delta}(s)$ is holomorphic in ${\rm Re\,} s > 0$ and satisfies $\widehat{f}_{\Delta}(s) \ll_{j_1} \Delta^{-j_1} |s|^{-j_1}$ in $1/10 < {\rm Re \,} s < 2$. \\
{\rm (iii)} For any $\pi\in S_3$ and $\delta>0$, the 
function $\widehat{f}_{\pi, \Delta, \delta}$ is holomorphic in ${\rm Re \,} s_j > 0$ and satisfies 
\begin{displaymath}
  \widehat{f}_{\pi, \Delta, \delta}(\textbf{s}) \ll_{j_1, j_2, j_3} \frac{\Delta^{-j_1-j_2-j_3}}{|s_1|^{j_1}|s_2|^{j_2}|s_3|^{j_3}} \quad \text{in }\frac{1}{10} < {\rm Re\,} s_j < 2.
\end{displaymath}
{\rm (iv)} The function 
$(s_1, s_2) \mapsto s_1s_2 \widehat{k}^{\pm}_{\Delta}(s_1, s_2)$ admits an 
analytic continuation to ${\rm Re \,} s_j > -1$ and satisfies
\begin{displaymath}
  \widehat{k}_{\Delta}^{\pm}(s_1, s_2) \ll_{j_1, j_2} \frac{\Delta^{-j_1-j_2}}{|s_1|^{j_1}|s_2|^{j_2}} \quad \text{in } -1/2 \leq {\rm Re \,} s_j \leq 2, \quad |s_j | \geq 1/10. 
\end{displaymath}
 \end{lem}
 
 {\it Proof.} This is repeated integration by parts in combination with  \eqref{difff}, \eqref{difffpi} and \eqref{diffk}. \\
  
%{\it Proof.} This follows immediately after repeated integration by parts.\\
%The holomorphy in (i) follows from the fact that $q$ vanishes in a neighbourhood of the origin; the claimed estimate follows after $j_1$ integrations by parts. The proof of (ii) is similar: we integrate by parts $j_1$ times and use \eqref{supp}. In (iii) we integrate by parts $j_1$ times with respect to $x_1$, $j_2$ times with respect to $x_2$ and $j_3$ times with respect to $x_3$ and recall \eqref{diffh}. Similarly,  we integrate by parts $j_1$ times with respect to $x_1$ and $j_2$ times with respect to $x_2$ to prove (iv). The resulting integral converges absolutely and locally uniformly in $\Re s_j > -1$, and \eqref{bounds} follows on using once again \eqref{supp}. \\

More precise statements are possible for the unsmoothed weight functions.
\begin{lem}\label{lem11}  {\rm (i)} One has $\widehat{f}_0(s) = 1/s$.\\
{\rm (ii)} For $\pi \in S_3$ one has
\begin{displaymath}
  \widehat{f}_{\pi, 0, 0}(s_1, s_2, s_3) = \frac{1}{s_{\pi(1)}   (s_{\pi(1)} + s_{\pi(2)} )(s_{\pi(1)} + s_{\pi(2)} + s_{\pi(3)})}.
\end{displaymath}
{\rm (iii)} 
 The function 
$(s_1, s_2) \mapsto s_1s_2 \widehat{k}^{\pm}_{0}(s_1, s_2)$ admits an 
analytic continuation to ${\rm Re \,} s_j > -1$. 
For fixed $s_2$ with ${\rm Re\,} s_2 > 0$, the function  
$s_1\mapsto \widehat{k}_0^{\pm}(s_1, s_2)$ is meromorphic in a 
neighbourhood of $s_1 = 0$. At $s_1 = 0$ there is a simple pole with 
\begin{displaymath}
  \underset{s_1=0}{\text{\rm res}}\widehat{k}_0^{\pm}(s_1, s_2) = \frac{1}{s_2}.
\end{displaymath}
Moreover, $\widehat{k}_0^{\pm}(s_1, s_2) = \widehat{k}_0^{\pm}(s_2, s_1)$. %For fixed $s_1$ with $\Re s_1 >0$, the function 
%$s_2 \mapsto\widehat{k}_0^{\pm}(s_1, s_2)$ is meromorphic in a 
%neighbourhood of $s_2 = 0$. At $s_2 = 0$ there is a simple pole with
%\begin{displaymath}
 %\underset{s_2=0}{\text{res}}\widehat{k}_0^{\pm}(s_1, s_2) = \frac{1}{s_1}.
%\end{displaymath} 
\end{lem}

{\it Proof.} The proof of (i) and (ii) is a straightforward calculation. 
Next, observe that the integral representation of the Euler 
Beta-function \cite[8.380.1]{GR} yields 
\begin{equation}\label{explicit}
  \widehat{k}^{+}_{0}(s_1, s_2) = \frac{\Gamma(s_1)\Gamma(s_2)}{\Gamma(1+s_1+s_2)}.  
\end{equation}
Hence the claims in (iii) for $k^+_0$ are a consequence of elementary properties
of the Gamma function. For  $k^-_0$, we need to work directly from the definition.
When $\Re s_j >0$ $(j=1,2)$, Fubini's theorem gives
$$
\widehat{k}^{-}_0(s_1, s_2)=\int_0^{10}   x_2^{s_2-1}  
\int_{\max(0, x_2-1)}^{\min(x_2+1, 10)}     x_1^{s_1-1}\,dx_1\,dx_2.
$$
Here, in the inner integral, we extend the integration over $[0,10]$,
and subtract the terms added in to correct the error. This artifice produces the
identity 
\begin{equation}\label{intbyparts}
\widehat{k}^{-}_0(s_1, s_2) =
 \frac{10^{s_1+s_2}}{s_1s_2} - \frac{1}{s_1} F(s_2,s_1) - \frac{1}{s_2} F(s_1,s_2) 
\end{equation}
in which
$$ F(s_1,s_2) =  \int_0^9 (x+1)^{s_1-1}x^{s_2}\, dx. $$
This integral defines $F$ as a holomorphic function in $s_1\in\CC$ and $\Re s_2 > -1$.
Thus, \rf{intbyparts} provides the desired continuation of  $s_1s_2 \widehat{k}^{-}_0(s_1, s_2)$. Also, when $\Re s_1 >0$ is fixed, it follows from \rf{intbyparts} that 
$\widehat{k}^{-}_0(s_1, s_2)$ has a simple pole at $s_2=0$. Its residue is the
value at $s_2=0$ of the function $s_2\widehat{k}^{-}_0(s_1, s_2)$, and hence
equals
$$ \frac{10^{s_1}}{s_1} - F(s_1,0) = \f{1}{s_1}.$$
The final symmetry statement is clear. %The same argument applies with $s_1,s_2$ interchanged. 
This completes the proof.\\

%By \eqref{defk}, it follows that
%\begin{equation}
%  \widehat{k}^{+}_{\Delta}(s_1, s_2) = \frac{\Gamma(s_1)\Gamma(s_2)}{\Gamma(1+s_1+s_2)}  \widehat{\rho}_{\Delta}(s_1)\widehat{\rho}_{\Delta}(s_2). 
%\end{equation}
 
%{\it Proof of Lemma \ref{lem11}.} The proof of For the proof of (i) we can assume $\pi = (1)$, and an elementary computation shows
%\begin{displaymath}
  %\underset{0 \leq x_1 \leq x_2 \leq x_3 \leq 1}{\int\int\int} \textbf{x}^{\textbf{s}- \textbf{1}}d\textbf{x} = \frac{1}{s_1(s_1+s_2)(s_1+s_2+s_3)}.
%\end{displaymath}

%\medskip

Our final lemma  will eventually estimate the error when we remove the smoothing at the end of the argument. The proof is rather long and technical and will occupy the rest of this section. %It can safely be skipped at a first reading.

\begin{lem}\label{lem12} Let $0 \leq \eta \leq 1$. \\
{\rm (i)} For $1/10 < {\rm Re \,} s < 2$ we have 
\begin{displaymath}
  \widehat{f}_{\Delta}(s) - \widehat{f}_0(s) \ll \Delta^{\eta} |s|^{\eta-1} \quad \text{and} \quad   \max(\widehat{f}_0(s), \widehat{f}_{\Delta}(s)) \ll |s|^{-1}.
  \end{displaymath} 
{\rm (ii)} For $1/10 < {\rm Re\,} s_j < 2$, $\delta > 0$  and $\pi \in S_3$ we have
\begin{displaymath}
  \widehat{f}_{\pi, \Delta, \delta}(\textbf{s}) - \widehat{f}_{\pi, 0, \delta}(\textbf{s}) \ll  \Delta^{3\eta}|s_1s_2s_3|^{\eta-1} \quad \text{and} \quad 
%\end{displaymath}
%and 
%\begin{displaymath}
  \max(\widehat{f}_{\pi, 0, \delta}(\textbf{s}), \widehat{f}_{\pi, \Delta, \delta}(\textbf{s})) \ll |s_1s_2s_3|^{-1}. 
\end{displaymath}  
{\rm (iii)} %The function $(s_1, s_2) \mapsto s_1s_2 K_0^{\pm}(s_1, s_2)$ is holomorphic in $\Re s_j > -1$. 
Fix   $-1/2 \leq \alpha \leq 0$. Then in the region
\begin{equation}\label{region}
  \alpha \le {\rm Re \,} s_j \le 2, \quad {\rm Re\,} s_1 + {\rm Re \,} s_2 \geq 1/10,  \quad |s_j| \geq 1/10
  \end{equation}
we have  
\begin{displaymath}
   \widehat{k}^{\pm}_{\Delta}(\textbf{s})  -\widehat{k}^{\pm}_0(\textbf{s})  \ll  \frac{\Delta^{\eta}}{\max(|s_1|, |s_2|)^{1+\alpha-\eta}\min(|s_1|, |s_2|)^{1/2-\alpha}}
\end{displaymath}
and 
\begin{displaymath}
  \max\bigl(\widehat{k}^{\pm}_0(\textbf{s}), \widehat{k}^{\pm}_{\Delta}(\textbf{s})\bigr) \ll \frac{1}{\max(|s_1|, |s_2|)^{1+\alpha}\min(|s_1|, |s_2|)^{1/2-\alpha}}.
\end{displaymath}   
The bounds in $\text{\rm (i), (ii), (iii)}$ 
remain valid if the functions on the left hand side of the displayed estimates are replaced by a partial derivative. 
\end{lem}

{\it Proof.}  The remark on the derivatives follows 
immediately from Cauchy's integral formula. 
We launch the proof of (i), (ii) and (iii) 
with a deduction of the auxiliary bound  
\begin{equation}\label{diff}
  \widehat{\rho}_{\Delta}(s) = 1 + O( (\Delta |s|)^{\eta})
\end{equation}
that holds for any fixed $0 \leq \eta \leq 1$ uniformly for $1/10 \leq \Re s \leq 2$. Indeed, it follows easily from \eqref{supprho} and \eqref{diffrho} that 
\begin{equation}\label{R1}
  \widehat{\rho}_{\Delta}(s) \ll 1.
\end{equation}
On the other hand,  from \eqref{deff} we conclude 
\begin{equation}\label{PsiR}
     \widehat{f}_{\Delta}(s) - \widehat{f}_0(s) = \widehat{\rho}_{\Delta}(s) \widehat{f}_0(s) - \widehat{f}_0(s)  = \frac{1}{s}(\widehat{\rho}_{\Delta}(s) - 1),
 \end{equation}    
 hence
\begin{equation}\label{R2}
  \widehat{\rho}_{\Delta}(s) - 1 = s \int_0^{\infty} (f_{\Delta}(x) - f_0(x)) x^{s-1} dx \ll |s| \Delta. 
\end{equation}
Clearly \eqref{R1} and \eqref{R2}  imply \eqref{diff}. 

We are now prepared for the main argument. We have already shown part (i), as an inspection of \eqref{diff} -- \eqref{PsiR} shows.
For the proof of part (ii), we  note that \eqref{deffpi} yields
\begin{displaymath}
  \widehat{f}_{\pi, \Delta, \delta}(\textbf{s}) = \widehat{f}_{\pi, 0, \delta}(\textbf{s}) \widehat{\rho}_{\Delta}(s_1) \widehat{\rho}_{\Delta}(s_2) \widehat{\rho}_{\Delta}(s_3).
\end{displaymath}
Hence by  \eqref{diff}  it is enough to show
\begin{equation}\label{needtoshow}
  \widehat{f}_{\pi, 0, \delta} \ll |s_1s_2s_3|^{-1}.
\end{equation}
This can be seen as follows. By \eqref{deffpi0delta} we have 
$$
  f_{\pi, 0, \delta}(\textbf{x}) = F_0(\textbf{x}) H_{\pi, \delta}(\textbf{x})
$$
where $F_0(\textbf{x}) = f_0(x_1) f_0(x_2)f_0(x_3)$ and $H_{\pi, \delta}(\textbf{x}) = h_{\pi, \delta}(\textbf{x}) f_1(x_1)f_1(x_2)f_1(x_3)$. It follows easily from \eqref{diffh} that the Mellin transform $\widehat{H}_{\pi, \delta}
(\textbf{s})$ is holomorphic in $\Re s_j > 0$ and rapidly decaying on vertical lines. We observe that
\begin{displaymath}
 \widehat{ f}_{\pi, 0, \delta}(\textbf{s}) = \frac{1}{(2\pi i)^3} \int_{(\textbf{c} )} \widehat{F}_0(\textbf{s} - \textbf{t}) \widehat{H}_{\pi, \delta}(\textbf{t}) d\textbf{t} =  \frac{1}{(2\pi i)^3} \int_{(\textbf{c} )} \prod_{j=1}^3 (s_j - t_j)^{-1} 
\widehat{H}_{\pi, \delta}( \textbf{t}) d\textbf{t}
\end{displaymath}
for $0  < c_j < \Re s_j$, %any $\textbf{c} \in (0,\infty)^3$, 
and \eqref{needtoshow} follows easily.

Finally we prove (iii). The convolution formula \eqref{defk} yields 
%\begin{displaymath}
  $\widehat{k}^{\pm}_{\Delta}(\textbf{s}) = \widehat{k}^{\pm}_0(\textbf{s}) \widehat{\rho}_{\Delta}(s_1)\widehat{\rho}_{\Delta}(s_2)$, 
hence
%\end{displaymath}
 (iii) follows from \eqref{diff}, once we can show that the bound
 \begin{equation}\label{K2}
  \widehat{k}_0^{\pm}(\textbf{s}) \ll  \frac{1}{\max(|s_1|, |s_2|)^{1+\alpha} \min(|s_1|, |s_2|)^{1/2-\alpha}}
\end{equation}  
%whenever $\alpha \leq  \Re s_j \leq 2$ and $ |s_j| \geq 1/10$. 
holds uniformly for all $\mathbf s$ in the region defined by
\eqref{region}. For $\widehat{k}^{+}_0$,  \rf{K2} follows 
from \eqref{explicit} in combination with Stirling's formula. 
In the absence of such an explicit formula for $\widehat{k}^{-}_0$, 
we return to \rf{intbyparts} and observe that for $\mathbf s$ 
in accordance with \rf{region}, the first summand on the
right hand side of this identity does not exceed $10^4 |s_1s_2|^{-1}$, which in
turn is bounded by the right hand side of \rf{K2}. By symmetry, it therefore
suffices to establish that in the region described by \rf{region}, one has
\be{001}
\f1{s_2} F(s_1,s_2) \ll \frac{1}{\max(|s_1|, |s_2|)^{1+\alpha} \min(|s_1|, |s_2|)^{1/2-\alpha}}.
\ee

We write $s_j=\sigma_j+i t_j$ with real numbers $\sigma_j$, $t_j$, and 
define the functions $\phi, \omega: (0,\infty) \to \Bbb R$ by
$$ \omega(x) = (x+1)^{\sigma_1-1} x^{\sigma_2}, \quad \phi(x) 
=  t_1 \log(x+1)+  t_2\log x. $$
Then, the definition of $F$ may be rewritten as
\be{002} F(s_1,s_2) = \int_0^9 \omega(x) e^{i\phi(x)} \,dx, \ee
and we observe that
\be{003} \overline{F(s_1,s_2)} = F(\bar s_1, \bar s_2). \ee

The bound \rf{001} is trivial in the compact part of \rf{region} described by
$|t_1|\le 2$, $|t_2|\le 2$. In the range $|t_1|\le 2$, $|t_2|\ge 2$ one has
$\min(|s_1|,|s_2|) \asymp 1$ and $\max(|s_1|,|s_2|) \asymp |s_2|$. Also, using
$\al\le\sigma_j \le 2$, we see that
$$ |F(s_1,s_2)| \le \int_0^9 \omega(x)\,dx \le \int_0^9 (x+1) x^{\sigma_2} dx
\ll 1, $$
and \rf{001} follows. 

It remains to discuss the case where $|t_1|\ge 2$. Here the treatment 
must  be based
on integration by parts, enhanced by a stationary phase argument 
when necessary. Whenever
$0<a\le b\le 9$ and $\phi'$ does not vanish on the interval $[a,b]$,
partial integration with respect to the phase $\phi$ yields the identity 
\be{004} 
i \int_a^b \omega(x) e^{i\phi(x)}\, dx
= \frac{\omega(b)}{ \phi'(b)} e^{i\phi(b)} - 
\frac{\omega(a)}{ \phi'(a)} e^{i\phi(a)}
-  \int_a^b \left( \frac{\omega'(x)}{\phi'(x)} - 
\frac{\omega(x)\phi''(x)}{\phi'(x)^2}  \right)  e^{i\phi(x)} dx.
\ee

 Now consider the situation where $2\le |t_1|\le \f{11}{10} |t_2|$. By \rf{003}, we may assume that $t_2\ge 2$. 
Then, for $0<x\le 9$, one has 
\be{phistrich} x\phi'(x)= t_2 + \f{t_1x}{x+1}\ge \f{t_2}{100}.\ee
In particular, $x\phi'(x)$ is continuous at $x=0$, and bounded below. 
By \rf{002} and \rf{004} with $a\to 0$, we then have 
 \begin{equation}\label{operator}
    i F(s_1,s_2) =  
\frac{\omega(9)}{ \phi'(9)} e^{i\phi(9)}  - % \int_0^9  \left(\frac{
%  \omega'(x)}{\phi'(x)} + \frac{\omega(x)\phi''(x)}{(\phi'(x))^2}  \right)  e^{i\phi(x)} dx = 
  \int_0^9 \left( \frac{x\omega'(x)}{x\phi'(x)} - \frac{\omega(x) \cdot x^2\phi''(x)}{(x\phi'(x))^2}  \right)  e^{i\phi(x)} dx. 
\end{equation}
However, for $0<x\le 9$, $\al\le \sigma_j\le 2$
  one has $\omega(x) \ll x^{\sigma_2}$, and a short calculation also confirms that
$$ x\omega'(x) \ll x^{\sigma_2} , \quad x^2\phi''(x) = -t_2-t_1\Bigl(\f{x}{x+1}\Bigr)^2 \ll t_2.$$ 
It is now immediate that the integral on the right hand side of \rf{operator} 
is $O(t_2^{-1})$, and consequently, the left hand side of \rf{001} is 
$O(|s_2|^{-2})$, which is more than is required.

By \rf{003}, the only case that remains to be discussed is when $t_1>2$
and $|t_2| \le \f{10}{11} t_1$. We write
$$ \eta = \f{1+ |t_2|^{1/2}}{8t_1}, \quad \xi = \f{-t_2}{t_1+t_2}, $$
and note that $0<\eta <\f14$. If in addition to the 
current assumptions one has $t_2\ge 0$, then by \rf{phistrich} we have
$\phi'(x) \ge t_1/10$ throughout the range $0< x\le 9$. Hence, we may take
$a=4\eta$, $b=9$ in \rf{004} and, subject to \rf{region}, 
estimate by brute force to establish the bound
\be{006}  \int_{4\eta}^9 \omega(x) e^{i\phi(x)}\, dx
\ll \f{\eta^\alpha}{t_1} + \f{1}{t_1^2} \int_{4\eta}^9 \omega(x) |\phi''(x)|\,
dx \ll  \f{\eta^\alpha}{t_1} +  \f{\eta^{\alpha-1}t_2}{t_1^2}. \ee
Also, by direct estimates,
\be{007}  \int_0^{4\eta} \omega(x) e^{i\phi(x)}\, dx \ll 
\int_0^{4\eta} x^{\sigma_2}\,dx
\ll \eta^{1+\alpha}. \ee
By \rf{002}, this combines to 
$$ F(\mathbf s) \ll \eta^\alpha \Bigl( \f{1}{t_1} + \f{t_2}{\eta t_1^2}
+ \eta\Bigr) \ll \Bigl( \f{1+t_2^{1/2}}{t_1} \Bigr)^{1+\alpha}, $$
and \rf{001} is immediate.

In the case where $t_2<0$ the function $\phi$ has exactly one zero
on $[0,\infty)$, at $x=\xi$. In the peculiar situation where $\eta>\f12\xi$
the previous argument needs little adjustment: \rf{007} remains valid as it 
stands, and the identity
\be{xphistrich}  x\phi'(x) = \f{t_1+t_2}{1+x} (x-\xi) \ee
shows that for $x\ge 4\eta$ one has $x\phi'(x) \ge \f1{110} t_1(x-\xi) \ge 
\f{1}{220} t_1 x$, which yields \rf{006}, and \rf{001} follows as before.

We may now suppose that  $\eta\le\f12\xi$. Also, the currently active
constraints that $0\le -t_2\le \f{10}{11} t_1$ and $t_1\ge 2$ imply that
$0\le \xi\le 10$. We partition the interval $[0,9]$ into the three disjoint
set  
$$ J_1 = [0,9]\cap (-\infty, \xi-\eta], \quad
J_2 = [0,9]\cap (\xi-\eta,\xi+\eta), \quad
J_3 = [0,9]\cap [\xi+\eta,\infty).
$$
Then, one routinely finds that
$$ \int_{J_2} \omega(x) e^{i\phi(x)}\, dx \ll \int_{J_2} x^{\sigma_2}\, dx 
\ll \eta\xi^\alpha \ll \f{|t_2|^{1/2+\alpha}}{t_1^{1+\alpha}}. $$
Next, we note that for $x\in J_1$, one finds from \rf{xphistrich} the 
lower bound
$x|\phi'(x)|\ge \f1{110} t_1(\xi-x)$. Thus, by \rf{004} with $a\to 0$ and
 $b=\min(9,\xi-\eta)$,
\begin{displaymath}
\begin{split}
\int_{J_1} \omega(x) e^{i\phi(x)}\, dx & \ll 
\left|\frac{b\omega(b)}{b \phi'(b)}\right|  
+  \int_0^b \left|\frac{x\omega'(x)}{x\phi'(x)}\right| \,dx 
+ \int_0^b
\frac{x^2\omega(x)|\phi''(x)|}{(x\phi'(x))^2}  \, dx.
\\
& \ll \f{\xi^{\alpha+1}}{t_1\eta} + \f{1}{t_1} \int_0^b  \f{x^{\sigma_2}}{\xi-x}\,dx
+ \f{1}{t_1^2} \int_0^b \f{x^{\sigma_2}(t_1x^2 + |t_2|)}{(\xi-x)^2}\,dx \ll  \f{|t_2|^{1/2+\alpha}}{t_1^{1+\alpha}}.
\end{split}
\end{displaymath}
A very similar calculation estimates the integral with $J_1$ replaced by $J_3$
to the same precision. By \rf{002}, the desired estimate \rf{001} is now 
immediate. This completes the proof of Lemma \ref{lem12}.\\

It should be observed that our method of estimation for $\widehat{k}^{-}$
applies equally well to $\widehat{k}^{+}$. Thus, the precise information
that is contained in \rf{explicit} is, strictly spreaking, not required
in this paper. On the other hand, Stirling's formula may be used to show
that our bound for $\widehat{k}^{+}$ is essentially the best possible. Thus
the decay of $\widehat{k}^{+}$ as $|s_j|$ gets large 
is far too weak to be in $L^1$. This will
cause serious technical difficulties in Chapter \ref{remove}.

%\textbf{Remark:} While $\widehat{f}_{\pi, 0, \delta}(\textbf{c} + i\textbf{t})$ fails to be in $L^1(\Bbb{R}^3, d\textbf{t})$ only by an $\varepsilon$, the  kernels $\widehat{k}^{\pm}_0$ fail  to be in $L^1$ much more severely. As long as we are sheltered by the smoothing parameter $\Delta$, this is not a problem, but sending $\Delta$ to $0$ will cause some technical difficulties in Section \ref{remove}. The occurrence of such Beta-type kernels seems to be a new feature in analytic methods in the context of Manin's conjecture.  

\section{Analytic Methods}\label{analytic}

\subsection{Counting with weights}\label{sec81} We are now prepared to prove \eqref{asympV}.  
The counting function $V(P)$ will   be ``smoothed'' in several steps to
facilitate its evaluation by Mellin inversion and Dirichlet series techniques
in the following sections. Let $0 < \Delta \leq 1/10$. We begin by writing the definition of $V(P)$
in the form
$$ V(P) = \sum_{(\mathbf x,\mathbf y) \in\cal W} \prod_{j=1}^3 
f_0\left(\frac{|x_j|}{P}\right)
f_0\left(\frac{|y_j|}{P}\right),
$$
and then replace all $f_0(|y_j|/P)$ by $f_\Delta(|y_j|/P)$. This will
produce an error which in view of \rf{suppfsimple} and the 
notation introduced in the
preamble to Lemma \ref{lem3} is bounded by the quantity $V^*((1+\Delta)P,[P,(1+\Delta)P])$.
We now invoke the first conclusion in Lemma \ref{lem3} and apply \eqref{sym}
to replace the product $f_\Delta(|y_1|/P)f_\Delta(|y_2|/P)f_\Delta(|y_3|/P)$.
For any fixed $\delta\in(0,1/10)$, this yields
\be{mein1}
V(P ) =  \sum_{\pi \in S_3}\sum_{(\textbf{x}, \textbf{y})\in\cal W} 
f_{\pi, \Delta, \delta}\left(\frac{|y_1|}{P}, \frac{|y_2|}{P}, \frac{|y_3|}{P}\right)f_{0}\left(\frac{|x_1|}{P}\right)f_{0}\left(\frac{|x_2|}{P}\right)f_{0}\left(\frac{|x_3|}{P}\right) + O(\Delta P^{3+\varepsilon}). 
\ee
For any $\pi\in S_3$, one has $(\mathbf x, \mathbf y)\in\cal W$ if and only
if  $(\pi\mathbf x, \pi\mathbf y)\in\cal W$. Hence, the inner sum on the
right hand side of \rf{mein1} is independent of $\pi$, so that
it suffices to consider henceforth the contribution from the 
identical permutation, id. Now insert the parametrization obtained in Lemma
\ref{lem2}, and note that the contributions 
to \rf{mein1} do not depend on the sign of the $w_j$. This produces 
the estimate
   \be{NEU}
 \begin{split}
 V(P ) =  48 &\underset{(u, \textbf{u}, \textbf{w}) \in \Bbb{N}^7}{\left.\sum\right.^{\ast}} f_{\text{id}, \Delta, \delta}\left(\frac{uu_2u_3w_1}{P}, \frac{uu_1u_3w_2}{P}, \frac{uu_1u_2w_3}{P}\right) \Upsilon(\textbf{u}, \textbf{w})  +  O(\Delta P^{3+\varepsilon})
   \end{split}
\ee
in which $\sum^{\ast}$ indicates the coprimality conditions \eqref{2.2} and  \eqref{2.5}, and we wrote
\begin{displaymath}
  \Upsilon(\textbf{u}, \textbf{w}) = \sum_{ r_1 \in \mathcal{S}(d_1) } \sum_{r_2, r_3 \in \Bbb{Z}}   f_{0}\left(\frac{|u_2r_3-u_3r_2|w_1}{P}\right)f_{0}\left(\frac{|u_3r_1-u_1r_3|w_2}{P}\right)f_{0}\left(\frac{|u_1r_2-u_2r_1|w_3}{P}\right).
\end{displaymath} 

Next, we smooth out the sum $ \Upsilon(\textbf{u}, \textbf{w})$. The inner sum
over $r_2,r_3$ in the definition of this sum depends only on $r_1\bmod u_1$.
This follows from the concluding remark  prior to Lemma \ref{lem2}, for example.
Hence, by \rf{trick}, one infers that
 \begin{displaymath}
  \Upsilon(\textbf{u}, \textbf{w}) = \sum_{ r_1=1 }^\infty \sum_{r_2, r_3 \in \Bbb{Z}}
q\left(\f{r_1}{u_1}\right)   f_{0}\left(\frac{|u_2r_3-u_3r_2|w_1}{P}\right)f_{0}\left(\frac{|u_3r_1-u_1r_3|w_2}{P}\right)f_{0}\left(\frac{|u_1r_2-u_2r_1|w_3}{P}\right).
\end{displaymath}
Let
 $  \Upsilon^{(1)}(\textbf{d}, \textbf{z})$ be the contribution 
to this sum from terms
with $r_2r_3\not= 0$, let  $   \Upsilon^{(2)}(\textbf{d}, \textbf{z}) $ be the contribution from terms with $r_2 = 0$, $r_3 \not= 0$, and let $   \Upsilon^{(3)}(\textbf{d}, \textbf{z})$ be the contribution with $r_2=r_3 = 0$.
Then, by symmetry, 
 $$\Upsilon(\textbf{u}, \textbf{w}) =   \Upsilon^{(1)}(\textbf{u}, \textbf{w})+  2 \Upsilon^{(2)}(\textbf{u}, \textbf{w})+   \Upsilon^{(3)}(\textbf{u}, \textbf{w}).$$

 %This induces a corresponding decomposition 
%\begin{equation}\label{decompV}
 % V_{\Delta}(P ) = V^{(1)}_{\Delta}(P ) + 2V^{(2)}_{\Delta}(P ) + V^{(3)}_{\Delta}(P ). 
%\end{equation}  

We rewrite the sums defining $\Upsilon^{(j)}(\textbf{d}, \textbf{z})$ by sorting
the sum according to the signs of $r_2,r_3$. This yields the identities
\begin{displaymath}
\begin{split}
    \Upsilon^{(1)}(\textbf{u}, \textbf{w}) & = \sum_{\mathbf r \in {\Bbb N}^3} \sum_{\sigma_2, \sigma_3 \in \{\pm 1\}  }q\left(\frac{r_1}{u_1}\right)    f_{0}\left(\frac{|u_2\sigma_3r_3-u_3\sigma_2r_2|w_1}{P}\right)
      f_{0}\left(\frac{|u_3r_1-u_1\sigma_3r_3|w_2}{P}\right)
\\ & \hspace{25em}\times f_{0}\left(\frac{|u_1\sigma_2r_2-u_2r_1|w_3}{P}\right),\\
        \Upsilon^{(2)}(\textbf{u}, \textbf{w}) & =  \sum_{ r_1, r_3 = 1}^\infty\sum_{\sigma_3 \in \{\pm 1\}}  q\left(\frac{r_1}{u_1}\right)    f_{0}\left(\frac{u_2r_3w_1}{P}\right)f_{0}\left(\frac{|u_3r_1-u_1\sigma_3r_3|w_2}{P}\right)f_{0}\left(\frac{ u_2r_1w_3}{P}\right),\\
        \Upsilon^{(3)}(\textbf{u}, \textbf{w}) & =  \sum_{ r_1= 1}^\infty q\left(\frac{r_1}{u_1}\right)     f_{0}\left(\frac{u_3r_1 w_2}{P}\right)
f_{0}\left(\frac{u_2r_1 w_3}{P}\right).
   \end{split}
\end{displaymath}
The support conditions of $q$ and $f_{\text{id}, \Delta,\delta}$ in \eqref{suppf} 
imply that non-zero contributions to the sum \rf{NEU} only arise  from summands with
$1 \leq r_1 \leq 2u_1$  and  $uu_2u_3w_1\leq \frac{5}{4}uu_1u_3w_2 \leq \frac{25}{16}uu_1u_2w_3$. 
Similarly, since all $(\mathbf x, \mathbf y)\in\cal W$ that occur in \rf{mein1}
with a non-zero weight satisfy $|x_j|\le P$, $|y_j|\le (1+\Delta)P$, we deduce from Lemma \ref{lem2A} that $r_iu_jw_k \le 8P$ holds for any choice of
$\{i,j,k\}=\{1,2,3\}$. Therefore, by \eqref{defkpm}, we may rewrite 
expressions like $f_0(|u_2\sigma_3 r_3 - u_3\sigma_2 r_2|/P)$ in terms of
$k_0^\pm(u_2r_3/P, u_3r_2/P)$. With
\begin{displaymath}
  E = \{(-, -, -), (+,+, -), (+, -, +), (-, +, +)\}
\end{displaymath}
this produces
 \begin{displaymath}
\begin{split}
    \Upsilon^{(1)}(\textbf{u}, \textbf{w}) & = \sum_{ \mathbf r \in \Bbb N^3} \sum_{{\bm \epsilon} \in E} q\left(\frac{r_1}{u_1}\right)    k^{\epsilon_1}_{0}\left(\frac{u_2r_3w_1}{P}, \frac{u_3r_2w_1}{P}\right)  k^{\epsilon_2}_{0}\left(\frac{u_3r_1w_2}{P}, \frac{u_1r_3w_2}{P}\right)k^{\epsilon_3}_{0}\left(\frac{u_1 r_2w_3}{P}, \frac{u_2r_1w_3}{P}\right)
 \end{split}
 \end{displaymath}   
and similarly, 
$$
\Upsilon^{(2)}(\textbf{u}, \textbf{w})  =  \sum_{ r_1, r_3 = 1}^\infty
\sum_{\epsilon \in \{\pm\}}  q\left(\frac{r_1}{u_1}\right)    
f_{0}\left(\frac{u_2r_3w_1}{P}\right)
k^{\epsilon}_{0}\left(\frac{u_3r_1w_2}{P}, \frac{u_1r_3w_2}{P}\right)
f_{0}\left(\frac{ u_2r_1w_3}{P}\right).
$$
 We now smooth the sums $\Upsilon^{(j)}(\textbf{u}, \textbf{w})$ by replacing
$f_0$ with $f_\Delta$ and $k^\pm_0$ with $k^\pm_\Delta$ where appropriate. 
Thus, we define the sums
\begin{equation}\label{b}
\begin{split}
    \Upsilon_{\Delta}^{(1)}(\textbf{u}, \textbf{w}) & = \sum_{\mathbf r \in\Bbb N^3} \sum_{{\bm \epsilon} \in E} q\left(\frac{r_1}{u_1}\right)    k^{\epsilon_1}_{\Delta}\left(\frac{u_2r_3w_1}{P}, \frac{u_3r_2w_1}{P}\right)  k^{\epsilon_2}_{\Delta}\left(\frac{u_3r_1w_2}{P}, \frac{u_1r_3w_2}{P}\right)k^{\epsilon_3}_{\Delta}\left(\frac{u_1 r_2w_3}{P}, \frac{u_2r_1w_3}{P}\right),\\
        \Upsilon_{\Delta}^{(2)}(\textbf{u}, \textbf{w}) & =  \sum_{ r_1, r_3 = 1}^\infty \sum_{\epsilon \in \{\pm\}}  q\left(\frac{r_1}{u_1}\right)    f_{\Delta}\left(\frac{u_2r_3w_1}{P}\right)k^{\epsilon}_{\Delta}\left(\frac{u_3r_1w_2}{P}, \frac{u_1r_3w_2}{P}\right)f_{\Delta}\left(\frac{ u_2r_1w_3}{P}\right),\\
             \Upsilon_{\Delta}^{(3)}(\textbf{u}, \textbf{w}) & =  \sum_{ r_1\geq 1} q\left(\frac{r_1}{u_1}\right)     f_{\Delta}\left(\frac{u_3r_1 w_2}{P}\right)
 f_{\Delta}\left(\frac{u_2r_1 w_3}{P}\right)
       \end{split}
\end{equation}
and, in accordance with \rf{NEU}, 
\be{defVjD}
  V^{(j)}_{\Delta}(P ) = 48 \underset{(u, \textbf{u}, \textbf{w}) \in \Bbb{N}^7}{\left.\sum\right.^{\ast}} f_{\text{id}, \Delta, \delta}\left(\frac{uu_2u_3w_1}{P}, \frac{uu_1u_3w_2}{P}, \frac{uu_1u_2w_3}{P}\right) \Upsilon^{(j)}_{\Delta}(\textbf{u}, \textbf{w}).  
\ee
Then, on recalling \rf{suppf}, \rf{suppk+} and \eqref{suppk-} and reversing the above 
analysis, one finds that the sum
$ V^{(1)}_{\Delta}(P ) + 2V^{(2)}_{\Delta}(P ) + V^{(3)}_{\Delta}(P )$
differs from the first term on the right hand side of \rf{NEU} by at most
$V_*((1-10\Delta) P, [(1-10\Delta)P, (1+10\Delta)P])$. Here, we have applied the notation
used in Lemma \ref{lem3}, and this lemma now shows that
\begin{equation}\label{smooth}
  V(P ) = V^{(1)}_{\Delta}(P ) + 2V^{(2)}_{\Delta}(P ) + V^{(3)}_{\Delta}(P ) + O(\Delta P^{3+\varepsilon}).
\end{equation}
This completes the preparatory transformation of $V(P)$.

\subsection{Contour integration}\label{cont} 
The sums in \rf{b} and \rf{defVjD} are ready for treatment by Mellin inversion.
Shifts of  complex contour integrals will ultimately yield asymptotic 
formulae for $V^{(j)}_{\Delta}(P )$.  The hardest case to analyze 
will be $j=1$, and we consider this  one first. %For $\eta \in \Bbb{R}$ let
%\begin{displaymath}
 % \textbf{c}_{\eta} = \left(\frac{1}{3}+\eta, \frac{1}{3}+\eta, \frac{1}{3}+\eta; \frac{1}{3}+\eta, \eta, \frac{2}{3}+\eta, \frac{2}{3}+\eta, \eta, 1 + \eta\right) \in \Bbb{R}^{10}.
%\end{displaymath}
%For $\eta > 0$ 
Let  $\textbf{c} \in (0,\infty)^{10}$. Then, by \rf{b}, \rf{defVjD}  
and \rf{mellininv},  one finds that
\begin{displaymath}
\begin{split}
  V_{\Delta}^{(1)}(P ) =&  \underset{(u, \textbf{u}, \textbf{w}, \textbf{r}) \in \Bbb{N}^{10}}{\left.\sum\right.^{\ast}}   \sum_{{\bm \epsilon} \in E} \frac{48}{(2\pi i)^{10}} \int_{(\textbf{c})} \frac{P^{s_1+s_2+s_3}\widehat{f}_{\text{id}, \Delta, \delta}(s_1, s_2, s_3)}{w_1^{s_1}w_2^{s_2}w_3^{s_3}u_1^{s_2+s_3}u_2^{s_1+s_3} u_3^{s_1+s_2}u^{s_1+s_2+s_3}}\\
  & \quad \frac{P^{s_4+s_5}\widehat{k}_{\Delta}^{\epsilon_1}(s_4, s_5) }{(r_3u_2w_1)^{s_4}(r_2u_3w_1)^{s_5}}\, \frac{P^{s_6+s_7}\widehat{k}^{\epsilon_2}_{\Delta}(s_6, s_7) }{(r_1u_3w_2)^{s_6}(r_3u_1w_2)^{s_7}}\,  \frac{P^{s_8+s_9}\widehat{k}^{\epsilon_3}_{\Delta}(s_8, s_9)}{(r_2u_1w_3)^{s_8}(r_1u_2w_3)^{s_9}}\,  \frac{\widehat{q}(s_{10}) u_1^{s_{10}}}{r_1^{s_{10}}}\, d\textbf{s}.  
\end{split}  
\end{displaymath}
Note that the integral on the right is absolutely convergent. Recall also that 
$\sum^{\ast}$ denotes the coprimality conditions \eqref{2.2} and \eqref{2.5}. 

The last identity may be rewritten in more balanced form. With this end in 
view, we define the function $ \Phi_{\Delta} (\textbf{s})$ as the sum 
\begin{equation}\label{defphi}
\begin{split}
   \sum_{{\bm \epsilon} \in E}
\textstyle  \widehat{f}_{\text{id}, \Delta, \delta}\left(\frac{1}{3} + s_1, \frac{1}{3} + s_2, \frac{1}{3} + s_3\right)\widehat{k}_{\Delta}^{\epsilon_1}\left(\frac{1}{3} + s_4, \frac{1}{3} + s_5\right) 
     \widehat{k}^{\epsilon_2}_{\Delta}\left(s_6, \frac{2}{3}+ s_7\right)   \widehat{k}^{\epsilon_3}_{\Delta}\left(\frac{2}{3} + s_8, s_9\right)  \widehat{q}(1+s_{10}). 
\end{split}    
\end{equation}
The dependence of $\Phi_{\Delta}$ on $\delta$ is suppressed, as on earlier occasions. By Lemma \ref{lem10}(i), (iii) and Lemma \ref{lem11}(iii),
the function $s_6s_9 \Phi_{\Delta}(\textbf{s})$   is holomorphic in 
$\Re s_j > - 1/3$. 
We introduce the linear forms $\ell_j=\ell_j(\mathbf s)$ by
\begin{displaymath}
\begin{split}
& \ell_1 = s_6, \quad \ell_2 = s_9, \quad
  \ell_3= s_6+s_9+s_{10}, \quad \ell_4 = s_5+s_8, \quad \ell_{5} = s_4+s_7,  \\
& 
\ell_6  = s_1+s_4+s_5,  
\quad  \ell_7 = s_2+s_6+s_7, \quad \ell_8 = s_3 + s_8+s_9,\quad
\ell_{9} = s_2+s_3 + s_7+s_8-s_{10}, \\
&  \ell_{10} = s_1+s_3 + s_4 + s_9, \quad \ell_{11} = s_1+s_2 + s_5 + s_6,\quad \ell_{12} = s_1+s_2+s_3, 
 \end{split}
\end{displaymath}
and may then recast the previous expression for $ V_{\Delta}^{(1)}(P )$ in the form
\begin{displaymath}
\begin{split}
  V_{\Delta}^{(1)}(P ) = 
\underset{(u, \textbf{u}, \textbf{w}, \textbf{r}) \in \Bbb{N}^{10}}{\left.\sum\right.^{\ast}}     \frac{48}{(2\pi i)^{10}} \int_{(\textbf{c})} 
\frac{P^{3+ s_1+\ldots + s_9}    s_6s_9\Phi_{\Delta}(\textbf{s}) }
{\ell_1\ell_2 
   r_1^{1+\ell_3}r_2^{1+\ell_4}r_3^{1+\ell_5}w_1^{1+\ell_6}w_2^{1+\ell_7}w_3^{1+\ell_8}u_1^{1+\ell_9}u_2^{1+\ell_{10}}u_3^{1+\ell_{11}}u_{}^{1+\ell_{12}}}\,  d\textbf{s}.
 \end{split}  
\end{displaymath}
Here we may sum inside the integral and are then
in a position to apply the theory developed in Chapter \ref{Diri}. Let $G$ be the graph defined in \eqref{graph}, and define
\begin{equation}\label{defZ}
  Z(s) = \zeta(1+s)s,
\end{equation}
\begin{equation}\label{x}
  \Theta_{\Delta}(\textbf{s}) = s_6s_9\Phi_{\Delta}(\textbf{s})\Xi_G(1+\ell_6(\textbf{s}), 1+\ell_7(\textbf{s}), \ldots, 1+\ell_{11}(\textbf{s})) \prod_{j=3}^{12}Z(\ell_j(\textbf{s})). 
  \end{equation}
Note that 
 $\Theta_{\Delta}(\textbf{s})$ is  holomorphic in $|\Re s_j| < 1/10$ by Theorem \ref{thm3}. 
Hence, when $0<c_j<1/10$, 
the previous expression for $V_{\Delta}^{(1)}(P )$ may be rewritten in the form
$$
  V_{\Delta}^{(1)}(P ) =  \frac{48}{(2\pi i)^{10}} \int_{(\textbf{c})} \frac{P^{3+ s_1+\ldots + s_9}    \Theta_{\Delta}(\textbf{s}) }{\ell_1(\textbf{s})\cdot \ldots \cdot \ell_{12}(\textbf{s}) } d\textbf{s}. 
$$

      Although a more careful analysis is needed later, for the moment we 
content ourselves with recording the crude bound
  \begin{equation}\label{crude}
     \Theta_{\Delta}(\textbf{s}) \ll  \Delta^{-99}\prod_{j=1}^{10}(1+|s_j|)^{-2}
  \end{equation}
that is available uniformly in  $ |\Re s_j| \le 1/20 $. This  follows 
from Lemma \ref{lem10} and the convexity bound for Riemann's zeta function
in the critical strip,  and ensures absolute convergence of all 
integrals that occur in the following discussion. 
  
  We observe that the   linear forms $\ell_j$, $1 \leq j \leq 12$,  
span an 8-dimensional vector space. 
%\footnote{This is   a precursor that the degree of the polynomial in Theorem \ref{thm1} is $12-8 = 4$.} 
   One checks that the linear forms $\ell_1, \ldots, \ell_8$ are linearly independent, while
  \begin{displaymath}
  \begin{split}
   & \ell_9 = \ell_7 + \ell_8 - \ell_3, \quad \ell_{10} = \ell_6 + \ell_8 - \ell_4, \quad \ell_{11} = \ell_6 + \ell_7 - \ell_5,\quad  \ell_{12} = \ell_6 + \ell_7 + \ell_8 - \ell_1 - \ell_2 - \ell_4 - \ell_5.
  \end{split}  
  \end{displaymath}
 We prepare a linear change of variables $A\textbf{z} =  \textbf{s}$ where $A \in GL_{10}(\Bbb{R})$ is such that
 \begin{displaymath}
\begin{split}
   z_j & = \ell_j(\textbf{s}), \quad 1 \leq j \leq 7, \\
   z_8 & = \ell_8(\textbf{s}) + \ell_7(\textbf{s}) + \ell_6(\textbf{s}) = s_1 + \ldots + s_9,\\
   z_9 & = s_1, \quad z_{10} = s_2.
  \end{split} 
\end{displaymath}
The definition of $z_9$ and $z_{10}$ is fairly arbitrary, and these variables will play no role in the following computations.  A straightforward computation shows that  $\det A =1$. 

Now consider the
 additional linear forms $\lm_j= \lm_j(\mathbf z)$ defined via
\begin{displaymath}
\lambda_1 =  z_6+z_7 - z_5,  \quad \lambda_2 = z_6 + z_7, \quad  
\lambda_3 = z_6 + z_3, \quad   \lambda_4 = z_7 + z_4,  \quad  
\lambda_5 = z_1+z_2+z_4+z_5. 
\end{displaymath}
We are ready for the change of variable. Let $\eta_1 = 10^{-6}$ and put
$c_8 = 5 \eta_1$, $c_j=\eta_1$ for $1\le j \le 10$, $j\neq 8$. Then $A\mathbf c$
has non-negative coordinates only, as one readily checks, and after a 
modest computation we  conclude that 
\begin{equation}\label{equation}
 V_{\Delta}^{(1)}(P )=  \frac{48}{(2\pi i)^{10}} \int_{(\textbf{c})}  
\frac{P^{3+z_8}\Theta_{\Delta}(A\textbf{z})  d\textbf{z}}{z_1   \cdots    z_7   \lambda_1(\textbf{z})\cdot  \prod_{\nu=2}^5 (z_8 - \lambda_{\nu}(\textbf{z})) %(u_8 - \lambda_2(\textbf{u}))(u_8 -   \lambda_3(\textbf{u})) (u_8- \lambda_5(\textbf{u}) ) 
  } .  
\end{equation}
For later purposes we unfold the rather compact notation and
spell out the integrand in more detail: as one readily checks, one  has
\be{explicitres} \frac{ \Theta_{\Delta}(A\textbf{z})  }{z_1  \cdots   z_7\lambda_1(\textbf{z}) \prod_{\nu=2}^5 (z_8 - \lambda_{\nu}(\textbf{z}))  }\hspace{10em}
\ee \vspace{-3mm}
\begin{displaymath}
\begin{split}
   = &\sum_{{\bm \epsilon} \in E} \widehat{f}_{\text{id}, \Delta, \delta}
\textstyle   \left(\frac{1}{3} + z_9, \frac{1}{3} + z_{10}, \frac{1}{3} + z_8 - z_1 -z_2 - z_4 - z_5 - z_9 - z_{10}\right) \\ & \textstyle
   \times   \widehat{k}^{\epsilon_2}_{\Delta}\left(z_1, \frac{2}{3}+ z_7 - z_1 - z_{10}\right)   \widehat{k}^{\epsilon_3}_{\Delta}\left(\frac{1}{3} + z_1 + z_4 + z_5 +  z_9 + z_{10} - z_6 - z_7, z_2\right) \\ \textstyle
& \times   \widehat{k}_{\Delta}^{\epsilon_1}\textstyle\left(\frac{1}{3} + z_1 + z_{10} + z_5 - z_{7}, \frac{1}{3} + z_6 + z_7 - z_1 - z_5 -z_9 - z_{10}\right)  \widehat{q}(1+z_3 - z_2 - z_1) H(\textbf{z}) \end{split}  
\end{displaymath} 
where   $H$ is a function independent of $\Delta$ that satisfies
\begin{equation}\label{boundH}
  H(\textbf{z}) \ll \prod_{j=1}^{10} (1+|z_j|)^{1/1000}
\end{equation}
uniformly in the region $|\Re z_j| \leq 10\eta_1$, $|z_j|\ge 1$. This again follows from 
standard bounds for Riemann's zeta function close to the line $\Re s = 1$.

\subsection{Multiple contour shifts}\label{mult}
To extract a main term from  \eqref{equation}, we apply the common routine
and shift  the lines of integration of $z_1, \ldots, z_8$ to the left. We will 
do this step by step. The factor $P^{3+z_8}$ in \rf{equation} will not be affected until the last step when the line for $z_8$ will be moved. 
The following notation will be helpful in describing 
the manoevre. Let $\textbf{z} \in \Bbb{C}^{10}$ and 
$I \subset \{1, \ldots, 7\}$. Integration over all variables $u_i$ with $i \not\in I$ will be denoted by  $d\textbf{z}_{I}$. We also denote by
 $\textbf{z}_{I} \in \Bbb{C}^{10}$ the 10-tuple whose $i$-th entry is $u_i$ if $i \not \in I$, and equal to 0 if $i \in I$. Furthermore, let 
\begin{displaymath}
\langle \textbf{z}_I \rangle =  \prod_{\substack{1 \leq i \leq 7\\ i \not\in I}} z_i. 
 \end{displaymath}

%\smallskip

 We begin by   shifting the contours of the 
variables $z_1, \ldots, z_5$ to the left up to $\Re z_j = -7\eta_1$,
one after the other. The  shift of $\Re z_i = \eta_1$ to $\Re z_i = -7\eta_1$ 
 passes through exactly one pole at $z_i = 0$, and this pole is simple. 
In order to describe the outcome of the residue theorem, we  define  
for  $I \subset \{1, \ldots, 5\}$ a vector 
$\textbf{c}^{I} \in \Bbb{R}^{10 - |I|}$ by $c_i = \eta_1$ if 
$i \in \{6, 7, 9, 10\}$, $c_8 = 5\eta_1$ and $c_i = -7\eta_1$ if $ i \leq 5$ and $i \not \in I$. Then, by \eqref{equation}, we infer that  
\begin{displaymath}
V_{\Delta}^{(1)}(P ) =
  \sum_{I \subset \{1, \ldots, 5\}} \frac{48}{(2\pi i)^{10-|I|}} \int_{(\textbf{c}^{I})}  \frac{ P^{3+z_8}{\Theta}_{\Delta}(A\textbf{z}_{I}) \, d\textbf{z}_I }{ \langle \textbf{z}_I \rangle  \lambda_1(\textbf{z}_I) \prod_{\nu=2}^5 (z_8 - \lambda_{\nu}(\textbf{z}_I)) }.
\end{displaymath}

 Next we shift the variables $z_6$ and $z_7$ to $\Re z_i = - 3\eta_1$. 
First consider the case where $5 \not\in I$. Then $\lambda_1(\textbf{z}_I) = z_6+z_7 - z_5$.  Each of the two contour shifts passes through 
exactly one  simple pole at $z_j = 0$. Much as before, 
for a subset $J \subset \{1, 2, 3, 4, 6, 7\}$ let 
$\tilde{\textbf{c}}^J \in \Bbb{R}^{10 - |J|}$ be defined by 
$\tilde{c}_j = -7\eta_1$ for $j \leq 5$, $j \not\in J$, by $\tilde{c}_j = -3\eta_1$ for $j \in \{6, 7\}$,  $j \not\in J $, and finally by  
$\tilde{c}_8 = 5\eta_1$ and $\tilde{c}_j = \eta_1$ for $j \in \{9, 10\}$. Then the contribution of the terms with $5 \not\in I$  can be written as
\begin{equation}\label{simple}
\begin{split}
&  \sum_{J \subset \{1, 2, 3, 4, 6, 7\}} \frac{48}{(2\pi i)^{10-|J|}} \int_{(\tilde{\textbf{c}}^J)}  \frac{P^{3+z_8}{\Theta}_{\Delta}(A\textbf{z}_J)\,  d\textbf{z}_J }{
 \langle \textbf{z}_J \rangle \lambda_1(\textbf{z}_J) \prod_{\nu=2}^5    (z_8 - \lambda_{\nu}(\textbf{z}_J))  }   .  
\end{split}  
\end{equation}
 The case $5 \in I$ that we now consider,  is different. Here 
$\lambda_1(\textbf{z}_I) = z_6+z_7$. Hence, when we shift 
$z_6$ to  $\Re z_6 = - 3\eta_1$, we
pass through two simple poles at $z_6 = 0$ and $z_6 = -z_7$.
This generates three terms: the residue $\mathcal{R}_1$
from the pole at $z_6=0$, the residue $\mathcal{R}_2$
from the pole at $z_6=-z_7$, and a  remaining integral
 $\mathcal{I}$  over the line $\Re z_6 = - 3\eta_1$. In $\mathcal{R}_1$ 
the variable $z_6$ is no longer active, and hence, $\lambda_1 = z_7$. 
Consequently,
shifting $z_7$ to $\Re z_7 =-3\eta_1$ in $\mathcal{R}_1$   
passes through a double pole at $z_7 = 0$ which contributes the residue
\begin{equation}\label{harder1}
\begin{split}
   \sum_{\{5, 6, 7\} \subset  J  \subset \{1, \ldots, 7\}} \frac{48}{(2\pi i)^{10-|J|}}& \int_{(\tilde{\textbf{c}}^{J})}  \frac{P^{3+z_8} }{ \langle \textbf{z}_J \rangle(z_8 - \lambda_3(\textbf{z}_J))(z_8-\lambda_5(\textbf{z}_J))} \\ 
   &   \frac{\partial}{\partial z_7} \left(\frac{{\Theta}_{\Delta}(A\textbf{z}_{\tilde{J}})}{(z_8 - \lambda_2(\textbf{z}_{\tilde{J}}))(z_8 - \lambda_4(\textbf{z}_{\tilde{J}}))}\right)\Bigl|_{z_7 = 0}   d\textbf{z}_J. 
\end{split}  
\end{equation}
Here we have written $\tilde{J} = J \setminus \{7\}$ in the second line to reactivate the variable $z_7$ for differentiation. After the shift, an integral remains that, together with the contribution of $\mathcal{I}$, 
produces a term identical to \eqref{simple}, but where $J \subset \{1, \ldots 7\}$ is subject to the conditions $5, 6 \in J$, $7 \not \in J$ or $5 \in J$, $6 \not\in J$. Together with \eqref{simple}, these terms combine to 
\begin{equation}\label{harder2}
\begin{split}
&   \sum_{\substack{J \subset \{1, \ldots, 7\}\\ \{5, 6, 7\} \not\subset J}} \frac{48}{(2\pi i)^{10-|J|}} \int_{(\tilde{\textbf{c}}^J)}  \frac{P^{3+z_8}{\Theta}_{\Delta}(A\textbf{z}_J) \,   d\textbf{z}_J}{
  \langle \textbf{z}_J \rangle \lambda_1(\textbf{z}_J) \prod_{\nu=2}^5 (z_8 - \lambda_{\nu}(\textbf{z}_J))  } .  
\end{split}  
\end{equation}

The term $\mathcal{R}_2$ is more complicated. When we shift $z_7$, we pass 
through a double pole at $z_7=0$ that contributes the residue  
\begin{equation}\label{harder3}
\begin{split}
 -  \sum_{\substack{ 
 J  \subset \{1, \ldots, 7\}\\ \{5, 6, 7\} \subset J}} \frac{48}{(2\pi i)^{10-|J|}}& \int_{(\tilde{\textbf{c}}^{J})}  \frac{P^{3+z_8} }{\langle \textbf{z}_J \rangle z_8  (z_8-\lambda_5(\textbf{z}_J))}    \frac{\partial}{\partial z_7} \left(\frac{{\Theta}_{\Delta}(A\textbf{z}_{\tilde{J}})}{(z_8  -  \tilde{\lambda}_3(\textbf{z}_{\tilde{J}}))(z_8 - \lambda_4(\textbf{z}_{\tilde{J}}))}\right)\Bigl|_{z_7 = 0}   d\textbf{z}_J. 
\end{split}  
\end{equation}
Here we have written $\tilde{\lambda}_3(\textbf{z}) = -z_7+z_3$ and, as before, $\tilde{J} = J \setminus \{7\}$.

Moreover, there may be  an additional simple pole at $z_7 = -z_8$, and this is so if and only if $3 \in I$  
(that is, $\lambda_2(\textbf{z}_I) = z_6$). The residue contributes
\begin{equation}\label{harder4}
\begin{split}
  & -   \sum_{\{3, 5, 6, 7\} \subset J \subset \{1, \ldots, 7\}} \frac{48}{(2\pi i)^{10-|J|}} \int_{(\tilde{\textbf{c}}^J)}  \frac{P^{3+z_8}{\Theta}_{\Delta}(A\textbf{z}_J)}{z_8^3(2z_8 - \lambda_4(\textbf{z}_J))(z_8 - \lambda_5(\textbf{z}_j))}  d\textbf{z}_J .  
\end{split}  
\end{equation}
After the shift to $\Re z_7 = -3\eta_1$, remaining integral contributes
\begin{equation}\label{harder5}
\begin{split}
& -   \sum_{\substack { J \subset \{1, \ldots, 7\}\\ 5, 6 \in J, \, 7 \not\in J}} \frac{48}{(2\pi i)^{10-|J|}} \int_{(\tilde{\textbf{c}}^J)}  \frac{P^{3+z_8}{\Theta}_{\Delta}(A\textbf{z}_J)\,  d\textbf{z}_J }{\langle \textbf{z}_J \rangle   z_7z_8 (z_8 - \tilde{\lambda}_3(\textbf{z}_J))  (z_8-   \lambda_4(\textbf{z}_J))  (z_8- \lambda_5(\textbf{z}_J) ) }     
\end{split}  
\end{equation}
in which again $\tilde{\lambda}_3(\textbf{z}) = -z_7+z_3$. 

To summarize the above analysis, it is mandatory to record here that
$V^{(1)}_{\Delta}(P )$ is 
the sum of the five terms \eqref{harder1} -- \eqref{harder5}. In each 
term, we  finally shift the line of integration 
$\Re z_8= \eta_1$ to $\Re z_8 = -\eta_1$. 
The only pole that may occur is at  $z_8 = 0$. Its order depends on the particular set $J$ in the various sums, but it is immediate that no pole of order
higher than five does occur. Consequently, the residues of these poles will
combine to  $P^3  Q^{(1)}_{\Delta}(\log P)$
where $Q_{\Delta}$ is a polynomial depending on $\Delta$ 
of degree at most four. Also, the integral over $\Re z_8 = -\eta_1$ that remains after the shift may be estimated by \eqref{crude}, and it transpires that its contribution does not exceed
$ O(P^{3-\eta_1} \Delta^{-99})$. Thus, we have now shown that
\begin{equation}\label{error}
V^{(1)}_{\Delta}(P ) = P^3  Q^{(1)}_{\Delta}(\log P)+ O(P^{3-\eta_1} \Delta^{-99}).
\end{equation}

This expansion has an unfortunate defect: the polynomial $Q^{(1)}_{\Delta}$
has coefficients depending on $\Delta$, and this should not be the case. To
rectify this, one may replace ${\Theta}_{\Delta}(A\textbf{u}_J)$ by 
${\Theta}_0(A\textbf{u}_J)$ in the computation of all residues at $z_8=0$
that contribute to the coefficients of $Q^{(1)}_{\Delta}$. One then obtains a polynomial $Q^{(1)}_0$, say, that is independent of $\Delta$, yet one then has to control the error that is introduced by this procedure. We postpone a detailed discussion of this matter to the next chapter, and continue with the derivation of formulae of the type \rf{error} for  $V^{(2)}_{\Delta}(P )$ and $ V^{(3)}_{\Delta}(P )$.

\subsection{Multiple contour shifts II}

Having discussed $V_{\Delta}^{(1)}(P )$ in detail, it is now an easy exercise to analyze the less complicated terms $V^{(2)}_{\Delta}( P)$ and $V^{(3)}_{\Delta}(P )$. We can be brief here. 

We start with the analysis of $V^{(2)}_{\Delta}( P )$ and recall  \eqref{b}.  Following the argument at the beginning of Section \ref{cont}, we define
\begin{displaymath} 
  \Phi_{\Delta}(\textbf{s}) =  \sum_{\epsilon=\pm}
\textstyle \widehat{f}_{\text{id}, \Delta, \delta}\left(\frac{1}{3} + s_1, \frac{1}{3} + s_2, \frac{1}{3} + s_3\right)\widehat{f}_{\Delta}\left(\frac{2}{3} + s_4\right) \widehat{k}^{\epsilon}_{\Delta}\left(\frac{1}{3} + s_5, \frac{1}{3}+ s_6\right)   \widehat{f}_{\Delta}\left(\frac{2}{3} + s_7\right)  \widehat{q}(s_{8}),   
\end{displaymath}  
the eight linear forms
\begin{displaymath}
\begin{split}
& \ell_1 = s_1 + s_2 + s_3, \quad \ell_2 = s_2 + s_3 + s_6 - s_8,  
\quad \ell_3 = s_1 + s_2 + s_5,  \quad \ell_4 = s_1 + s_4,\\
 & \ell_5 = s_3 + s_7, \quad \ell_6 = s_2 + s_5 + s_6, \quad   
\ell_7 = s_5 + s_7 + s_8, \quad \ell_8 = s_4 + s_6 
\end{split}
\end{displaymath}
and the linear polynomial $l(\textbf{s}) = s_1 + s_3 + s_4 + s_7 + 2$.
Note that this overwrites the definitions of $\Phi$ and $\ell_j$
used in the preceding section; confusion should not arise.  
%and the 4 linear polynomials
%\begin{displaymath}
%\begin{split}
%& l_1(\textbf{s}) = \frac{9}{5} + s_1 + s_2 + s_3, \quad l_2(\textbf{s}) = \frac{11}{10} + s_2 + s_3 + s_6-s_8, \\
% & l_3(\textbf{s}) = \frac{9}{5} + s_1 + s_3 + s_4 + s_7, \quad l_4(\textbf{s}) = \frac{11}{10} + s_5 + s_7 + s_8.
%\end{split} 
%\end{displaymath}
Then, for any $\textbf{c} \in (0,1/4)^{8}$, we have
\begin{displaymath}
\begin{split}
  V^{(2)}_{\Delta}(P ) = & \underset{\substack{\textbf{u}, \textbf{w} \in \Bbb{N}^3\\ u, r_1, r_3 \in \Bbb{N}}}{\left.\sum\right.^{\ast}}     \frac{48}{(2\pi i)^8} \int_{(\textbf{c})} \frac{P^{3 + s_1+\ldots + s_7}\Phi_{\Delta}(\textbf{s}) \, d \textbf{s}}{u^{1+\ell_1(\textbf{s})}u_1^{1+ \ell_2(\textbf{s})}u_2^{ l(\textbf{s})}u_3^{1+\ell_3(\textbf{s})}w_1^{1+\ell_4(\textbf{s})}w_2^{1+\ell_6(\textbf{s})} w_3^{1+\ell_5(\textbf{s})}r_1^{1+ \ell_7(\textbf{s})}r_3^{1+\ell_8(\textbf{s})}}
  \end{split}
\end{displaymath}
where again $\left.\sum\right.^{\ast}$ denotes a summation subject to the coprimality conditions \eqref{2.2} and \eqref{2.5}. We pull the multiple sum inside the integral and obtain the identity
\begin{displaymath}
\begin{split}
  V^{(2)}_{\Delta}(P ) =    \frac{48}{(2\pi i)^8}  \int_{(\textbf{c})} \frac{P^{3 + s_1+\ldots + s_7}\Theta_{\Delta}(\textbf{s}) }{\ell_1(\textbf{s}) \cdot \ldots \cdot \ell_8(\textbf{s}) } d \textbf{s}
  \end{split}
\end{displaymath}
where
\begin{displaymath}
  \Theta_{\Delta}(\textbf{s}) = \Phi_{\Delta}(\textbf{s})\zeta(l(\textbf{s}) ) \Xi_G\bigl(1+ \ell_2(\textbf{s}), l(\textbf{s}), 1+\ell_3(\textbf{s}), \ldots, 1+\ell_6(\textbf{s})\bigr)\prod_{j=1}^8 Z(\ell_j(\textbf{s})),  \end{displaymath}
and where $Z$ is still defined by \eqref{defZ}. Lemma \ref{lem10} shows that 
  $\Theta_{\Delta}$ is holomorphic in $|\Re s_j|\le 1/4$ where the bound  
  \begin{displaymath}
     \Theta_{\Delta}(\textbf{s}) \ll  \Delta^{-99}\prod_{j=1}^{8} (1+|s_j|)^{-2}
  \end{displaymath}
holds. This follows as in \eqref{crude}. 
The eight linear forms $\ell_1(\textbf{s}), \ldots, \ell_8(\textbf{s})$ 
span a space of dimension 6, and we define the matrix $A \in GL_8(\Bbb{R})$ by 
 $A\textbf{z} = \textbf{s}$ and
\begin{displaymath}
  z_j  = \ell_j(\textbf{s}) \quad (1 \leq j \leq 5), 
\quad z_6 = s_1 + \ldots + s_7, \quad  z_7 =s_2,\quad z_8 =s_1.
  \end{displaymath}
One checks that $\det A = 1$. Now choose $\textbf{c} \in \Bbb{R}^8$ with
$c_6 = 3\eta_1$ and $c_j = \eta_1$ for $j \not= 6$. Then, $A\mathbf c$
has all entries positive, and a change of variable produces
 \begin{equation}\label{resV2}
   V^{(2)}_{\Delta}(P ) =   \frac{48}{(2\pi i)^8}  \int_{(\textbf{c})}  \frac{P^{3 + z_6}{\Theta}_{\Delta}(A\textbf{z}) }{z_1\cdots z_5(z_6 - z_5-z_4)(z_6 -z_4 - z_2)(z_6 - z_5 - z_3) } d\textbf{z}.
 \end{equation}
 %We proceed as before and define
 %\begin{displaymath}
  %\lambda_1(\textbf{u}) = u_4 + u_5, \quad  \lambda_2(\textbf{u}) = u_2 + u_4, \quad  \lambda_1(\textbf{u}) = u_3 + u_5.
 %\end{displaymath}
 %Then shifting contours we obtain with the same notation as above
 We are ready to shift contours to the left. The new lines of integration are
 $\Re z_j = -2\eta_1$ for $1\le j\le 5$, and $\Re z_6 = -\eta_1$.  Then,
 much as before, one obtains  
\begin{equation}\label{prelim1}
  V^{(2)}_{\Delta}( P) = P^3 Q^{(2)}_{\Delta}(P ) + O(P^{3 - \eta_1} \Delta^{-99})
\end{equation}
where $Q_\Delta^{(2)}$ is a polynomial depending on $\Delta$ of degree at most 2. \\

%\medskip

The analysis of $V^{(3)}(P )$ is along the same lines. We (re-)define
 \begin{displaymath} 
  \Phi_{\Delta}(\textbf{s}) =  \widehat{f}_{\text{id}, \Delta, \delta}
\textstyle\left(1 + s_1, \frac{1}{2} + s_2,\frac12 + s_3\right)
\widehat{f}_{\Delta}\left(\frac{1}{2} + s_4 \right) 
\widehat{f}_{\Delta}\left(\frac{1}{2} + s_5 \right)   
\widehat{q}\left(s_{6}\right),  
\end{displaymath}
  the five linear forms
\begin{displaymath}
\begin{split}
& \ell_1(\textbf{s}) = s_2 + s_3 - s_6, \quad \ell_2(\textbf{s}) = s_1,  \quad \ell_3(\textbf{s}) = s_2 + s_4, \quad  \ell_4(\textbf{s}) = s_3+s_5,  \quad \ell_5(\textbf{s}) = s_4 + s_5 +s_6
\end{split}
\end{displaymath}
and the three linear polynomials
\begin{displaymath}
   l_1(\textbf{s}) = 2 + s_1 + s_2 + s_3, \quad 
l_2(\textbf{s}) = 2+ s_1 + s_3 +s_5, \quad l_3(\textbf{s}) = 2+ s_1+s_2+s_4.
\end{displaymath}   
Then, for $\mathbf c \in(0,1/4)^6$, one has 
\begin{displaymath}
  V^{(3)}_{\Delta}(P ) = \underset{\substack{\textbf{u}, \textbf{w} \in \Bbb{N}^3\\ u, r_1 \in \Bbb{N}}}{\left.\sum\right.^{\ast}}  \frac{48}{(2\pi i)^{5}} \int_{(\textbf{c})} \frac{P^{3 + s_1 + \ldots + s_5}   \Phi_{\Delta}(\textbf{s}) \, d\textbf{s}}{u^{l_1(\textbf{s})}  u_1^{1+ \ell_1(\textbf{s})}u_2^{ l_2(\textbf{s})}u_3^{ l_3(\textbf{s})}w_1^{1+\ell_2(\textbf{s})}w_2^{1+\ell_3(\textbf{s})} w_3^{1+\ell_4(\textbf{s})}r_1^{1+ \ell_5(\textbf{s})} }
\end{displaymath}
In $|\Re s_j|\le 1/4$, a holomorphic function is defined by
\begin{displaymath}
  \Theta_{\Delta}(\textbf{s}) = \Phi_{\Delta}(\textbf{s}) 
\Xi_G\bigl(1+ \ell_1(\textbf{s}), l_2(\textbf{s}), l_3(\textbf{s}), 1+\ell_2(\textbf{s}),1+\ell_3(\textbf{s}), 1+\ell_4(\textbf{s})\bigr)\prod_{j=1}^5 Z(\ell_j(\textbf{s}))  \prod_{j=1}^3\zeta\bigl(l_j(\textbf{s}) \bigr),  
  \end{displaymath}
and the previous formula for $V^{(3)}_{\Delta}(P )$ can be rewritten as
\begin{displaymath}
\begin{split}
  V^{(3)}_{\Delta}(P ) =     \frac{48}{(2\pi i)^6}  \int_{(\textbf{c})} \frac{P^{3 + s_1+\ldots + s_5}\Theta_{\Delta}(\textbf{s}) }{\ell_1(\textbf{s}) \cdot \ldots \cdot \ell_5(\textbf{s}) } d \textbf{s}. 
  \end{split}
\end{displaymath}
The linear forms $\ell_1, \ldots, \ell_5$ span a space of dimension 4, and we define $A \in GL_6(\Bbb{R})$ by $A \textbf{z} = \textbf{s}$ and
\begin{displaymath}
  z_j = \ell_j(\textbf{s}) \quad (1 \leq j \leq 3),\quad  z_4 = s_1 + \ldots +  s_5,\quad z_5 = s_2, \quad z_6 = s_6.
\end{displaymath}
Then $\det A=1$. Choose $\textbf{c}\in \Bbb{R}^6$ with $c_4 = 3\eta_1$ and
$c_j = \eta_1$ for $j \not= 4$. Then $ A\mathbf c$ has positive coordinates 
only, and a change of variable yields
\begin{equation}\label{resV3}
   V^{(3)}_{\Delta}(P ) =     \frac{48}{(2\pi i)^6}  \int_{(\textbf{c})} \frac{P^{3+z_4} {\Theta}_{\Delta}(A\textbf{z})}{z_1z_2z_3(z_4 - z_3 - z_2)(z_4- z_2 - z_1)} d\textbf{z}.
\end{equation}
One may now move the lines of integration to $\Re z_j = -2\eta_1$ for $1\le j\le3$, and then to $\Re z_4 = -\eta_1$. 
An argument similar to the one used above now readily  provides the expansion 
 \begin{equation}\label{prelim2}
   V^{(3)}_{\Delta}( P ) = P^3 Q_{\Delta}^{(3)}(\log P) + O(P^{3-\eta_1}\Delta^{-99})
 \end{equation}
 in which $Q^{(3)}_{\Delta}$ is a polynomial depending on $\Delta$, 
of degree at most 1.  Combining  \eqref{smooth}, \eqref{error}, \eqref{prelim1} and \eqref{prelim2}, we have shown
\begin{equation}\label{haveshown}
  V(P ) = P^3 Q_{\Delta}(\log P) + O(\Delta P^{3+\varepsilon} + 
P^{3-\eta_1}\Delta^{-99})
\end{equation}
where 
\begin{equation}\label{decompQ}
  Q_{\Delta} = Q^{(1)}_{\Delta} + 2 Q^{(2)}_{\Delta} + Q^{(3)}_{\Delta}.
\end{equation}  

\section{Removing the smoothing parameter}\label{remove}

\subsection{A useful lemma}
The principal goal in this chapter is to remove the dependence on $\Delta$
in the leading term on the right hand side of \rf{haveshown}. The discussion
of this theme generates certain multiple integrals, and we begin with a
lemma that ensures the existence of these integrals.

\begin{lem}\label{lem13} Let $m, n \in \Bbb{N}$, and let $\nu= 1/(8n)$.
Let $\ell_1, \ldots, \ell_m \in \Bbb{R}[x_1, \ldots, x_n]$ be $m$ linear 
forms, and suppose that for any $1\le i\le n$  the coefficient of $x_i$
is non-zero in at least one of the linear forms $\ell_j$. Then the function
$$
  \prod_{j=1}^n (1+|x_j|)^{\nu-1} \prod_{j=1}^m (1+|\ell_j(\mathbf{x})|)^{-1/3} 
$$
 is integrable over ${\Bbb R}^n$.
\end{lem}

Note that the hypotheses on the linear forms $\ell_j$ cannot be relaxed.
For example, if all $\ell_j$ would be independent of $x_1$, then a divergent
 integral
over $x_1$ would factor off.

For a proof of Lemma \ref{lem13}, write
$$ {\cal L}(\mathbf x) =  \prod_{j=1}^m (1+|\ell_j(\mathbf{x})|). $$
For $1\le i\le n$, the integral
$$  I_i =  \int_{\Bbb{R}^n} {\cal L}(\mathbf x)^{-4/3} 
\multprod{j=1}{j\neq i}^n (1+|x_j|)^{(3-4n)/(4n-4)} \, d\mathbf{x} $$
certainly exists, because by hypothesis there is a linear form $\ell_l$
for which the substitution $x_i\mapsto \ell_l$ is non-singular, and 
the inequality $ {\cal L}(\mathbf x) \ge 1+|\ell_l(\mathbf x)|$ shows that
the integral 
$$  \int_{\Bbb{R}^n} (1+ |\ell_l|)^{-4/3}\multprod{j=1}{j\neq i}^n (1+|x_j|)^{(3-4n)/(4n-4)}\, d\ell_l dx_1 \cdots \widehat{dx_i}\cdots dx_n$$
(with integration against $x_i$ omitted) is a majorant. Now, by H\"older's
inequality, 
$$ \int_{\Bbb{R}^n} {\cal L}(\mathbf x)^{-1/3}\prod_{j=1}^n (1+|x_j|)^{\nu-1} \, d\mathbf x
 \le (I_1 I_2 \cdots I_n)^{1/(4n)}
\Bigl( \int_{\Bbb{R}^n} \prod_{j=1}^n (1+|x_j|)^{-1-1/(12n)} \, d\mathbf x\Bigr)^{3/4},$$
which demonstrates the lemma.

\subsection{The generic case} We now turn to the main task 
in this chapter, and  derive the desired asymptotic formula 
\rf{asympV} from \rf{haveshown}. It will be necessary  to compare
the polynomial $Q_{\Delta}(\log P)$ with one that is independent of $\Delta$.
We have alluded to this problem already in commentary following \rf{error},
and the strategy suggested there will now be worked out in detail, 
in separate sections for the portions $ V^{(j)}_{\Delta}(P )$ contributing to the leading term in \rf{haveshown}. The goal is to show that 
there exists a polynomial $Q_{0}$
 such that
\begin{equation}\label{claim}
P^3Q_{\Delta}(\log P)   =   
   P^3Q_{0}(\log P)  +  O(P^{3+\varepsilon} \Delta^{1/200}). 
\end{equation}
 Taking this for granted, it follows from \eqref{haveshown}   that 
\begin{equation}\label{asymp}
  V(P ) = P^3 Q_0(\log P) + O\bigl(P^{\varepsilon}(\Delta^{1/200}  P^3 + P^{3 - \eta_1} \Delta^{-99})\bigr).
\end{equation}
One may take $
  \Delta = P^{- \eta_1/100 }
$ to infer 
 \eqref{asympV}  with $\tau = \eta_1/20000$, as required.\\

%\medskip

We now return to \rf{error}, and examine the origins of the polynomial
$Q^{(1)}_{\Delta}$. Its coefficients may be computed from the residues at
$z_8=0$, for each summand in the sums
\eqref{harder1} -- \eqref{harder5}. Summands corresponding to a set $J$ for which the integrand has no pole at $u_8=0$ do not contribute to $Q^{(1)}_{\Delta}$ and 
may therefore be ignored. 

As an illustrative example, we now examine in full detail the sum \rf{harder2}.
In this sum, the integrand has a pole of order $\nu$ at $z_8=0$ if and only if exactly $\nu$ of the four linear forms $\lm_2(\mathbf z _J), \ldots,
\lm_5(\mathbf z _J)$ vanish identically. But $\lm_l(\mathbf z _J)$ vanishes identically if and only if the condition $(l)$ in the list
$$  \text{(II)} \;\;\; \{6,7\} \subset J, \quad
\text{(III)} \;\;\; \{3,6\} \subset J, \quad
\text{(IV)} \;\;\; \{4,7\} \subset J, \quad
 \text{(V)} \;\;\; \{1,2,4,5\} \subset J 
$$
holds. The condition that $\{5,6,7\} \not\subset J$ implies that (II) and (V)
cannot hold simultaneously, so that the maximum order of the pole is $\nu=3$.

We begin with the summand $J=\{3,6\}$ in \rf{harder2}, and write $\hat J = J\cup
\{8\}$. In this case, the integrand in \rf{harder2} has a simple pole at
$z_8=0$, the residue of which is
\be{coeff1} P^3 \int_{(\tilde{\textbf{c}}^{\hat J})}  \frac{{\Theta}_{\Delta}(A\textbf{z}_{\hat J}) \,   d\textbf{z}_{\hat J}}{
  \langle \textbf{z}_J \rangle \lambda_1(\textbf{z}_J) \prod_{\nu=2}^5 (z_8 - \lambda_{\nu}(\textbf{z}_J))  }.
\ee
Thus, it transpires that the integral here contributes to the constant 
coefficient in $Q^{(1)}_{\Delta}$. As suggested earlier, one replaces 
$\Theta_\Delta$ by $\Theta_0$ and estimates the error. In the interest of brevity, we write 
$$ \Psi_\Delta(\textbf{z}_{\hat J}) = \frac{{\Theta}_{\Delta}(A\textbf{z}_{\hat J}) \,   }{
  \langle \textbf{z}_J \rangle \lambda_1(\textbf{z}_J) \prod_{\nu=2}^5 (z_8 - \lambda_{\nu}(\textbf{z}_J))  }
$$
and take $z_3=z_6=z_8=0$ in \rf{explicitres} to infer the alternative representation
\begin{equation}\label{Psi1}
  \begin{split} 
\Psi_\Delta &(\textbf{z}_{\hat J}) =  \textstyle \widehat{f}_{\text{id}, \Delta, \delta}\left(\frac{1}{3} + z_9, \frac{1}{3} + z_{10}, \frac{1}{3}  - z_1 -z_2 - z_4 - z_5 - z_9 - z_{10}\right)\\
   &\textstyle \widehat{k}_{\Delta}^{\epsilon_1}\left(\frac{1}{3} + z_1 + z_{10} + z_5 - z_{7}, \frac{1}{3} +  z_7 - z_1 - z_5 - z_9 - z_{10}\right) \widehat{k}^{\epsilon_2}_{\Delta}\left(z_1, \frac{2}{3}+ z_7 - z_1 - z_{10}\right)   \\
    &\textstyle  \widehat{k}^{\epsilon_3}_{\Delta}\left(\frac{1}{3} + z_1 + z_4 + z_5 +  u_9 + z_{10}  - z_7, z_2\right) \widehat{q}(1 - z_1 - z_2) H(\textbf{z}_{\hat J}).
 \end{split}  
\end{equation}
Here $H$ is independent of $\Delta$ and $\delta$. We may formally take 
$\Delta=0$ in \rf{Psi1} to define $\Psi_0$. Then, the error that arises
from replacing $\Delta$ by 0 in \rf{coeff1} is given by the integral
\be{star}
 P^3 \int_{(\tilde{\textbf{c}}^{\hat J})} \big(\Psi_\Delta(\textbf{z}_{\hat J})
- \Psi_0(\textbf{z}_{\hat J})\big)\, d\mathbf z_{\hat J}
\ee

To estimate this integral, we parametrize the lines of integration via
$z_j=c_j+it_j$, with $\mathbf t_{\hat J}\in\Bbb R^7$. Then, by \rf{Psi1}, 
Lemma \ref{lem10} (i), Lemma \ref{lem12} (with $\alpha=-50\eta_1$) and \rf{boundH}, for any $0<\eta<1/2$
one finds that
\be{star1} \Psi_\Delta(\textbf{u}_{\hat J})
- \Psi_0(\textbf{z}_{\hat J}) \ll \Delta^\eta {\cal I}(\mathbf t_{\hat J})
\ee
where
\begin{displaymath}\begin{split}
{\cal I}(\mathbf t_{\hat J}) = & \bigl((1 + |z_9|)(1+| z_{10}|)(1+|z_1+z_2+z_4+z_5+z_9 + z_{10}|)\bigr)^{\eta - 1} \\
& (1+|z_1  +z_5 + z_{10}- z_7|)^{\eta - 1} (1 + |z_7 - z_1 - z_5- z_9 - z_{10}|)^{\eta- 1/2}\\
 & (1+|z_2|)^{\eta-1}(1+|z_1+z_4+z_5 + z_9+z_{10}- z_7|)^{\eta-1/2} \\
 &(1+|z_7 - z_1-z_{10}|)^{\eta - 1} (1+|z_1|)^{\eta-1/2} (1+| z_2 + z_1|)^{-2}
 \prod_{j\in\{1, 2, 4, 5, 7, 9, 10\}} (1+|z_j|)^{1/1000}. 
 \end{split}  
\end{displaymath}
It remains to check that for $\eta=1/200$, the function 
${\cal I}(\mathbf t_{\hat J})$ is integrable over $\Bbb R^7$. Once this is
established, it follows from \rf{star1} that the expression in \rf{star}
is bounded by $O(P^3 \Delta^\eta)$, as is required for the verification
of \rf{claim}. 

The linear change of variable
\begin{displaymath}
\begin{split}
& x_1=t_9, \quad x_2=t_{10},\quad x_3 = t_1+t_2+t_4+t_5+t_9+t_{10},
\quad x_4= t_1+t_5+t_{10}-t_7\\
& x_5=t_2, \quad x_6=t_7-t_1-t_{10}, \quad x_7= t_2+t_1
 \end{split}  
\end{displaymath}
has determinant 1, and in the new coordinates, the bound 
for ${\cal I}(\mathbf t_{\hat J})$ now implies the cruder inequality
$$ {\cal I}(\mathbf t_{\hat J}) \ll 
\Big(\prod_{j=1}^7 (1+|x_j|)\Big)^{\eta-99/100}
\big( (1+|x_1+x_4|)(1+|x_2-x_3+x_6+x_7|)(1+|x_5-x_7|)\big)^{\eta-1/2}.
$$
The special case $n=7$, $m=3$ of  Lemma \ref{lem13} now shows that
when $\eta= 1/200$ the function ${\cal I}(\mathbf t_{\hat J})$ is indeed 
integrable.\\

%\medskip

We now consider the remaining terms in \rf{harder2}, and begin with a summand
corresponding to a set $J$ with $\{3,6\}\subset J$, but such that none of the conditions (II), (IV), (V) is met. Then, the integrand in \rf{harder2} still has a simple
pole at $z_8=0$ with residue given by \rf{coeff1} where now 
$\hat J = J\cup \{8\}$. The argument following  \rf{coeff1} remains valid
if one puts $z_j=0$ for $j\in J$, and the error crucial to our present discussion
is still given by \rf{star}. This integral being a lower-dimensional version
of the original \rf{star}, it is apparent from the above analysis that 
the present \rf{star} is again bounded by $O(P^3\Delta^{1/200})$, as required.
Next, consider a set $J$ that occurs in \rf{harder2} with $\{3,6\}\subset J$,
and such that at least one further condition among (II), (IV) and (V)
is satisfied. Then, the integrand in \rf{harder2} has a pole of order 2 or 3
at $u_8=0$, and its residue is a certain linear combination of
\be{coeff2}
 P^3 \int_{(\tilde{\textbf{c}}^{\hat J})}  \frac{\partial^l}{\partial z_8^l} 
\, \frac{{\Theta}_{\Delta}(A\textbf{z}_{ J}) \,  
}{
  \langle \textbf{z}_J \rangle \lambda_1(\textbf{z}_J) \prod_{\nu=2}^5 (z_8 - \lambda_{\nu}(\textbf{z}_J))  }\Big|_{z_8=0}  d\textbf{z}_{\hat J}
\ee
with coefficients that are rational polynomials in $\log P$. In order to complete
the current programme, we again have to replace $\Delta$ by 0 in \rf{coeff2}
and estimate the error. Following the previous appoach, this error is
\be{star2}
  P^3 \int_{(\tilde{\textbf{c}}^{\hat J})} 
 \frac{\partial^l}{\partial z_8^l}
\big(\Psi_\Delta(\textbf{z}_{J})
- \Psi_0(\textbf{z}_{J})\big)\Big|_{z_8=0} \, d\mathbf z_{\hat J}
\ee
Since the relevant bounds in Lemma \ref{lem12} hold for partial derivatives as well, we can estimate the integrand here to
the same precision as in \rf{star1}, and it then transpires that again 
the relevant error is bounded by $O(P^{3+\varepsilon} \Delta^{1/200})$.\\

%\medskip

This completes the discussion of all terms in \rf{harder2} where (III) holds.
Next, we consider the case $J=\{6,7\}$ and put $\hat J =J\cup \{8\}$. Then again,
the integrand in \rf{harder2} has a simple pole at $z_8=0$ with residue
given by \rf{coeff1}. The dependence on $\Delta$ is then removed by the 
argument following \rf{coeff1}, but now with the function
\begin{displaymath}
\begin{split} \Psi_\Delta(\mathbf z_{\hat J}) =
 & \textstyle \widehat{f}_{\text{id}, \Delta, \delta}\left(\frac{1}{3} + z_9, \frac{1}{3} + z_{10}, \frac{1}{3}  - z_1 -z_2 - z_4-z_5 -   z_9 - z_{10}\right)
\\  & \textstyle\widehat{k}^{\epsilon_3}_{\Delta}\left(\frac{1}{3} + z_1 + z_4 +   z_5 +  z_9 + z_{10} , z_2\right)  \textstyle\widehat{k}_{\Delta}^{\epsilon_1}\left(\frac{1}{3} + z_1 + z_5 + z_{10} , \frac{1}{3} - z_1 - z_5 - z_9 - z_{10}\right) \\
   & \textstyle\widehat{k}^{\epsilon_2}_{\Delta}\left(z_1, \frac{2}{3}  - z_1 - z_{10}\right)   \widehat{q}(1 + z_3 - z_2 - z_1) H(\textbf{z}_{\hat J}). 
  \end{split}  
\end{displaymath}
that is obtained from \eqref{explicitres} with $z_6 = z_7 =z_8 = 0$. With this new meaning of $ \Psi_\Delta$, one has to estimate the integral \rf{star}.
Proceeding as before, one finds that \rf{star1} still holds with 
${\cal I}(\mathbf t_{\hat J})$ now redefined as
\begin{displaymath}
\begin{split}
   & \bigl((1 + |z_9|)(1+|z_{10}|)(1+|z_1+z_2+z_4+z_5 + z_9+z_{10}|)\bigr)^{\eta - 1} (1+|z_1|)^{\eta-1}(1+|z_1+z_{10}|)^{\eta - 1/2} \\
& (1+|z_1  + z_5 + z_{10}|)^{\eta - 1}(1 + |z_1 +z_5 +  z_9 + z_{10}|)^{\eta- 1/2}  (1+|z_2|)^{\eta-1}\\
&(1+|z_1+z_4+z_5 + z_9+z_{10}|)^{\eta-1/2} 
(1+|z_1+z_2-z_3|)^{-2}
 \prod_{j\in\{1, 2,  4,5,  9, 10\}} (1+|z_j|)^{1/1000}  \end{split}  
\end{displaymath}
We take $\eta=1/200$ and conclude that the integral in \rf{star} is $O(P^3\Delta^\eta)$ provided that the function in the previous display is integrable 
over $\Bbb R^7$. This follows from Lemma \ref{lem13}, because the linear
transformation
\begin{displaymath}
\begin{split}
&x_1 = t_9,\quad x_2 = t_{10}, \quad x_3 = t_1+t_2+t_4+t_5+t_9+t_{10}, \\
& x_4=t_1, \quad  x_5 = t_1+t_5 + t_{10},\quad x_6 = t_2 \quad x_7 = t_1+t_2-t_3 
\end{split}  
\end{displaymath}
has determinant 1  and shows  the above product bounded by
$$ (1+ |x_7|)^{-2}  \bigl((1+|x_1|) \cdots (1+|x_6|)\bigr)^{\eta - 99/100}\bigl((1+|x_2 + x_4|)(1+|x_1 + x_5|)(1+|x_3 - x_6|)\bigr)^{\eta-1/2}. $$ 
For sets $J$ with $\{6,7\}\subset J$ we can modify this argument in much the 
same way as in the case $\{3,6\}\subset J$, and obtain the same error estimate.\\

%\medskip
The next case that we consider is $J=\{4,7\}$. In this situation, one takes
$z_4 = z_7 =z_8 = 0$ in \eqref{explicitres}, and reconsiders the previous
error analysis with 
\begin{displaymath}
\begin{split} 
\Psi_\Delta(\mathbf z_{\hat J}) =
 &\textstyle \widehat{f}_{\text{id}, \Delta, \delta}\left(\frac{1}{3} + z_9, \frac{1}{3} + z_{10}, \frac{1}{3}  - z_1 -z_2 -z_5 -   z_9 - z_{10}\right)\\ &\textstyle
 \widehat{k}^{\epsilon_3}_{\Delta}\left(\frac{1}{3} + z_1 +     z_5- z_6 +  z_9 + z_{10} , z_2\right) \widehat{k}_{\Delta}^{\epsilon_1}\left(\frac{1}{3} + z_1 + z_5 + z_{10} , \frac{1}{3} - z_1 - z_5+z_6 - z_9 - z_{10}\right) \\
   &\textstyle  \widehat{k}^{\epsilon_2}_{\Delta}\left(z_1, \frac{2}{3}  - z_1 - z_{10}\right) \widehat{q}(1 + z_3 - z_2 - z_1)   H(\textbf{z}_{\hat J}). 
  \end{split}  
\end{displaymath}
As before, one verifies \rf{star1} with 
\begin{displaymath}
\begin{split}
{\cal I}(\mathbf t_{\hat J}) =& \bigl((1 + |z_9|)(1+|z_{10}|)(1+|z_1+z_2+z_5 + z_9+z_{10}|)\bigr)^{\eta - 1} (1+|z_1|)^{\eta-1}(1+|z_1+z_{10}|)^{\eta - 1/2} \\
& (1+|z_1  + z_5 + z_{10}|)^{\eta - 1/2}(1 + |z_1 +z_5-z_6 +  z_9 + z_{10}|)^{2\eta- 3/2}  (1+|z_2|)^{\eta-1}\\
& 
(1+|z_1+z_2-z_3|)^{-2}
 \prod_{j\in\{1, 2, 5,6,  9, 10\}} (1+|z_j|)^{1/1000} 
\end{split}  
\end{displaymath}
Here the encounter the new phenomenon that the argument $z_1 +z_5-z_6 +  z_9 + z_{10}$
occurs twice in the definition of $\Psi_\Delta$. The change of variable
\begin{displaymath}\begin{split}
&  x_1 = t_9, \quad x_2 = t_{10}, \quad x_3 = t_1 +t_2 +t_5 +   t_9 + t_{10}, \quad x_4 = t_2,
\\ &   x_5 = t_1,  \quad x_6 = t_1 + t_5 -t_6 +t_9+ t_{10}, \quad x_7 = t_1+t_2-t_3
\end{split}\end{displaymath}
has determinant 1 and shows the previous product bounded by
$$\Big(  \prod_{j=1}^5 (1+|x_j|)\Big)^{\eta - 99/100} (1+|x_6|)^{2\eta-5/4}(1+|x_7|)^{-2}
(1+|x_2+x_4|)^{\eta-1/2} (1+|x_3-x_6-x_1|)^{\eta-1/2}. $$
For $\eta=1/200$, this product is integrable, as one finds from Lemma \ref{lem13} with $n=5$,
and the relevant error is therefore bounded as before. Also, it transpires that the more general case where $\{4,7\}\subset J$ is covered by this argument and the discussion 
towards the end of case (III) above.\\
 
%\medskip

Finally, we turn to case (V) and examine the situation where $J=\{1,2,4,5\}$. Here
\eqref{explicitres} with $z_1 = z_2 = z_4 = z_5 =z_8 = 0$ reduces to
\begin{displaymath}
\begin{split}
\Psi_\Delta(\mathbf z_{\hat J}) = &\textstyle \widehat{f}_{\text{id}, \Delta, \delta}\left(\frac{1}{3} + z_9, \frac{1}{3} + z_{10}, \frac{1}{3}   -   z_9 - z_{10}\right)\underset{u=0}{\text{res}} \widehat{k}^{\epsilon_3}_{\Delta}\left(\frac{1}{3}  - z_6 - z_7+  z_9 + z_{10} , z\right) \\
   &\textstyle \widehat{k}_{\Delta}^{\epsilon_1}\left(\frac{1}{3} + z_{10} - z_7 , \frac{1}{3} +z_6 + z_7 - z_9 - z_{10}\right)  \underset{u=0}{\text{res}}\widehat{k}^{\epsilon_2}_{\Delta}\left(z, \frac{2}{3}  + z_7 - z_{10}\right)   \widehat{q}(1+z_3) H(\textbf{z}_{\hat J}). 
  \end{split}  
\end{displaymath}
Following the previous argument, enhanced by Lemma \ref{lem11}(iii), once more, one confirms \rf{star1} with
\begin{displaymath}
\begin{split}
{\cal I}(\mathbf t_{\hat J}) =& \bigl((1 + |z_9|)(1+|z_{10}|)(1+|z_9+z_{10}|)\bigr)^{\eta - 1} (1+|z_6+z_7-z_9-z_{10}|)^{2\eta-3/2}\\
& (1+|z_7-z_{10}|)^{\eta-2} (1+|z_3|)^{-2} 
\prod_{j\in\{3,6,7, 9, 10\}} (1+|z_j|)^{1/1000}. 
\end{split}  
\end{displaymath}
This product is integrable, as one finds using the substitution
$$ x_1=t_9,\quad x_2=t_{10},\quad x_3=t_6+t_7-t_9-t_{10}, \quad x_4= t_7-t_{10},\quad x_5 = t_3. $$
Thus the error analysis can be performed as before, and it is again immediate that the
more general case $\{1,2,4,5\} \subset J$ is covered by this approach.

This completes  the discussion of the term \eqref{harder2}. % and we have now proved that
%the sum \rf{harder2} equals $P^3R(\log P) + O(P^{3+\varepsilon} \Delta^{1/200} +P^{3-\eta_1}\Delta^{-99})$ in which $R$ is a certain polynomial. 
An inspection of the terms \eqref{harder4} and \eqref{harder5} shows that a very similar
treatment is possible, and that one does not encounter integrals that need to be checked
for existence, other than those examined above. The sums \rf{harder1} and \rf{harder3}
contain a derivative in the integrand, and one can handle this in a fashion identical to
the treatment of the derivative in \rf{star2}. Then again, one may appeal to the discussion
above to conclude, as desired, that the dependence of $\Delta$ can be removed with acceptable error. One then arrives at the formula
$$
%\be{v1fertig}
P^3Q_{\Delta}^{(1)}(\log P) = P^3Q_0^{(1)} (\log P) + O(P^{3+\varepsilon} \Delta^{1/200})
%\ee
$$
where $Q_0^{(1)}$ is a certain polynomial.

\subsection{The analysis of $V^{(2)}(P )$ and $V^{(3)}(P )$} The discussion of $V^{(2)}(P )$ is similar to the work in the previous section, but rather less complex. We begin with an inspection of the transition from \eqref{resV2} to \eqref{prelim1}. For $J \subset \{1, \ldots, 6\}$ and $\textbf{z} \in \Bbb{C}^8$, let $\textbf{z}_J \in \Bbb{C}^8$ be the vector that is obtained from $\textbf{z}$ on replacing $z_j$ by 0 for all $j \in J$. Note that this is in accord with a similar convention in Section \ref{mult}. Also, for $J \subset \{1, \ldots, 5\}$ define the vector $\textbf{c}^J \in \Bbb{R}^{8-|J|}$ by $c_j = -2\eta_1$ for $1 \leq j \leq 5$, $j \not\in J$, and by $c_6 = 3\eta_1$, $c_7 = c_8 = \eta_1$, and put
\begin{displaymath}
  \langle \textbf{z}_j \rangle = \prod_{\substack{1 \leq j \leq 5\\ j \not\in J}} z_j, \quad L(\textbf{z}) = (z_6 - z_5 - z_4)(z_6 - z_4 - z_2)(z_6 - z_5 - z_3). 
\end{displaymath}
Then in analogy with \eqref{simple}, one finds from \eqref{resV2} that
\begin{equation}\label{star3}
  V^{(2)}_{\Delta}(P ) = \sum_{J \subset \{1, \ldots, 5\}} \frac{48}{(2\pi i)^{8 - |J|}} \int_{(\textbf{c}^J)} \frac{P^{3+z_6}\Theta_{\Delta}(A\textbf{z}_J)}{\langle \textbf{z}_J\rangle L(\textbf{z}_J)} d\textbf{z}_J. 
\end{equation}
To derive \eqref{prelim1}, one now shifts the line $\Re z_6 = 3\eta_1$ to $\Re z_6 = -\eta_1$. The integral over $\Re z_6 = -\eta_1$ contributes $O(P^{3-\eta_1}\Delta^{-99})$ by straightforward estimates. The only pole that may occur in the shift is at $z_6 = 0$, and this will be the case if and only if $z_6$ devides the polynomial $L(\textbf{z}_J)$, that is, if one of the conditions
\begin{displaymath}
  {\rm(I)} \,\,\, \{4, 5\} \subset J, \quad\quad  {\rm (II)} \,\,\, \{2, 4\} \subset J, \quad \quad {\rm(III)} \,\,\, \{3, 5\} \subset J
\end{displaymath}
holds. Summands in \eqref{star3} with this property generate residues at $z_6 = 0$ that assemble to $P^3 Q^{(2)}_{\Delta}(\log P)$. Following the line of attack in the previous section, we replace $\Delta$ by 0 in these residues and estimate the resulting error.

First suppose that $J= \{4, 5\}$ and put $\hat{J} = J \cup \{6\}$. Then the integrand in \eqref{star3} has a simple pole at $z_6 = 0$ with residue
\begin{equation}\label{star4}
  P^3 \int_{(\tilde{\textbf{c}}^{\hat{J}})} \frac{\Theta_{\Delta}(A\textbf{z}_{\hat{J}})}{z_1z_2^2z_3^2 } d\textbf{z}_{\hat{J}};
\end{equation}
here $(\tilde{\textbf{c}}^{\hat{J}})$ is the line $\Re z_j = -2\eta_1$ for $1 \leq j \leq 5$, $j \not\in J$, and $\Re z_7 = \Re z_8 = \eta_1$. The integrand in \eqref{star4} that we now abbreviate with $\Psi_{\Delta}(\textbf{z}_{\hat{J}})$ admits the alternative expression
\begin{equation}\label{star5}
\begin{split}
  \Psi_{\Delta}(\textbf{z}_{\hat{J}}) = &\widehat{f}_{\text{id}, \Delta, \delta}\left(\textstyle\frac{1}{3}+z_8 , \frac{1}{3} + z_7, \frac{1}{3}+ z_1 - z_7 - z_8 \right) \widehat{k}^{\pm}_{\Delta}\left(\textstyle\frac{1}{3} + z_3 - z_7 - z_8, \frac{1}{3}- z_3 + z_8\right) \\
& \widehat{f}_{\Delta}\left(\textstyle\frac{2}{3} - z_8\right)\widehat{f}_{\Delta}\left(\textstyle\frac{2}{3}- z_1 + z_7 + z_8\right)  \widehat{q}(z_1 - z_2 - z_3 )  H(\textbf{z}_{\hat{J}})
\end{split}  
\end{equation}
where $H$ is holomorphic in $|\Re u_j| \leq 1/4$ and satisfies the bound \eqref{boundH}. This much is obtained in analogy to the argument leading to \eqref{Psi1}. By Lemma   \ref{lem10} (i) and Lemma \ref{lem12}, on the lines of integration parametrized by $z_j = c_j + it_j$, one has 
\begin{displaymath}
  \Psi_{\Delta}(\textbf{z}_{\hat{J}}) - \Psi_0(\textbf{z}_{\hat{J}}) \ll \Delta^{\eta} \mathcal{I}(\textbf{t}_{\hat{J}})
\end{displaymath}
 where $\eta = 1/200$ and 
 \begin{displaymath}
 \begin{split}
 \mathcal{I}(\textbf{t}_{\hat{J}})=&\bigl(  (1+|z_8|)(1+|z_7|)(1+|z_1 - z_7 - z_8|)\bigr)^{\eta - 1}(1 + |z_3 - z_7 - z_8|)^{\eta-1}(1+|z_3 - z_8|)^{\eta - 1/2}\\
 &(1+|z_8|)^{\eta - 1}(1+ |z_1- z_7 - z_8|)^{\eta - 1}(1+|z_1-z_2-z_3|)^{-2} \prod_{j \in \{1, 2, 7, 8\}} ( 1+ |z_j|)^{1/1000}. 
\end{split}
\end{displaymath}
The change of variables
\begin{displaymath}
  x_1 = t_8, \quad x_2 = t_7, \quad x_3 = t_1 - t_7 - t_8, \quad  x_4 = t_1 - t_2 - t_3, \quad x_5 = t_3 - t_7 - t_8, 
\end{displaymath}
has determinant 1, and in the new coordinates one can bound $\mathcal{I}(\textbf{t}_{\hat{J}})$ as
\begin{displaymath}
\begin{split}
 \bigl( (1+|x_1|) (1+|x_3|)\bigr)^{2\eta-199/100} (1+|x_4|)^{-2}  \bigl( (1+|x_2|) (1+|x_5|)\bigr)^{\eta-99/100}(1+|x_5+x_2|)^{\eta-1/2}. % (1+|x_1|)^{\eta_2-1} (1+|x_3|)^{\eta_2-1}.
%\bigl( (1+|x_1|)(1+|x_2|)(1+|x_3|)(1+|x_5|)\bigr)^{\eta_2-99/100} (1+|x_4|)^{-2} (1+|x_5+x_2|)^{\eta_2-1/2} (1+|x_1|)^{\eta_2-1} (1+|x_3|)^{\eta_2-1}.
\end{split}
\end{displaymath}
Hence by Lemma \ref{lem13} one finds that $\mathcal{I}(\textbf{t}_{\hat{J}})$ is integrable, and consequently that 
\begin{displaymath}
  \int_{(\tilde{\textbf{c}}^{\hat{J}})} \bigl( \Psi_{\Delta}(\textbf{z}_{\hat{J}}) - \Psi_0(\textbf{z}_{\hat{J}}) \bigr) d\textbf{z}_{\hat{J}}\ll \Delta^{\eta}.
  \end{displaymath}
Therefore one may indeed replace $\Delta$ by 0 in \eqref{star4} at the cost of an error not exceeding $O(\Delta^{\eta} P^{3})$. 

Little change is necessary for the treatment of summands corresponding to other sets $J$ with $\{4, 5\} \subset J \subset \{1, \ldots, 5\}$. For some of these sets, poles of order 2 or 3 at $z_6 = 0$ occur, and it will then be necessary to consider certain partial derivatives with respect to $z_6$, similar to the occurrence of derivatives in \eqref{star2}. An inspection of the deliberations following \eqref{star2} shows that the current situation is fully covered by the above treatment, and it transpires that for all sets $J$ with $\{4, 5\} \subset J \subset \{1, \ldots, 5\}$, the computation of the residue at $z_6 = 0$ of the integrand in \eqref{star3} may be performed with $\Delta = 0$ at the cost of a total error not exceeding $O(P^{3+\varepsilon} \Delta^{\eta})$.\\ 

Next, consider the case $J = \{2, 4\}$ and follow the same pattern as before. The right hand side of \eqref{star5} is now to be replaced by 
\begin{displaymath}
\begin{split}
  &\widehat{f}_{\text{id}, \Delta, \delta}\left(\textstyle\frac{1}{3}+z_8 , \frac{1}{3} + z_7, \frac{1}{3}+ z_1 - z_7 - z_8 \right) \widehat{k}^{\pm}_{\Delta}\left(\textstyle\frac{1}{3} + z_3 - z_7 - z_8, \frac{1}{3}- z_3 - z_5+ z_8\right) \\
& \widehat{f}_{\Delta}\left(\textstyle\frac{2}{3} - z_8\right)\widehat{f}_{\Delta}\left(\textstyle\frac{2}{3}- z_1  + z_5 + z_7 + z_8\right)  \widehat{q}(z_1 - z_3 - z_5)H(\textbf{z}),
\end{split}
\end{displaymath}
and one is then led to check integrability of the product
\begin{displaymath}
\begin{split}
&\bigl(  (1+|z_8|)(1+|z_7|)(1+|z_1 - z_7 - z_8)\bigr)^{\eta - 1}(1 + |z_3 - z_7 - z_8|)^{\eta-1}(1+|z_3+z_5 - z_8|)^{\eta - 1/2}\\ &(1+|z_8|)^{\eta - 1}(1+ |z_1 - z_5- z_7 - z_8|)^{\eta - 1} (z_1 - z_3 - z_5)^{-2} \prod_{j \in \{1, 3, 5, 7, 8\}} ( 1+ |z_j|)^{1/1000}
\end{split}
\end{displaymath}
which is provided by Lemma \ref{lem13} after the substitution
\begin{displaymath}
   x_1 = t_8, \quad x_2 = t_7, \quad x_3 = t_1 - t_7 - t_8, \quad x_4 = t_3  - t_7 - t_8, \quad x_5 = t_1 - t_3 - t_5.
\end{displaymath}

For the other cases with $\{2, 4\} \subset J \subset \{1, \ldots, 5\}$ this argument may be modified along the lines suggested in the penultimate paragraph.\\

Now consider $J = \{3, 5\}$ where the right hand side of \eqref{star5} should read
\begin{displaymath}
\begin{split}
  &\widehat{f}_{\text{id}, \Delta, \delta}\left(\textstyle\frac{1}{3}+z_8 , \frac{1}{3} + z_7, \frac{1}{3}+ z_1 - z_7 - z_8 \right) \widehat{k}^{\pm}_{\Delta}\left(\textstyle\frac{1}{3}  -z_7 - z_8, \frac{1}{3}- z_4+ z_8\right) \\
& \widehat{f}_{\Delta}\left(\textstyle\frac{2}{3} +z_4- z_8\right)\widehat{f}_{\Delta}\left(\textstyle\frac{2}{3}- z_1  +  z_7 + z_8\right)   \widehat{q}(z_1 - z_2 - z_4)H(\textbf{z}), 
\end{split}
\end{displaymath}
and one has to consider the product
\begin{displaymath}
\begin{split}
&\bigl(  (1+|z_8|)(1+|z_7|)(1+|z_1 - z_7 - z_8)\bigr)^{\eta - 1}(1 + |  z_7 + z_8|)^{\eta-1}(1+|z_4 - z_8|)^{\eta - 1/2}\\
&(1+|z_4 - z_8|)^{\eta - 1}(1+ |z_1- z_7 - z_8|)^{\eta - 1} (1 + |z_1 - z_2 - z_4|)^{-2} \prod_{j \in \{1, 3, 5, 7, 8\}} ( 1+ |z_j|)^{1/1000}. 
\end{split}
\end{displaymath}
The substitution
\begin{displaymath}
   x_1 = t_8, \quad x_2 = t_7, \quad x_3 = t_1 - t_7 - t_8, \quad x_4 = t_4   - t_8, \quad x_5 = t_1 - t_2 - t_4
\end{displaymath}
shows that Lemma \ref{lem13} again yields the desired integrability of $\mathcal{I}(\textbf{t}_{\hat{J}})$. As before, other cases with $\{3, 5\} \subset J \subset \{1, \ldots, 5\}$ are very similar, and we may now conclude that there exists a polynomial $Q_0^{(2)}$ of degree at most 2 and such that
\begin{displaymath}
 P^3Q_{\Delta}^{(2)} (P )= P^3 Q_0^{(2)}(\log P ) + O(P^{3+\varepsilon}\Delta^{1/200}  ).
\end{displaymath}

Finally, we turn our attention to $V_{\Delta}^{(3)}(P )$. Here a brief inspection of \eqref{prelim2} suffices to confirm that the now familiar procedure again yields a polynomial of degree at most 1 with 
\begin{displaymath}
 P^3Q_{\Delta}^{(3)}(\log P) = P^3 Q_0^{(3)}(\log P ) + O(P^{3+\varepsilon}\Delta^{1/200} ).
\end{displaymath}
We may leave the details to the reader. This completes the proof of \eqref{claim}.

  \section{The Peyre constant}\label{peyre}

In remains to determine the leading coefficient of the polynomial $Q_0$ in \eqref{asymp}  in order to complete the proof of \eqref{asympV} and hence of Theorem \ref{thm1}. %In the decomposition \eqref{decompQ} 
It suffices to consider $Q_0^{(1)}$, since $Q_0^{(2)}$ and $Q_0^{(3)}$ have smaller degree. An inspection of the terms \eqref{harder1} -- \eqref{harder5} shows that only the terms \eqref{harder1} and \eqref{harder4} contribute, and only if $J = \{1, \ldots, 7\}$. Thus    the leading coefficient is given by  
\begin{displaymath}
\begin{split}
& \sum_{{\bm \epsilon} \in E} \frac{48}{(2\pi i)^2}  \int_{(\eta)}  \int_{(\eta)}   \frac{2 - \frac{1}{2}}{24}\Theta_{0}\bigl(A (0, \ldots, 0, z_9, z_{10})^T\bigr) dz_9 \, dz_{10}\\
= &  \sum_{{\bm \epsilon} \in E} \frac{3 }{(2\pi i)^2}  \int_{(\eta)}  \int_{(\eta)}   %\\
%& \quad\quad\quad\quad
 \Theta_0(z_9, z_{10}, -z_9 -z_{10}; z_{10}, -z_9 - z_{10}; 0, -z_{10}; z_9 + z_{10}, 0; 0)   dz_9 \, dz_{10}.
\end{split}
\end{displaymath}
 By \eqref{defphi} --  \eqref{x}  and Lemma  \ref{lem11} (iii), this equals
\begin{displaymath}
\begin{split}
 &  \sum_{\epsilon \in \{\pm\}} \frac{6\, \Xi_G(1, \ldots, 1) }{(2\pi i)^2}  \int_{(\eta)}  \int_{(\eta)}   \frac{ \widehat{f}_{\pi, 0, \delta}\left(\frac{1}{3} + z_9, \frac{1}{3} + z_{10}, \frac{1}{3} - z_9 - z_{10}\right)\widehat{k}_0^{\epsilon}(\frac{1}{3} + z_{10}, \frac{1}{3} - z_9 - z_{10})}{(\frac{2}{3}-z_{10})(\frac{2}{3} + z_9+z_{10})}     dz_9 \, dz_{10}.
 \end{split}
 \end{displaymath}
We call this constant $C_{\delta}$ and claim that $C_{\delta} = C_0$ for all sufficiently small $\delta$. This can be checked by direct computation, but we can also argue as follows. It is clear from \eqref{deffpi0delta} and \eqref{defh} that $C_{\delta} \rightarrow C_0$ as $\delta \rightarrow 0$. We have already shown that 
\begin{displaymath}
   \frac{V(P )}{P^3 (\log P)^4} = C_{\delta} + O_{\delta}\left(\frac{1}{\log P}\right)
\end{displaymath}
for any $\delta > 0$. 
Combining these two asymptotic relations we find that
\begin{displaymath}
  \lim_{P \rightarrow \infty} \frac{V(P )}{P^3 (\log P)^4} = C_{0}
\end{displaymath}
 as claimed. By Lemma \ref{lem11}(ii) it therefore remains to compute
\begin{displaymath} 
C_0 =       \frac{6\, \Xi_G(1, \ldots, 1)}{(2\pi i)^2}  \int_{(\frac{1}{3} + \eta)}  \int_{(\frac{1}{3} + \eta)}   \frac{  \widehat{k}_0^+(  z_{10}, 1 - z_9 - z_{10}) + \widehat{k}_0^-(  z_{10}, 1 - z_9 - z_{10})}{ z_9   (z_9 + z_{10})^2  (1-z_{10}) }     dz_9 \, dz_{10}.
 \end{displaymath}
We make a change of variables  $v_1 = 1 - z_{10}$, $v_2 = z_9 + z_{10}$   and insert the definition of $\widehat{k}^{\pm}_0$ as a double Mellin transform. In this way we see 
\begin{displaymath}
C_0 = 6  \, \Xi_G(1, \ldots, 1) \int_{\mathcal{Q}(2)} \bigl(k^{+}_0(\textbf{x}) + k^{-}_0(\textbf{x})\bigr) \frac{1}{(2\pi i)^2}  \int_{(\frac{2}{3})}\int_{(\frac{2}{3})}   \frac{ x_1^{-v_1} x_2^{-v_2}  d\textbf{v}  }{(v_1+v_2 - 1)v_1v_2^2}    d\textbf{x}.
\end{displaymath}
Shifting the $v_1, v_2$-contours to the right, we see that the inner integral vanishes unless $x_1, x_2 \leq 1$. In this region, however,  $k^-_0$ is constantly 1. In particular we see that  the non-canonical choices in the definitions \eqref{defkpm},   \eqref{defk} of $k^{-}_0$ and of $q$ play no role, as it should be. Shifting the $v_1, v_2$-contours to the left, one readily computes the inner integral. The constant $\Xi_G(1, \ldots, 1)$ has been computed in \eqref{A7}, and  the rest is a  straightforward evaluation of elementary integrals: 
 \begin{displaymath}
 \begin{split}
&6  \, \Xi_G(1, \ldots, 1) \int_0^1\int_0^1 \bigl(k^{+}_0(\textbf{x}) + k^{-}_0(\textbf{x})\bigr)  \left(\log (x_2) - 1 + \frac{1 + \max(0, \log(x_1/x_2))}{\max(x_1, x_2)}  \right)  d\textbf{x}\\
=& 6  \, \Xi_G(1, \ldots, 1) \left(\Bigl(- \frac{5}{4} + \frac{\pi^2}{12} + 2\log 2 \Bigr) + 1\right) = \frac{1}{2}(\pi^2 + 24 \log 2 - 3)  \prod_p \Big(1-\f{9}{p^2}+\f{16}{p^3}-\f{9}{p^4}+\f{1}{p^6}\Big)
\end{split}
\end{displaymath}
as claimed in \eqref{lead}.

\end{document}